\documentclass[12pt,amstex]{article}

\usepackage{amsmath}
\usepackage{amsfonts}

\newtheorem{theo}{Theorem}

\newtheorem{lemm}{Lemma}
\newtheorem{defn}{Definition}

\newcommand{\lbl}{\label}

\newcommand{\eq}[1]{$(\ref{#1})$}

\def\N{\mathbb{N}}
\newcommand{\nat}{\mathbb{N}}
\def\bx{{\bf x}}
\def\bX{{\bf X}}
\def\by{{\bf y}}
\def\bw{{\bf w}}

\def\bv{{\bf v}}
\def\bY{{\bf Y}}
\def\bW{{\bf W}}
\def\bz{{\bf z}}

\def\E{\mathbb{E}}
\def\0{{\bf 0}}

\def\R{\mathbb{R}}

\def\ch{{F}}

\renewcommand{\E}{\mathbb E \,}

\newcommand{\tod}{\stackrel{{\cal D}}{\longrightarrow}}
\newcommand{\toD}{\stackrel{{\cal D}}{\longrightarrow}}
\newcommand{\eqd}{\stackrel{{\cal D}}{=}}
\newcommand{\toP}{\stackrel{{P}}{\longrightarrow}}
\newcommand{\supp}{{\rm supp}}
\newcommand{\ska}{{\rm supp}(\kappa)}
\newcommand{\skala}{\Gamma_\la}

\newcommand{\inLL}{\stackrel{{L^2}}{\longrightarrow}}
\newcommand{\inL}{\stackrel{{L^1}}{\longrightarrow}}

\newcommand{\inLq}{\stackrel{{L^q}}{\longrightarrow}}

\newcommand{\eqco}{\setcounter{equation}{0}}
\newcommand{\thco}{\setcounter{theo}{0}}
\newcommand{\prco}{\setcounter{prop}{0}}
\newcommand{\laco}{\setcounter{lemm}{0}}
\newcommand{\coco}{\setcounter{coro}{0}}
\newcommand{\cjco}{\setcounter{conj}{0}}

\newcommand{\deco}{\setcounter{defn}{0}}
\newcommand{\allco}{\eqco  \thco \prco \laco \coco \cjco \deco}
\newcommand{\qed}{\hfill{\rule[-.2mm]{3mm}{3mm}}}

\newcommand{\Po}{{\cal P}}
\newcommand{\kamax}{{\|\ka\|_\infty}}
\newcommand{\MM}{{\cal M}}
\newcommand{\tPo}{\tilde{{\cal P}}}
\newcommand{\tH}{\tilde{{\cal H}}}
\newcommand{\tg}{\tilde{g}}
\newcommand{\tB}{\tilde{B}}

\newcommand{\X}{{\cal X}}
\newcommand{\LL}{{\cal L}}
\newcommand{\Y}{{\cal Y}}

\renewcommand{\H}{{\cal H}}
\renewcommand{\P}{{\cal P}}
\renewcommand{\SS}{{\cal S}}
\newcommand{\A}{{\cal A}}

\newcommand{\Var}{{\rm Var}}

\newcommand{\const}{{\rm const.} \times}
\newcommand{\x}{{\bf x}}

\newcommand{\limn}{\lim_{n \to \infty} }

\newcommand{\YY}{{\cal Y}}

\newcommand{\F}{{\cal F}}

\newcommand{\NN}{{\cal N}}

\newcommand{\eps}{\varepsilon}
\def\bdm{\begin{displaymath}}
\newcommand{\edm}{\end{displaymath}}
\def\benu{\begin{enumerate}}
\def\eenu{\end{enumerate}}
\def\beqn{\begin{equation}}
\def\eeqn{\end{equation}}
\def\be{\begin{equation}}
\def\ee{\end{equation}}
\def\bea{\begin{eqnarray}}
\def\eea{\end{eqnarray}}
\newcommand{\bean}{\begin{eqnarray*}}
\newcommand{\eean}{\end{eqnarray*}}
\newcommand{\bear}{\begin{eqnarray}}
\newcommand{\eear}{\end{eqnarray}}

\newcommand{\1}{{\bf 1}}

\newcommand{\y}{{\bf y}}

\renewcommand{\epsilon}{\varepsilon}

\def\inter{{\rm int}}
\def\dist{{D}}

\def\R{\mathbb{R}}

\def\A{{\cal A}}

\def\qed{\hfill\hbox{${\vcenter{\vbox{
    \hrule height 0.4pt\hbox{\vrule width 0.4pt height 6pt
    \kern5pt\vrule width 0.4pt}\hrule height 0.4pt}}}$}}

\def\la{{\lambda}}
\def\ka{{\kappa}}

\begin{document}

\title{\bf Convergence of random measures in geometric probability}

\author{Mathew D. Penrose \\
\\
{\normalsize{\em University of Bath 
}} }
\date{August 2005}
\maketitle

%\footnotetext{$~^1$ Department of Mathematical Sciences, University
%of Bath, Bath BA2 7AY, United Kingdom: {\texttt m.d.penrose@bath.ac.uk} }



\begin{abstract}
Given $n$ independent random marked $d$-vectors $X_i$ with a common density, 
define the measure $\nu_n = \sum_i \xi_i $, where
 $\xi_i$ is a measure 
(not necessarily a point measure)
determined by the (suitably rescaled)  
set of points near $X_i$. Technically, this means here 
 that $\xi_i$  stabilizes with a suitable power-law
decay of the tail of the radius of stabilization.
For bounded test functions $f$ on $R^d$, 
we give a law of large numbers and central limit theorem
for $\nu_n(f)$. The latter implies weak convergence of
 $\nu_n(\cdot)$, suitably scaled and centred,
to a Gaussian field acting on bounded
test functions.
% (not only continuous test functions). 
%The measure $\mu_n$ does not need to be a point measure.
% (this was previously known for continuous test functions). 
%The result
%is given in general terms to allow for marked point processes. 
The general result is illustrated with 
applications  including the volume and surface measure of germ-grain
models with unbounded grain sizes.

\end{abstract}

\newpage

\section{Introduction }
%and main results}
\allco

This paper is concerned with the study of the
limiting behaviour of {\em random measures}
based on marked Poisson or binomial point processes in $d$-dimensional space,
arising as the sum of contributions from each point of
the point process.
% with the property that the
%contribution from a point
%is determined by the local configuration of nearby points.
Many random spatial measures can be described in these terms, 
%and there is a large literature on 
and general limit theorems, including
 laws of large numbers, central limit theorems, and large deviation
principles,
 are known for the total measure of such measures,
based on a notion of {\em stabilization} (local dependence); see
 \cite{PY1,PY2,PY4,SY}.

Recently, attention has turned to the asymptotic behaviour of
the measure itself (rather than only its total measure),
 notably in \cite{BY1,BY,HM,Pe04,PY5,SY}.
It is of interest to determine when one can
show weak convergence of this
 %finite-dimensional distributions of  this
 measure to a Gaussian random field.
As in Heinrich and Molchanov \cite{HM}, Penrose \cite{Pe04}, one can consider
a limiting regime where  a homogeneous Poisson
process is sampled over an expanding window.
%and on what class of test functions, under 
In an alternative limiting regime, the intensity  
of the point process becomes large 
and
% locally scale 
the point process
 is locally  scaled 
to keep the average density of points bounded; the latter approach
% type of limiting regime 
allows for point processes with non-constant densities
and is the one adopted here.

%We consider the case where 
A random measure 
is said to be {\em exponentially stabilizing}
when the contribution of an inserted point
is determined by the configuration
of (marked) Poisson points within a finite (though in general random) distance,
known as a radius of stabilization,
having a uniformly exponentially
decaying tail after scaling of space.
This concept was introduced by
 Baryshnikov and Yukich \cite{BY},
who proved  general results on  
weak convergence to a limiting Gaussian field 
for exponentially stabilizing measures.
A variety of random measures are exponentially stabilizing,
 including those concerned with
nearest neighbour graph, Voronoi and Delaunay graph,
germ-grain models with bounded grains, and packing; see \cite{BY}.

In the present work we 
%refine the second moment calculations
extend the results 
of \cite{BY}
% broadening 
% which enables us to
%broaden the class of random measures, underlying density functions
%and test functions  considered in \cite{BY} 
in several directions.
Specifically, in \cite{BY} attention is restricted 
%we relax the  restriction in \cite{BY}
%Also, unlike in \cite{BY} we do not 
%restrict attention 
to the case where 
the random measure
is concentrated at the points of the underlying 
%from each point of the
 point process, and 
%the restriction in \cite{BY}  
%is  atomic at that point.
to continuous test functions;  we relax both of these restrictions,
and so 
%the latter  
%to point measures
% we relax certain continuity
% assumptions made in \cite{BY}  on test functions so as 
are able to include
indicator functions of Borel sets as test functions. 
Moreover, we relax the condition of
exponential stabilization to {\em power-law stabilization}.
%give results 

% in fact, 
%our results require only a weaker {\em power-law stabilization}
%condition.  

Along with the central limit theorems, we give a law
of large numbers  for random measures (under weaker stabilization
conditions) which is 
proved using similar ideas.
  In Section \ref{secexamples}, we illustrate our general results
with applications
 to germ-grain models with unbounded grains,
 and to questions
of Voronoi coverage raised by Khmaladze and Toronjadze \cite{KT}.
Many other fields of application have 
been discussed elsewhere \cite{BY1,BY,PY1,PY4} and we
do not attempt to review all of these.

We state our general results in Section \ref{secresults}.
Our approach to proof may be summarized as follows. 
In the case where the underlying point process is Poisson, we 
obtain the  covariance structure of our
limiting random field
(and also to prove the law of large numbers)    
by
 a refinement of the approach to second moment calculations 
used in \cite{BY} and \cite{PY4}, often known as the
  {\em objective method}  \cite{AS,Stebk}, which is discussed in Section
\ref{secobjective}. 
To show that the limiting random field is Gaussian, we
borrow normal approximation results from \cite{PY5} which
were proved there using Stein's method (in contrast, \cite{BY}
uses the method of moments). Finally, to de-Poissonize
the central limit theorems (i.e., to extend them to binomial
point processes with a non-random number of points), in
Section \ref{secdepo} we
perform further second moment calculations using a version of
the objective method. This approach entails an annoyingly large
number of similar calculations 
%  is rather tedious to implement
(see Lemmas \ref{lemweak3} and \ref{fordepolem}) but  avoids
 the necessity of introducing  a notion of `external
 stabilization' (see Section
\ref{secresults}) which was used to deal with the second moment
calculations for de-Poissonization in \cite{BY}.
This, in turn, seems to be necessary to include germ-grain models
with unbounded grains.
We give our proofs in the general 
setting of {\em marked} point process, which is the context
for many of the applications.

%in a different manner
%from that used for the equivalent calculations in \cite{BY}

% prove the actual convergence to a 

% using the method of moments.

%In the case where the measure in question 
%is a sum of atomic contributions at each point of the underlying point process,
% In contrast with the results
%in \cite{HM,Pe04}, they did not require homogeneity of
%the underlying point processes.
%Later, Penrose and Yukich 
% \cite{PY5} showed how to derive similar results 
%via dependency graphs and Stein's method.
%However, the latter approach still requires
%second moment calculations.
%of \cite{BY}.

%Combining this second moment calculation with
%the arguments presented in \cite{PY5}, we
%are able to  strengthen the results
%given in \cite{BY} on convergence of measures.

\section{Notation and results}
\lbl{secresults}
\allco
Let $(\MM,\F_\MM,\mu_\MM)$ be a probability space (the {\em mark space}). 
Let $d \in \N$.
Let $\xi(\bx,t;\X,A)$ be a Borel-measurable $\R$-valued function  defined for
all 4-tuples $(\bx,t,\X,A)$, where $\X \subset \R^d\times \MM$
 is finite and where $(\bx,t) \in \X$ (so $\bx \in \R^d$ and $t \in \MM$),
and $A$ is a Borel set in $\R^d$. 
We assume that $\xi(\bx,t;\X):= \xi(\bx,t;\X,\cdot)$ is a
% signed
 measure on $\R^d$ with finite total measure.
(Our results actually hold when  
 $\xi(\bx,t;\X)$
%:= \xi(\bx,t;\X,\cdot)$
 is a
{\em signed} measure with finite total variation; see the remarks at 
the end of this section.) 
%Adapting the proofs to this more general case
%takes a little extra work and we omit this to save space.)
% variation (i.e., the sum of the total measure
%of the positive and negative parts of $\xi(\bx,t;\X,\cdot)$)
% denoted $|\xi(\bx,t;\X)|$. 

Suppose $(\bx,t)\in \R^d \times \MM$
and $\X \subset \R^d \times \MM$ is finite.
If $(\bx,t) \notin \X$, we abbreviate notation and write $\xi(\bx,t; \X)$
instead of $\xi(\bx,t; \X \cup \{(\bx,t)\})$. We also write
$\X^{\bx,t}$ for $\X \cup \{(\bx,t)\}$.

Given $a > 0$ and $\y \in \R^d$, we let $\y+ a\X:=
\{(\y+ a\bx,s): (\bx,s) \in \X\}$;
% and $y+ \X := \{(\y+\bx,s): (\bx,s) \in \X \}$;
in other words,
scalar multiplication and translation act on
only the first component of elements of
 $\R^d \times \MM$. 
We say $\xi$ is {\em translation invariant}  if
 $$
\xi(\bx,t; \X,A) = \xi(\y+ \bx,t ; \y+ \X,\y+ A )
$$
 for all $\y \in \R^d$, all finite $\X \subset \R^d \times \MM$ and $\bx \in
\X$, and all Borel $A \subseteq \R^d$.
 Some of the general concepts defined in the sequel can
be expressed more transparently when $\xi$ is translation invariant.

Let $\kappa $ be a probability  density function
on $\R^d$. Abusing notation slightly,  we also let $\ka$ denote the 
 corresponding  probability measure on $\R^d$,
i.e. we write  $\ka(A) $ for $\int_A \kappa (\bx) d\bx$,
for Borel $A \subseteq \R^d$.
% $\ka$ with compact support $\Gamma \subset \R^d$, 
 We shall assume throughout that the density function
 $\ka$ is bounded with supremum
denoted $\|\kappa\|_\infty$, 
%Also we assume (**perhaps) that
and that $\kappa$ is Lebesgue-almost everywhere continuous. Let
$\supp(\kappa)$ denote the support of $\kappa$, i.e., the
smallest closed set $B$ in $\R^d$ with $\ka(B)=1$.

For all $\la > 0$, let $\P_{\la}$ denote a
 Poisson point process in $\R^d \times \MM$ with intensity measure 
$\la \ka \times \mu_\MM$.
 For $n \in \N$, let $\X_n$ be the point process
consisting of $n$ independent identically distributed
random elements of $\R^d \times \MM$ with common distribution
given by $\kappa \times \mu_\MM$.
Let $\H_\la$ denote a  Poisson point process in
$\R^d \times \MM$ with intensity $\la$  times the product of
$d$-dimensional Lebesgue measure 
and $\mu_\MM$. 
%For $(\bx,t) \in \R^d \times \MM$ we write
%$\Po_\la^{\bx,t}$ for $\Po_\la \cup \{(\bx,t)\}$ and
%$\H_\la^{\bx,t}$ for $\H_\la \cup \{(\bx,t)\}$.

Suppose we are given a family of Borel
subsets $\Gamma_\la$ of $\R^d$, indexed by $\la \geq 1$.
Assume the sets $\Gamma$ are nondecreasing in $\la$,
i.e. $\Gamma_\la \subseteq \Gamma_{\la'}$ for $\la < \la'$. 
Denote by $\Gamma $ the limiting set, i.e. set $\Gamma \:= \cup_{\la \geq 1}
\Gamma_\la$.
%We assume that there exists a limiting set $\Gamma \subseteq \R^d$,
%such that as $\la \to \infty$, we have
%$\Gamma_\la \to \Gamma$, in the sense that
% ${\bf 1}_{\Gamma_\la} \to {\bf 1}_{\Gamma} $ pointwise.
We assume that the topological boundary of $\Gamma$
(i.e., the intersection of the closure of $\Gamma$ with that of its
complement) has zero Lebesgue measure, and that 
$\inter(\Gamma) = \cup_{\la \geq 1} \inter(\Gamma_\la)$, where
 $\inter(\cdot)$ denotes interior.
%The simplest special case has $\Gamma_\la =\ska$ for all $\la$.
The simplest special case has $\Gamma_\la =\R^d$ for all $\la$.

For $\la > 0$, and for finite  $\X \subset \R^d \times \MM$ with 
$(\bx,t) \in \X$, and Borel $A \subset \R^d$, let
$$
\xi_{\la}(\bx,t; \X,A) : = \xi (\bx,t; \bx+ \la^{1/d}(-\bx +\X),
\bx + \la^{1/d} (-\bx +A ) )
{\bf 1}_{\skala}(\bx).
$$
  When $\xi$ is
 translation invariant,  the rescaled measure $\xi_\la$  simplifies to
\bea
\xi_{\la}(\bx,t; \X,A) = \xi (\la^{1/d}\bx,t; \la^{1/d}\X,\la^{1/d}A )
\1_{\Gamma_\la}(\bx).
\lbl{TIxila}
\eea

Our principal objects of interest are the
random 
%signed
  measures ${\mu}_{\la }^{\xi}$ and  $\nu_{\la,n}^{\xi}$
 on $\R^d$, defined for $\la >0$ and $n \in \N$, by
$$
{\mu}_{\la }^{\xi}:= 
\sum_{(\bx,t) \in {\cal P}_{\la}
%\cap (\skala \times \MM)
}
 \xi_{\la}(\bx,t;{\cal
P}_{\la }) ; 
%$$
~~~~~
%$$
{\nu}_{\la,n }^{\xi}:= 
\sum_{(\bx,t) \in {\X}_{n}
%\cap (\Gamma_\la \times \MM)
} \xi_{\la}(\bx,t;{\cal X}_{n }) .
$$
%where $\delta_\bx$ denotes a unit point mass at $\bx$.
We are also interested in
 the centred  versions of these measures 
 $\overline{\mu}^{\xi}_{\la } :=
{\mu}_{\la }^{\xi} - \E[{\mu}_{\la}^{\xi}]$ and
 $\overline{\nu}^{\xi}_{\la,n } :=
{\nu}_{\la,n }^{\xi} - \E[{\nu}_{\la,n }^{\xi}]$
(which are signed measures).
We study these measures via their action on
test functions in the space
 $B(\R^d)$
of bounded Borel-measurable functions on $\R^d$.
We let $\tB(\R^d)$ denote the subclass of $B(\R^d)$ consisting
of those functions that are Lebesgue-almost everywhere
continous.

Given
$f \in B(\R^d)$, set $\langle f, \xi_\la(\bx,t;\X) \rangle := \int_{\R^d}
f(\bz) \xi_\la(\bx,t;\X,d\bz)$. Also, set
$$
\langle f, {\mu}^{\xi}_{\la } \rangle := \int_{\R^d} f d
{\mu}^{\xi}_{\la } = \sum_{(\bx,t) \in \Po_\la} \langle f,\xi_\la
(\bx,t;\Po_\la)
\rangle
$$
and set
$\langle f, \overline{\mu}^{\xi}_{\la } \rangle := \int_{\R^d} f d
\overline{\mu}^{\xi}_{\la }$, so that 
$\langle f, \overline{\mu}^{\xi}_{\la } \rangle = 
\langle f, {\mu}^{\xi}_{\la } \rangle - 
\E \langle f, {\mu}^{\xi}_{\la } \rangle$.
Similarly, let 
$\langle f, {\nu}^{\xi}_{\la ,n} \rangle := \int_{\R^d} f d \nu_{\la,n}^\xi$
and 
$\langle f, \overline{{\nu}}^{\xi}_{\la ,n} \rangle := 
\langle f, {\nu}^{\xi}_{\la, n } \rangle - 
\E \langle f, {\nu}^{\xi}_{\la, n } \rangle.
$

% and their centred versions
%$\bar{{\mu}}_{\la }^{\xi}$ and  $\bar{\nu}_{\la,n}$,

Let $|\cdot|$ denote the Euclidean norm on $\R^d$, and
 for $\bx \in \R^d$ and $r > 0$, define the ball
 $B_r(\bx):= \{\by \in \R^d:|\by-\bx| \leq r\}$.
% denote the closed Euclidean ball
%centred at $\bx$ of radius $r$. 
We denote by $\0$  the origin of
$\R^d$ and abbreviate $B_r(\0)$ to $B_r$.

We say a set $\X \subset \R^d \times \MM$ is 
{\em locally finite} if $\X \cap (B \times \MM)$ is
finite for all bounded  $B \subset \R^d$.
For $(\bx,t) \in \R^d \times \MM$ and 
Borel $A \subseteq \R^d$,
we extend the definition of 
%$\xi(\bx,t;\bx+\X,\bx+ A)$ to locally finite infinite point
$\xi(\bx,t;\X,A)$ to locally finite infinite point
sets $\X $ by setting
$$
\xi(\bx,t;\X,A) := \limsup_{K \to \infty} \xi(\bx,t ; 
 [\X \cap (B_K \times \MM)],  A ).
%\xi_\infty^{\bx,t}(\X,A) = \limsup_{K \to \infty} \xi(\bx,t ; 
%\bx+ [\X \cap (B_K \times \MM)], \bx + A ),
$$
Also, we define the $\bx$-shifted version $\xi_\infty^{\bx,t}(\cdot,\cdot)$
 of $\xi(\bx,t;\cdot,\cdot)$ by 
$$
\xi_\infty^{\bx,t}(\X,A) = 
%\limsup_{K \to \infty} 
\xi(\bx,t ; \bx + \X , \bx + A ).
$$
Note that if $\xi$ is translation-invariant then
$\xi_\infty^{\bx,t}(\X,A)= \xi_\infty^{\0,t}(\X,A)$ for all
$\bx \in \R^d$, $t \in \MM$, and Borel $A \subseteq \R^d$.

Notions of {\em stabilization}, introduced in
 \cite{PY1,PY4,BY}, play a central role in all that follows.

\begin{defn}
\lbl{RSdef}
For any locally finite $\X \subset \R^d \times \MM$, 
any $(\bx,t) \in \R^d \times  \MM$,
and any Borel region $A \subset \R^d$, define 
$R(\bx,t;\X,A)$ 
(the {\em radius of stabilization of $\xi$} at $(\bx,t)$ with respect to $\X$
and $A$) to be the smallest
integer-valued $r$ such that 
$r \geq 0$ and
$$
\xi(\bx,t; \bx+ ([ \X \cap (B_r \times \MM)]
 \cup \Y ),B ) = \xi(\bx,t; \bx + [\X \cap (B_r \times \MM)], B )
$$
for all finite $\YY  \subseteq (A \setminus B_r)\times \MM$
and all Borel $B \subseteq \R^d$.
If no such $r$ exists, we set  $R(\bx,t;\X,A) = \infty$.

When $A$ is the whole of $\R^d$
we abbreviate the notation $R(\bx,t;\X,\R^d)$ to $R(\bx,t;\X)$.
\end{defn}

In the case where $\xi$ is translation-invariant, 
 $R(\bx,t;\X) =R(\0,t;\X)$ so that 
$R(\bx,t;\X)$ does not depend on $\bx$.
Of particular importance to us will be radii of stabilization
with respect to the homogeneous Poisson processes
$\H_{\la}$ and  with respect to the  non-homogeneous
Poisson process $\Po_\la$, suitably scaled.

We assert that 
 $R(\bx,t;\X,A)$ is a measurable function of $\X$,
and hence, when $\X$ is a random point set such as $\H_\la$ or $\Po_\la$,
 $R(\bx,t;\X,A)$ is an $\N \cup \{\infty\}$-valued random variable.
To see the assertion, observe that
 by Dynkin's pi-lambda lemma, for any $k\in \N$ the event 
$\{R(\bx,t;\X,A) \leq k \}$  equals the event $\cap_{B \in {\cal B}} \{
s(\X,B) = i(\X,B)\},$ where
${\cal B}$ is the $\Pi$-system consisting 
of the rectilinear hypercubes in $\R^d$ whose
 corners have
 rational coordinates, and for $B \in {\cal B}$ we set
\bean
 %\{ 
s(\X,B) := \sup \{ \xi(\bx,t; \bx+ ([\X \cap ( B_k \times \MM)] \cup \YY ),B):
\Y \subseteq (A \setminus B_k )\times \MM \}, 
\\
i(\X,B) := \inf \{ \xi(\bx,t; \bx+ ([\X \cap ( B_k \times \MM)] \cup \YY ),B):
\Y \subseteq (A \setminus B_k )\times \MM \}. 
\eean
Also, $s(\X,B)$ is a measurable function of $\X$ because  we assume
$\xi$ is Borel-measurable and for any $b$ we have
$$
\{ \X: s(\X)  > b \} = \pi_1 (\{(\X,\Y): \xi 
(\bx,t; \bx+ [\X \cap (B_k \times \MM) ] \cup  
[\Y \setminus (B_k \times \MM) ], B) > b \})
$$  
where $\pi_1$ is projection onto the first component,
acting on pairs $(\X,\Y)$ with $\X$ and $\Y$ finite sets in
$\R^d \times \MM$. Similarly, $i(\X,B)$ is a measurable function of $\X$.

%If $R(\bx,t;\X,\R^d)$ is finite, then $\xi_\infty^{\bx,t}(\X,\cdot)$
%is a signed measure with finite total variation, and in this case,
%for Borel $B \subseteq \R^d$,  we write
%$|\xi_\infty^{\bx,t}|(\X,B)$
%for the total variation (i.e. the sum of the positive and negative
%parts) of the restriction
%of the measure $\xi_\infty^{\bx,t}(\X,\cdot)$ to $B$.

\begin{defn} 
Let $T$, $T'$, $T''$ and $T'''$ denote generic random
elements of $\MM$ with distribution $\mu_\MM$, independent
of each other and of all other random objects
we consider. Similarly, let 
 $\bX$ and $\bX'$ denote generic random $d$-vectors
 with  distribution $\kappa$, independent
of each other and of all other random objects
we consider. 
\end{defn} 
For $p>0$, we consider $\xi$ satisfying the moments
conditions
\begin{equation} 
\lbl{mom}
  \sup_{\la \geq 1, \ \bx \in
\ska} 
%\E [ | \xi_{\la}(\bx,T; {\cal P}_{\la },\R^d ) |^p ]
\E [  \xi_{\la}(\bx,T; {\cal P}_{\la },\R^d )^p ]
< \infty.
\end{equation}
and
\begin{equation} \lbl{mom2}
  \sup_{\lambda \geq 1, \ \bx, 
%\in \ska,
 \by \in \ska}
 %\E [ | \xi_{\la}(\bx,T; {\cal P}_{\la } \cup \{(\y,T')\} \R^d) |^p ]
 \E [  \xi_{\la}(\bx,T; {\cal P}_{\la } \cup \{(\y,T')\}, \R^d)^p ]
< \infty.
\end{equation}

Let $\xi^*(\bx,t;\X,\cdot)$ be the point measure
at $\bx$ with the same total measure as 
$\xi(\bx,t;\X,\cdot)$, 
% be the point measure signed measure
%consisting of a single point measure  at $\bx$
%with total value equal to $\xi(\bx,t;\X,\R^d)$, 
i.e.,
for Borel $A \subseteq \R^d$  let
\bea
\xi^*(\bx,t;\X,A) := \xi(\bx,t;\X,\R^d) {\bf 1}_A(\bx).
\lbl{eq2.2a}
\eea
%Finally we make some further assumptions about the
%functional $\xi$ and test function $f$. 
We  consider measures $\xi$ and 
test functions $f \in B(\R^d)$ satisfying one of
the following assumptions:
\begin{enumerate}
\item[A1:] $\xi = \xi^*$, i.e., $\xi(\bx,t;\X,\cdot)$ is a point mass at
$\x$ for all $(\bx,t,\X)$.
\item[A2:] 
$\xi(\bx,t;\X,\cdot)$ 
is absolutely continuous
with respect to Lebesgue measure on $\R^d$, 
 with Radon-Nikodym derivative
denoted 
$\xi'(\bx,t;\X,\by)$ for $\by \in \R^d$,  
satisfying  
%$|\xi'(\bx,t;\X,\by)|\leq K_0$ for all (signed)    
$\xi'(\bx,t;\X,\by)\leq K_0$ for all    
$(\bx,t;\X,\by)$, where $K_0$ is a finite positive constant.
\item[A3:] $f$ is almost everywhere continuous, i.e., $f \in \tB(\R^d)$.
\end{enumerate}
Note that Assumption A2 will hold if $\xi(\bx,t,\X,\cdot)$ 
is Lebesgue measure on some random subset of $\R^d$ determined
by $\bx,t,\X$.

Our first result is a law of large numbers
for
 $\langle f, \mu_\la^\xi \rangle$ and
 $\langle f, \nu_{\la,n}^\xi \rangle$, for
$f \in B(\R^d)$. This 
 extends 
a result in \cite{PY4} where only the case where $f$ is a constant
is considered.
 We require almost surely finite radii of stabilization
with respect to homogenous Poisson processes, along with
moments conditions.

\begin{theo}
\lbl{thlln1}
Suppose that  
$R(\bx,T;\H_{\ka(\bx)})$ is almost surely finite for
 $\ka$-almost all $\bx \in \Gamma$. Suppose also that
$f \in B(\R^d)$,
 and one or more of assumptions A1, A2, A3 holds. 
Let $q =1$ or $q=2$.
Then: 
%\bea
%\int_{\R^d} P[R(\bx,T;\H_{\ka(\bx)})= \infty] \ka(\bx) d\bx =0.
%\lbl{homstab0}
%\eea
%Then:

(i) If
% Assumption A1, A2 or A3 holds and 
 there exists $p>q$ such that
 \eq{mom} and \eq{mom2} hold, then
%for $q=2$ 
 %for all $f  \in B(\R^d)$ 
as $\la \to \infty$ 
we have 
\bea
\la^{-1} \langle f, \mu_\la^\xi \rangle \inLq \int_{\Gamma} f(\bx) \E
\xi_\infty^{\bx,T}(\H_{\ka(\bx)},\R^d) \ka(\bx) d\bx
\lbl{LLNpo}
\eea
and
\bea
\la^{-1} \sum_{(\bx,t) \in \Po_\la} |\langle f,\xi_\la(\bx,t;\Po_\la)
- \xi_\la^*(\bx,t;\Po_\la) \rangle | \inL 0.
\lbl{LLN2po}
\eea

%(ii) If A1 or A2 holds and there exists $p>1$ such that 
% \eq{mom} and \eq{mom2} hold, then \eq{LLNpo} holds for $q=1$.

(ii)
 %If A1, A2, or A3 holds, and 
If
$\la(n)/n \to 1$ as $n \to \infty$, and
%for $q=1$ or $q=2$ 
there exists $n_0>0$ and $p >q$ such that
\bea
\sup_{n \geq n_0} \E [ 
%|\xi_{\la(n)}(\bX,T;\X_{n-1},\R^d)|^p] <\infty ,  (signed)
\xi_{\la(n)}(\bX,T;\X_{n-1},\R^d)^p] <\infty ,
%\1\{\bX \in \Gamma_\la\}]< \infty,
\lbl{momlnb} 
\eea
then 
%for all  $f \in B(\R^d)$ 
%for $q=2$,
we have as $n \to \infty$ that 
\bea
n^{-1} \langle f, 
\nu_{\la(n),n}^\xi \rangle \inLq \int_{\Gamma} f(\bx) \E
\xi_\infty^{\bx,T}(\H_{\ka(\bx)},\R^d) \ka(\bx)d\bx.
\lbl{LLNbi}
\eea
and
\bea
n^{-1} \sum_{(\bx,t) \in \X_n} |\langle f,\xi_\la(\bx,t;\X_n)
- \xi_\la^*(\bx,t;\X_n) \rangle | \inL 0.
\lbl{LLN2bi}
\eea
%
%(ii)
%If A1 or A2 holds, and 
%$\la(n)/n \to 1$ as $n \to \infty$, and
%there exists $n_1>0$ and $p >1$ such that
%\eq{momlnb} holds for $n \geq n_1$,
%then \eq{LLNbi} holds for $q=1$.
%
\end{theo}

Our main results are 
central limit theorems to go with the laws of large numbers. 
For these we need further conditions.  The first of these requires
finite radii of stabilization with respect to homogeneous Poisson processes,
possibly with a point inserted,  and,
 in the non translation-invariant case,
also requires local tightness of 
%the probability distributions of
 these radii. 
%along with local tightness of the limiting measures
% $\xi_\infty^{\bx,T}$

%allowing 
%some  
% (\ref{homstab0})
%by requiring some local tightness of these radii as $\bx$ varies,  of  
%the radii of stabilization   at $\bx$ with respect to $\H_{\ka(\bx)}$.

\begin{defn}\lbl{stab0}
 For $\bx \in \R^d$ and $\la >0$,
we shall say  that $\xi$ is  {\em
$\la$-homogeneously stabilizing} at $\bx$
 if 
%there exists a neighbourhood $\NN_{\bx}$ of $\bx$ such that
 for all $\bz \in \R^d$,
\bea
P\left[ \lim_{\eps \downarrow 0} \sup_{\by \in B_\eps(\bx)}
\max (R(\by,T; \H_\la),R(\by,T;\H_{\la}^{\bz,T'})) < \infty \right] =1. 
%\lim_{t \to \infty}
%\sup_{\by \in \NN_{\bx}}
%P[ R(\by,T;\H_{\la}) >t] 
%=
%\lim_{t \to \infty}
%\sup_{\by \in \NN_{\bx}}
% P[R(\by,T;\H_{\la}^{\bz,T'})>t ]
% = 0,
\lbl{homstab2}
\eea
%and also, for any $\eps >0$, $\bz \in \R^d$, it is the case that
%\bea
%\lim_{t \to \infty} \sup_{\by \in \NN_\bx}
%P[|\xi_\infty^{\by,T}| (\H_\la,\R^d \setminus B_t) > \eps ] = 0;
%\lbl{homstab3}
%\eea
%\bea
%\lim_{t \to \infty} \sup_{\by \in \NN_\bx}
%P[|\xi_\infty^{\by,T}| (\H_\la^{\bz,T'},\R^d \setminus B_t) > \eps ] = 0.
%\lbl{homstab4}
%\eea
\end{defn}
%{\bf Remark.}
In the case where $\xi$ is translation-invariant, 
 $R(\bx,t;\X)$ does not depend on $\bx$, 
and $\xi_\infty^{\bx,T}(\cdot)$ does not depend on $\bx$ so that
 the simpler-looking condition
\bea
P[ R(\0,T;\H_{\la})< \infty] = 
 P[R(\0,T;\H_{\la}^{\bz,T'})< \infty]
 = 1,
\lbl{homstab}
\eea
suffices to guarantee condition
 (\ref{homstab2}).
% \eq{homstab3}, and \eq{homstab4}.
%Moreover, under assumption A1, that $\xi=\xi^*$, conditions \eq{homstab3}
%and \eq{homstab4} always hold if \eq{homstab2} holds.

We now introduce notions of exponential and power-law stabilization.
The terminology refers to the tails of the distributions
of  radii of stabilization 
with respect to
the non-homogeneous point processes $\Po_\la$ and $\X_n$.

For $k=2$ or $k=3$, let $\SS_k$ denote the set of all finite $\A
\subset \ska \times \MM$ with at most $k$ elements (including
the empty set), and for nonempty $\A \in \SS_k$, let $\A^*$
denote the subset  of $\ska \times \MM$ (also with $k$ elements)
obtained by equipping each element  of $\A$ with a $\mu_\MM$-distributed
mark; for example, for $\A= \{\bx,\by\} \in \SS_2$ set $\A^*=
\{(\bx,T'),(\by,T'')\}$.

\begin{defn}\lbl{stab}
For $\bx \in \R^d$, $\la >0$ and $n \in \N$,
and $\A \in \SS_2$, 
define the $[0,\infty]$-valued random variables $R_\la(\bx,T)$ and
$R_{\la,n}(\bx,T;\A)$  by
\bea
R_\la(\bx,T) = R(\bx,T;\la^{1/d}(-\bx+ \Po_\la),\la^{1/d}(-\bx+ \ska) ),
\lbl{5311a}
\\
R_{\la,n}(\bx,T;\A) = %\sup_{{\cal A}} 
R(\bx,T;\la^{1/d}(-\bx+ (\X_n \cup \A^*)) ,\la^{1/d}(-\bx+ \ska) ).
\lbl{5311b}
\eea
%where the supremum is over all sets 
%${\cal A} \subset \supp(\kappa)$ with
%$0,1$ or $2$ elements.
When $\A$ is the empty set $\emptyset$ we write $R_{\la,n}(\bx,t)$ for
$R(\bx,t;\emptyset)$.

For $s>0$ and $\eps \in (0,1)$ define the tail probabilities $\tau(s)$ 
and $\tau_\eps(s)$by
\bean
\tau(s) = \sup_{\la \geq 1,\ \bx \in \skala}
 P[R_\la(\bx,T) > s] ;
%~~~~~
\\
\tau_\eps(s) = \sup_{\la \geq 1, n \in \N \cap ((1-\eps)\la,(1+\eps)\la),
 \bx \in \skala,\A \in \SS_2}
 P[R_{\la,n}(\bx,T;\A^*) > s] .
\eean
Given $q>0$,
we say  $\xi$ is {\em power-law stabilizing of order} $q$ for $\ka$ if 
$\sup_{s \geq 1} s^{q} \tau(s)<\infty$.
We say  $\xi$ is {\em exponentially stabilizing} for $\ka$ if 
$\limsup_{s \to \infty} s^{-1} \log \tau(s) <0$.
We say  $\xi$ is {\em binomially power-law stabilizing of order}
 $q$ for $\ka$ if  there exists $\eps >0$ such that
$\sup_{s \geq 1} s^{q} \tau_\eps(s)<\infty$.
We say  $\xi$ is {\em binomially exponentially stabilizing} for $\ka$ if 
there exists $\eps >0$ such that
$\limsup_{s \to \infty} s^{-1} \log \tau_\eps(s) <0$.
\end{defn}
It  is easy to see that if $\xi$ is exponentially stabilizing
for $\ka$
then it is power-law stabilizing of all orders for $\ka$.
Similarly, if $\xi$ is binomially exponentially stabilizing
for $\ka$
then it is binomially power-law stabilizing of all orders for $\ka$.

In the non translation-invariant case, we shall also requre the following
continuity condition.

\begin{defn}
\lbl{defctsmeas}
For $\bx \in \R^d$,
we say $\xi $ has
 {\em almost everywhere continuous total measure} if 
there exists $K_1>0$ such that
for all $m\in \nat$ and
 Lebesgue-almost all $(\bx,\bx_1,\ldots,\bx_m) \in (\R^d)^{m+1}$, and
 $(\mu_\MM \times \cdots \times \mu_\MM)$-almost all
$(t,t_1,t_2,\ldots, t_m)\in \MM^{m+1}$, 
for $A=\R^d$ 
 or $A=B_K$ with  
 $K>K_1$, 
 the function
$$
(\by,\by_1,\by_2,\ldots,\by_m)
 \mapsto  \xi(\by,t;\{(\by_1,t_1),(\by_2,t_2),\ldots,(\by_m,t_m)\},\by +A)
$$
is continuous at $(\by,\by_1,\ldots,\by_m) = (\bx,\bx_1,\ldots,\bx_m)$.
%
%and also the function
%$$
%(\by,\by_1,\by_2,\ldots,\by_m)
% \mapsto |  \xi(\by,t;\{(\by_1,t_1),(\by_2,t_2),\ldots,(\by_m,t_m)\})|
%$$
%is continuous at $(\by,\by_1,\ldots,\by_m) = (\bx,\bx_1,\ldots,\bx_m)$.
\end{defn}
Subsequent results will require one of the follwing 
 formal assumptions to hold. 
\begin{enumerate}
\item[A4:]
$\xi$ is translation-invariant.
\item[A5:] 
$\xi $ has almost everywhere continuous total meaure.
\end{enumerate}

Our next result gives the asymptotic variance
of $\langle f, {\mu}_\la^\xi \rangle$ for $f \in B(\R^d)$. 
\begin{theo}
\lbl{thm1}
Suppose $\xi$ is 
%almost everywhere continuous, and is 
$\kappa(\bx)$-homogeneously stabilizing for $\ka$-almost all
$\bx \in \R^d$.
Suppose also that $\xi$ satisfies the moments conditions (\ref{mom})
and (\ref{mom2}) for some $p>2$, 
and is  power-law stabilizing for 
%$(\kappa,\Gamma)$.
$\kappa$ of order $q $ for some $q$  with $q>p/(p-2)$.
Let $f \in B(\R^d)$,  and assume that
  Assumption A1 or  A3 holds, and that  Assumption A4 or A5 holds.
Then
 \bea
\lbl{BYeq}
\lim_{\la \to \infty}
\la^{-1}  \Var [
\langle f, \mu^{\xi}_{\la } \rangle ]  =
\int_{\Gamma} f(\bx)^2 V^\xi(\bx,\kappa(\bx)) \kappa(\bx) d\bx,
\eea
with
$V^\xi(\bx,a)$
given by $V^\xi(\bx,0)=0$, and for $a >0$ by
\bea
V^\xi(\bx,a) = \E[\xi_\infty^{\bx,T}(\H_a,\R^d)^2] 
\nonumber \\ 
+ a \int_{\R^d} (
\E[\xi_\infty^{\bx,T}(\H_a^{\bz,T'},\R^d)\xi_\infty^{\bx,T'}
 (-\bz +\H^{\0,T}_{a},\R^d) ]
- (\E[\xi_\infty^{\bx,T}(\H_a,\R^d) ])^2
)d\bz.
\lbl{Vdef}
\eea 
Also, the integral in (\ref{Vdef}) converges
for $\ka$-almost all $\bx$,  and the right hand side
of (\ref{BYeq}) is finite.
\end{theo}
%The conclusions of Theorem \ref{thm1} 
%imply the conclusion of the
%laws of large numbers in Theorem \ref{thlln1} (ii) (with $q=2$).
%However, Theorem \ref{thlln1} is not redundant since
%the conditions in Theorem \ref{thm1} are stronger, in
%particular the fact that a further kind of stabilization is needed.

Our next result is a central limit theorem for the random
field  
$(\la^{-1/2}\langle f, \overline{\mu}_{\la }^{\xi}\rangle, f \in \tilde{B}(\R^d))$,
or 
$(\la^{-1/2}\langle f, \overline{\mu}_{\la }^{\xi}\rangle, f \in B(\R^d))$.
% obtained by combining Theorem \ref{thm1}
%with a normal approximation result
%for exponentially stabilizing functionals 
%given in \cite{PY5}.
We list some further assumptions.
\begin{enumerate}
\item[A6:] 
For some $p>2$,
 $\xi$ satisfies the moments conditions (\ref{mom})
 and (\ref{mom2}) 
 and
 is exponentially stabilizing for $\kappa$.  
\item[A7:]
For some $p>3$,
 $\xi$ satisfies the moments conditions (\ref{mom})
 and (\ref{mom2}) 
 and is power-law  stabilizing for 
$\kappa$ of order $q$ for some $q > d(150+6/p)$. 
\end{enumerate}

\begin{theo}
\lbl{thm2}
Suppose $\|\kappa\|_\infty < \infty$
 and $\kappa$ has bounded support.
Suppose $\xi$ 
%is almost everywhere continuous, and
 is
 $\ka(\bx)-$homogeneously stabilizing at $\bx$ for $\ka$-almost
all $\bx \in \R^d$,
% Suppose $\xi$
 satisfies either A$4$ or A$5$, and satisfies either
A$6$ or A$7$.
 %Suppose also {\em either} that
 %$\xi$ satisfies the moments conditions (\ref{mom})
 %and (\ref{mom2}) for some $p>2$ and
 %is exponentially stabilizing for $\kappa$,  
%{\em or} 
% that 
%$\xi$ satisfies the moments conditions (\ref{mom})
% and (\ref{mom2}) for some $p>3$,
% and is power-law  stabilizing for 
%$\kappa$ of order $q$ for some $q > d(150+6/p)$. 

 Then
 the finite-dimensional distributions of the random field
\linebreak
$(\la^{-1/2}\langle f, \overline{\mu}_{\la }^{\xi}\rangle, f \in \tB(\R^d))$
 converge weakly 
as $\la
 \to \infty$  to those of a mean-zero finitely additive
Gaussian field
with covariances given by
 $ \int_{\R^d} f_1(\bx) f_2(\bx) V^{\xi}(\bx,\ka(\bx))
\ka(\bx) d\bx,$ with $V^\xi(\bx,a)$ given by (\ref{Vdef}). 

If also $A1$ holds,
 the finite-dimensional distributions of the random field
\linebreak
$(\la^{-1/2}\langle f, \overline{\mu}_{\la }^{\xi}\rangle, f \in B(\R^d))$
 converge weakly 
%in law 
as $\la
 \to \infty$  to those of a mean-zero finitely additive
Gaussian field
with covariances given by
 $ \int_{\R^d} f_1(\bx) f_2(\bx) V^{\xi}(\bx,\ka(\bx))
\ka(\bx) d\bx$.
% with $V^\xi(\bx,\la)$ given by (\ref{Vdef}). 
\end{theo}
The corresponding results for the random
 measures $\overline{\nu}_{\la,n}^\xi$ require some further conditions.
These  extend the previous stabilization
and moments conditions to binomial point processes.
Our 
%first
 extra moments condition is
\bea
\inf_{\eps >0}
\sup_{\la \geq 1,\bx \in \R^d, \A \in \SS_3} 
%\sup_{\bx \in \skala} \sup_\A
 \sup_{(1-\eps) \la \leq m \leq (1+\eps) \la}
\E \left[ 
 % |\xi_\la(\bx,T;\X_m \cup \A )|^p 
 \xi_\la(\bx,T;\X_m \cup \A^*,\R^d )^p 
%\1_{\skala}(\bX) 
\right]
< \infty,
\lbl{mom4}
\eea
%where the supremum $\sup_\A$ is over sets
%$\A \subset \supp(\kappa)$ having $0,1,2,$ or $3$ elements,
%each element carrying 
%an independent mark with distribution $\mu_\MM$. 

%Our second extra moments condition requires there to be
%a constant $\beta_1$ such that
%\bea
%\E[ \xi_\la(\bX,T;\X_m,\R^d)^2]
%  \leq \beta_1( 1 + \la + m)^{\beta_1}, ~~~ 
%\forall \la \geq 1, m \in \{0,1,2,\ldots\} 
%\lbl{polybd}
%\eea

We give strengthened versions of A6 and A7 above, to include
condition \eq{mom4} and binomial stabilization.
\begin{enumerate}
\item[A$6'$:] 
For some $p>2$,
 $\xi$ satisfies the moments conditions (\ref{mom}),
(\ref{mom2}),  \eq{mom4}, 
%\eq{polybd},
 and
 is exponentially stabilizing for $\kappa$ and binomially exponentially
stabilizing for $\kappa$.  
\item[A$7'$:] 
For some $p>3$,
$\xi$ satisfies the moments conditions (\ref{mom}),
 (\ref{mom2}) and \eq{mom4},
% and \eq{polybd},
 and is power-law  stabilizing and 
 binomially power-law stabilizing
 for 
$\kappa$ of order $q$ for some $q > d(150+6/p)$.
%and is
% for $\kappa$ of order $q'$ for
%some $q'> 2dp/(p-2)$. 
\end{enumerate}

For $\bx \in \R^d$ and $a >0$, set 
\bea
\delta(\bx,a) := \E[ \xi_\infty^{\bx,T} (\H_{a},\R^d) ] +
a \int_{\R^d} \E [ 
\xi_\infty^{\bx,T}(\H_{a}^{\by,T'},\R^d)  
- \xi_\infty^{\bx,T}(\H_{a},\R^d) ] d\by.  
~~~
\lbl{deltadef}
\eea

%We give a limit theorem for $\overline{\nu}_{\la(n),n}^\xi$. 
\begin{theo}
\lbl{thm3}
Suppose $\|\kappa\|_\infty < \infty$
 and $\kappa$ has bounded support.
Suppose $\xi$ 
%is almost everywhere continuous, 
is $\ka(\bx)-$homogeneously stabilizing at $\bx$ for $\ka$-almost
all $\bx \in \R^d$, 
satisfies Asuumption A$4$ or A$5$, and also satisfies
%Suppose also that either 
A$6'$ or A$7'$.
%and that $\xi$ is polynomially bounded
%\eq{polybd}.
% Suppose also {\em either} that $\xi$
% satisfies the moments conditions (\ref{mom})
% and (\ref{mom2}) for some $p>2$ and
% is binomially exponentially stabilizing for 
%$\kappa$,  
%{\em or} 
% that $\xi$
% satisfies the moments conditions (\ref{mom})
% and (\ref{mom2}) for some $p>3$. 
%Assume
%also that $\xi$ is binomially power-law stabilizing 
%for $\kappa$ of order $q$ for some $q > d(150+6/p)$ with
%also $q >2 dp/(p-2)$. 
%Finally, suppose the further moments conditions \eq{mom4} holds.
% \eq{mom4}, \eq{mom5}, \eq{mom6}   hold.

Then for any sequence $(\la(n),n \in \N)$
taking values in $(0,\infty)$,
such that \linebreak
 $\limsup_{n \to \infty} n^{-1/2}|\la(n) -n|  < \infty$,
we have for $f \in \tB(\R^d)$ 
%under any of assumptions A1, A2 or A3  
that
\bea
\lim_{n \to \infty} n^{-1} \Var \langle f,\nu_{\la(n),n}^\xi \rangle 
= \int_{\Gamma} f(\bx)^2 V^\xi(\bx,\ka(\bx)) \ka(\bx) d\bx 
\nonumber \\
- \left( \int_{\Gamma} f(\bx) \delta(\bx,\ka(\bx)) \ka(\bx) d\bx \right)^2, 
\lbl{0824b}
\eea
and the finite-dimensional distributions of the random field
$(n^{-1/2}\langle f, \overline{\nu}_{\la(n),n }^{\xi}\rangle,
 f \in \tB(\R^d))$
 converge weakly 
as $n \to \infty$  to 
those of a mean-zero
 finitely additive Gaussian field
with covariances given by
\bea
 \int_{\Gamma} f_1(\bx) f_2(\bx) V^{\xi}(\bx,\ka(\bx))
\ka(\bx) d\bx
\nonumber \\
- \int_{\Gamma} f_1(\bx) 
\delta(\bx,\ka(\bx))
\ka(\bx)d\bx
 \int_{\Gamma} f_2(\by)
 \delta(\by,\ka(\by)) 
\ka(\by)d\by
\lbl{0209}
\eea
 with $V^\xi(\bx,\la)$ given by (\ref{Vdef}).

If in addition, assumption A$1$ 
 holds, then \eq{0824b} holds for $f \in B(\R^d)$ and
 the finite-distributions of the the random field
$(\la(n)^{-1/2}\langle f, \overline{\nu}_{\la(n),n }^{\xi}\rangle,
 f \in B(\R^d))$
 converge weakly 
as $n \to \infty$  to those of a mean-zero finitely additive Gaussian
field  with covariances given by the expression in \eq{0209}.
\end{theo}

{\bf Remarks.}
Theorems \ref{thm1}, \ref{thm2} and \ref{thm3} resemble
 the main results
of Baryshnikov and Yukich \cite{BY}, in that they provide central
limit theorems for random measures under general stabilization 
conditions. 
We indicate here some of the 
ways in
which our results extend those in \cite{BY}.

In \cite{BY}, attention is restricted to
cases where assumption  A1 holds, i.e., where 
the contribution from each point to the random measures is a point
mass at that point. It is often  natural   
to drop this restriction, for example when considering the volume
or surface measure associated with a germ-grain model, examples
we shall consider in detail in Section \ref{secgerm}.

Another difference is that under A1, we consider
 bounded test functions in $B(\R^d)$ whereas
in \cite{BY}, attention is restricted to {\em continuous}
 bounded test functions.
By taking test functions which are
indicator functions of arbitrary Borel sets $A_1,\ldots,A_m$ in $\R^d$,
%(or linear combinations thereof), 
we see from Theorem \ref{thm2} that under Assumption A1,
 the joint  distribution
of $(\la^{-1/2} \bar{\mu}_\la^\xi(A_i), 1 \leq i \leq m)$ converges 
to a multivariate normal with covariances given by
$\int_{A_i \cap A_j} V^\xi(\kappa(\bx)) \ka(\bx)d\bx$,
and likewise for $\overline{\nu}_{\la(n),n}^\xi$ by Theorem \ref{thm3}.
This desirable conclusion is not achieved from the results
of \cite{BY}, because indicator functions of Borel sets
are not continuous. When our assumption A1  fails, for
the central limit theorems we 
restrict attention to
%jneed (Assumption A3) to assume
 almost everywhere continuous
 test functions,
 which means we can still obtain the above conclusion
provided the sets $A_i$ have Lebesgue-null boundary.
 % is still weaker than the assumption of continuous
%test functions in \cite{BY}.

The de-Poissonization argument in \cite{BY} requires finiteness of
what might be called the radius of {\em external} stabilization;
see Definition 2.3 of \cite{BY}. 
Loosely speaking, an inserted point at $\bx$ is not affected by
 and {\em does not affect} points at a distance beyond
the radius of external stabilization; in contrast  
an inserted point at $\bx$ is unaffected by points
at a distance beyond the radius of stabilization,
but might affect other points beyond that distance.
Our approach does not require external stabilization, which
brings some examples within the scope of our results that
do not appear to be covered by the results of \cite{BY}.
See the example of germ-grain models, considered in Section \ref{secgerm}.

In the non-translation-invariant case, we require $\xi$ to
have almost everywhere continuous total measure,
 whereas in \cite{BY} the functional
$\xi$ is required to be in a class 
SV(4/3)
 of `slowly varying' functionals. The almost everywhere 
continuity condition on $\xi$
is usually easier to check. 
%In the translation invariant
%case, simple modifications of the proofs
%of our results (which we omit) show that
%they hold without requiring
%the almost eveywhere continuity condition.

We assume the underlying density function $\ka$
is almost everywhere continuous, and for Theorems
\ref{thm2} and \ref{thm3} but {\em not} for Theorem \ref{thm1},
we assume it has compact support. In contrast,
in \cite{BY} 
it is assumed that $\ka$ has compact convex support
and is continuous on its support (see the remarks just before Lemma
4.2 of \cite{BY}). 

Our  moments condition (\ref{mom2})
is simpler than the corresponding condition 
in \cite{BY} (eqn (2.2) of \cite{BY}). 
Using A$7$ and A$7'$ in Theorems \ref{thm2} and \ref{thm3},
 we obtain Gaussian limits for
random fields under polynomial stabilization of sufficiently
high order; the  corresponding results in \cite{BY} need 
 exponential stabilization.  

%Finally, our proof of
% Theorems \ref{thm2} and \ref{thm3} is  
%different from that of the corresponding results in \cite{BY}.
%and possibly somewhat simpler. 
%
%Unlike \cite{BY},
We spell out the statement and proof of
 Theorems \ref{thm2} and \ref{thm3}   
for the setting of {\em marked} point processes
setting (i.e. point processes in $\R^d \times \MM$ rather than
in $\R^d$), whereas the proofs in earlier works
\cite{BY,PY1} are given for the setting of unmarked point process
(i.e., point processes in $\R^d$). The marked point process setting includes 
many interesting examples such as germ-grain models
%, the lilypad model,
and  on-line packing, and
generalizes the unmarked point process setting
because we can always take $\MM$ to have a single element
and then identify $\R^d \times \MM$ with $\R^d$, to recover
results for unmarked point processes from the general
results for marked point processes.  

Other papers concerned with central limit theorems for random measures
include Henrich and Molchanov \cite{HM} and Penrose \cite{Pe04}.
% Both of these restrict
%their attention to homogeneous point processes. 
The setup of \cite{HM} is somewhat different from ours; 
the emphasis there is on measures associated with germ-grain models and
 the method for defining the measures from
the marked point sets (eqns (3.7) and (3.8) of \cite{HM}) is
more prescriptive than that used here. 
In \cite{HM} the underlying point processes are
 taken to be stationary point processes satisfying a mixing condition
and no  notion of stabilization is used,
whereas we restrict attention
to Poisson or binomial point processes but do not require any
spatial homogeneity.

The setup in \cite{Pe04} is closer to that used here
 (although the proof of central limit theorems is different)
but has the following notable differences.
The point processes considered  in \cite{Pe04}
 are assumed to have 
constant intensity on their support. The notion of  stabilization used in
\cite{Pe04} is a form of {\em external} stabilization.
For the multivariate central limit theorems  in \cite{Pe04}
to be applicable,
the radius of external stabilization needs to be almost surely
finite but, unlike in the present work, no
bounds on the tail of this radius of stabilization are required.
The test functions in \cite{Pe04} lie in a subclass
of $\tB(\R^d)$, not $B(\R^d)$.
The description of the limiting variances in \cite{Pe04} is
different from that given here.

In most examples, the sets $\Gamma_\la$ are all the same as $\Gamma$.
However, there are cases where
moments conditions such as \eq{momlnb} and \eq{mom4} 
hold for a sequence of sets $\Gamma_\la$ but would not hold if
we were to take $\Gamma_\la =\Gamma$ for all $\la$. See, e.g.
\cite{PW}.

%{\bf Non-marked point processes.} If $\MM$ has a single element,
%then we can identify $\R^d \times \MM$ with $\R^d$ and
%obtain a theory for stabilizing functionals on subsets of
%$\R^d$. 

%Fifth, our result  requires a condition of homogeneous
%stabilization which is not included as a condition in
%\cite{BY}. However, all functionals known to the
%author which satisfy exponential stabilization also
%satisfy homgeneous stabilization. Also, it can be shown (?)
%that exponential stabilization with respect to $(\kappa,\Gamma)$
%implies homoegeneous stabilization in the case where $\Gamma = \R^d$.
%Finally, it may be that in fact the proof of the results in
%\cite{BY} actually requires homogeneous stabilization as an
%implicit condition. 

Last but not least, we note that our results carry over
to the case where
 $\xi(\bx,t;\X,\cdot)$ 
is a {\em signed} measure with finite total variation.
The conditions for the theorems remain unchanged if we take signed
measures, except that if $\xi$ is a signed measure,
the moments conditions \eq{mom}, \eq{mom2}, \eq{momlnb} and
\eq{mom4}  need to hold for both the positive and the negative
part of $\xi$. The proofs need only minor modifications
to take signed measures into account.

\section{Weak convergence:
the objective method}
\lbl{secobjective}
\allco
In the proof of Theorems \ref{thlln1} \ref{thm1},
and \ref{thm3}, most of the work required is to
prove convergence of first and second moments. 
A 
%first
 step in this direction is to obtain certain
weak convergence results, namely Lemmas \ref{couplem0a},
\ref{couplem0c}, \ref{lemweak1}, \ref{lemweak2} and \ref{lemweak3} below.
It is noteworthy that in all of these lemmas, the
stabilization conditions used always refer to {\em homogeneous}
Poisson processes on $\R^d$; the notion of exponential stabilitization
with respect to a non-homogeneous point process is not used 
until later on.

To prove these lemmas,
we shall use a version of what is sometimes
called the  `objective method' \cite{AS,Stebk},
whereby convergence in distribution (denoted $\tod$)
for a functional defined on  a
sequence of finite probabilistic objects
(in this case, re-scaled marked point processes), is
established by showing that these probabilistic objects
themselves converge in distribution to an infinite probabilistic
object (in this case, a homogeneous marked Poisson process),
and that the functional of interest is continuous. 
We can then use the {\em Continuous Mapping Theorem}
(\cite{Bill}, Chapter 1, Theorem 5.1),
which says that if $h$ is a mapping from a metric space
$E$ to another metric space $E'$, and $X_n $ are $E$-valued
random variables converging in distribution to $X$ which
lies almost surely at a continuity point of $h$, then
$h(X_n)$ converges in distribution to $h(X)$.

A {\em point process} in $\R^d\times \MM$ is an $\LL$-valued
random variable, where $\LL$ denotes the space of
locally finite subsets of $\R^d \times \MM$.
We use the following metric on $\LL$:
\bea
\dist(\A , \A') =\left(  \max \{K \in \N: \A  \cap [ B_K \times \MM ]
 = \A' \cap [B_K \times \MM] \} \right)^{-1}.
\lbl{ppmet}
\eea 
With this metric, $\LL$ is a metric space which is complete but
not separable. In the unmarked case where $\MM$ has a single element,
 our choice
 of metric is {\em not} the same as the 
metric used in Section 5.3 of \cite{Stebk}.
Indeed, for one-point unmarked sets our metric
generates the discrete topology rather than the Euclidean topology.

To prove the weak convergence of point processes, we use
 a refinement of a coupling method
used in \cite{PY4}.
% to obtain a less general version
%of the law of large numbers in Theorem \ref{thlln1}.
In particular, we shall use a device which we
here call the {\em pivoted coupling}.

Let $\bx_0 \in \R^d$ and let $\la >0$.
% and $\mu >0$.  
Let $\H^+$ denote a homogeneous Poisson process of unit intensity 
in $\R^d \times \MM \times [0,\infty)$. Let $\Po'_\la$
 denote the image of the restriction of $\H^+$ to the set
 $\{(\bx,t,s)\in \R^d \times \MM \times [0,\infty): s \leq \la \ka(\bx)\}$,
under the mapping $(\bx,t,s)\mapsto (\bx,t)$. 
For $a >0$, let $\H'_{a}$,
 denote the image of the restriction of $\H^+$ to the set
 $\{(\bx,t,s)\in \R^d \times \MM \times [0,\infty): s \leq  \la a )\}$,
under the mapping
$$
(\bx,t,s)\mapsto (\la^{1/d}(\bx-\bx_0),s).
$$
Then by the Mapping Theorem \cite{Kin}, 
 $\Po'_\la$ has the same distribution as $\Po_\la$
while $\H'_{a}$ has the same distribution as $\H_{a}$.
We shall refer to $\Po'_\la$ and $\H'_a$
as {\em coupled realizations of $\Po_\la$
and $\H_{a}$ with pivot point $\bx_0$}. 

\begin{lemm}
\lbl{pivlem}
Suppose $\bx \in \R^d$ is a continuity point of $\kappa$, 
%suppose $(a(\la), \la >0)$ is a real-valued function 
%which tends to $\kappa(\bx)$ as $\la \to \infty$,
 and suppose $(\by(\la), \la >0) $ is an $\R^d$-valued function which
tends to $\bx$ as $\la \to \infty$. 
If $\Po'_\la$ and $\H'_{\ka(\bx)}$ are the coupled realizations
%of $\Po_\la$ and $\H_{a(\la)}$ with pivot point $\by(\la)$,
of $\Po_\la$ and $\H_{\ka(\bx)}$ with pivot point $\by(\la)$,
then for any $K\in (0,\infty)$,
 as $\la \to \infty$,
\bea
P[ \la^{1/d}(-\by(\la)+\Po'_{\la}) \cap ( B_K \times \MM)=  
\H'_{\ka(\bx)} \cap (B_K\times \MM) ]  \to 1.
%  ~~{\rm as }~ \la \to \infty.
\lbl{0426}
\eea  
\end{lemm}
{\em Proof.} The number of points of the point set 
$$
(\la^{1/d}(-\by(\la) + \Po'_\la) \triangle
 \H'_{\ka(\bx)} ) \cap (B_K \times \MM)
 %\tH_{a(\la)} ) \cap (B_K \times \MM)
$$
equals  the number of ponts $(\bX,T,S)$ of $\H^+$ with 
$\bX \in B_{\la^{-1/d}K}(\y(\la))$
and with either $\la \ka(\bx)< S \leq \la \ka(\bX)$ or 
 $\la \ka(\bX)< S \leq \la \ka(\bx)$. This is Poisson distributed with 
mean
$$
\la \int_{B_{\la^{-1/d} K} (\ka(\bx))} |\ka(\bz)- \ka(\bx)|d\bz
$$
which tends to zero by the assumed continuity of $\ka(\cdot)$ at $\bx$.
% and the assumption that $a(\la)\to \ka(\bx_0)$.
 $\qed$ \\

In second moment computations, we are interested
weak convergence, 
not only of a point process in $\LL$,
but also of a {\em pair} of point processes in $\LL \times \LL$.
The limiting object in this case will be a pair of independent homogeneous
Poisson processes, and we need notation for this. To this end,
for $a>0$ and $b>0$ let 
$\tH_{b}$ denote  a homogeneous Poisson process
on $\R^d \times \MM$, independent of $\H_a$.
% has the distribution 
%of $\H_{\ka(\by)}$ and is independent of $\H_{\ka(\bx)}$.

%In the statement of results in Section \ref{secresults},
%it was not necessary to assume that the  
Our first weak convergence lemma is concerned 
with the point processes $\X_m$. In this result, 
we assume these point processes are coupled togeteher in
the natural way
(this was not needed for the statement of results
 in Section \ref{secresults}).
To do this,
%In proving the law of large numbers (Theorem \ref{thlln1}), 
%and in de-Poissonizing the central limit theorem,
%we shall use the following lemma, which
%is proved by using the pivoted coupling.
%In this lemma, 
 we let $(\bX_1,T_1),(\bX_2,T_2), \ldots$ denote a sequence of independent
identically distributed random elements of $\R^d \times \MM$ with
distribution $\kappa \times \mu_\MM$, and assume the point   
processes $\X_m, m \geq 1$ are given by
%coupled by setting 
\bea
\X_m = \{(\bX_1,T_1), (\bX_2,T_2),\ldots, (\bX_m,T_m)\}
\lbl{050404}
\eea
 for each $m$.
In proving the lemma, we use the notation $\eqd$ for equality of
distribution,
%or for equality of joint distribution
%(so e.g. $(X,Y)\eqd (X',Y')$ means $X,Y$ have the
%same joint distribution 
%as $X',Y'$).

\begin{lemm}
\lbl{pplemweak}
Suppose $(\bx, \by) \in \R^d \times \R^d$
with $\ka(\bx)>0$, $\ka(\by)>0$, with $\ka(\cdot)$ continuous at $\bx$ and
at $\by$, and with $\bx \neq \by$.
Let $(\la(k),\ell(k),m(k))_{k \in \N}$
be a $((0,\infty) \times \N \times \N)$-valued sequence
satisfying $\la(k)\to \infty$, and $\ell(k)/\la(k) \to 1$
and $m(k)/\la(k)\to 1$ as $k \to \infty$. Then as $k \to \infty$,
\bea
%(\la(k)^{1/d}(-\bx + \X_{\ell_k}),
%\la(k)^{1/d}(-\bx + \X_{m_k}),
%\la(k)^{1/d}(-\by + \X_{m_k}), 
%%\nonumber \\
%(\la(k)^{1/d}(-\bx + \X_{\ell_k}^{\by,T'}),
%\la(k)^{1/d}(-\bx + \X_{m_k}^{\by,T'}),
%\la(k)^{1/d}(-\by + \X_{m_k}^{\bx,T}) )
%\nonumber \\
(\la(k)^{1/d}(-\bx + \X_{\ell(k)}),
\la(k)^{1/d}(-\bx + \X_{m(k)}),
\la(k)^{1/d}(-\by + \X_{m(k)}), 
\nonumber \\
\la(k)^{1/d}(-\bx + \X_{\ell(k)}^{\by,T'}),
\la(k)^{1/d}(-\bx + \X_{m(k)}^{\by,T'}),
\la(k)^{1/d}(-\by + \X_{m(k)}^{\bx,T}) )
\nonumber \\
\tod ( \H_{\ka(\bx)},\H_{\ka(\bx)},\tH_{\ka(\by)},
 \H_{\ka(\bx)},\H_{\ka(\bx)},\tH_{\ka(\by)}).
\lbl{050401a}
\eea
\end{lemm}
{\em Proof.}
In this proof  we ease notation by suppressing
 all mention of the marks in the notation,
and write for example $\bX$ for $(\bX,T)$.
Moreover, we suppress mention of the parameter $k$ and
write simply $\la$ for $\la(k)$, $\ell$ for $\ell(k)$,
and $m$ for $m(k)$.

We use the following coupling.
Suppose we are given $\la$.  On a suitable probability
space, let $\Po$ and $\tPo$ be independent Poisson processes
on $\R^d$ with intensity function $\la \ka(\cdot)$ with each point
carrying an independent $\MM$-valued mark with distribution
$\mu_\MM$.
 Let $\bY_1,\bY_2,\ldots$ be independent
random $d$-vectors with distribution $\ka$ in $\R^d$  carrying marks
with distribution $\mu_\MM$,
 independent of $\Po$ and
$\tPo$. 

Let $\Po'$ be the point process consisting of those
(marked) points
of $\Po$ which lie closer to $\bx$ than to $\by$ (in the Euclidean norm),
together with
 those (marked) points
of $\tPo$ which lie closer to $\by$ than to $\bx$.  Clearly $\Po' $
is a (marked) Poisson process of intensity  $\la \ka(\cdot)$ on $\R^d$.

%Let $\tilde{\H}_{\ka(\bx)}$ be a (marked) homogeneous Poisson point process
% in $\R^d$, of intensity $\ka(\bx)$. Similarly 
%let $\tilde{\H}'_{\ka(\by)}$ be a  (marked)
% homogeneous Poisson process in $\R^d$ of intensity $\ka(\by)$.
Let $\H'_{\ka(\bx)}$ and $\tH'_{\ka(\by)}$ 
 be  (marked) homogeneous Poisson point processes
 in $\R^d$, of intensity $\ka(\bx)$ and $\ka(\by)$ respectively.
Assume  $\H'_{\ka(\bx)}$ and
 $\tH'_{\ka(\by)}$ are independent of each other
and of $(\bY_1,\bY_2, \bY_3,\bY_3,\ldots)$. Assume also
that $\H'_{\ka(\bx)}$ is coupled to $\Po$ by the pivoted coupling
 with pivot at $\bx$ (see Lemma \ref{pivlem}), and
$\tH'_{\ka(\by)}$ is coupled to $\tPo$ by the pivoted coupling with
pivot at $\by$.

Let $N$ denote the number of points of $\Po'$
(a Poisson variable with mean $\la$). Choose an ordering on the points
of $\Po'$, uniformly at random from all $N!$
possible such orderings.
 Use this ordering to list the points of
$\Po'$ as $\bW_1,\bW_2,\ldots,\bW_N$. Also, set $\bW_{N+1} =
\bY_1, \bW_{N+2} = \bY_2,\bW_{N+3} = \bY_3$ and so on. Set
$$
\X'_\ell := \{ \bW_1,\ldots,\bW_\ell\}, ~~~
\X'_m := \{ \bW_1,\ldots,\bW_m\}.
$$
Clearly $(\X'_\ell,\X'_m) \eqd (\X_\ell,\X_m)$,
and $(\H'_{\ka(\bx)},\tH'_{\ka(\by)}) \eqd (\H_{\ka(\bx)},\tH_{\ka(\by)})$

Let $K\in \N$, and let $\delta >0$.
Let $\theta$ denote the volume of the unit ball in $d$ dimensions.
%Let $\eps >0$, and
%set $\delta := \eps /(20 \kamax \theta K^d)$.
Define the events
\bean
E:= 
\{ \X'_m  \cap B_{\la^{-1/d} K}(\bx) 
=
 \Po' \cap B_{\la^{-1/d} K}(\bx) \};
\\
F:= \{
(\la^{1/d}(-\bx + \Po))\cap B_K = \H'_{\ka(\bx)} %\H_{\ka(\bX)} 
\cap B_K \}. 
\eean
Event $E$ occurs unless, either,  one or more of the
$(N-m)^+$ ``discarded'' points of $\Po'$, or, one or more of the
$(m-N)^+$ ``added'' points of $\{\bY_1,\bY_2,\ldots\}$ lies in
$B_{\la^{-1/d} K}(\bx)$.  
For each added or discarded point, the probability of
lying in $B_{\la^{-1/d} K}(\bx)$ is at most $\theta \kamax K^d/\la$.
Thus, for $k$  large enough so that  $|m-\la|\leq \delta \la$, we have
\bean
P[E^c] \leq P[|N-\la| > \delta \la ] +
 (2\delta \la)\theta \|\ka\|_\infty K^d/\la
\eean
which is less than $
 3\delta \theta \|\ka\|_\infty K^d$ for large enough $k$.
Hence, $P[E^c] \to 0$ as $k \to \infty$.
Moroever,  by  \eq{0426} we also have $P[F^c] \to 0$ as $k \to \infty$.

Assuming $\la$ is so large that
$ | \bx -\by | > 2 \la^{-1/d} K, $
if $E \cap F$ occurs then
\bean
\H'_{\ka(\bx)} \cap B_K
= 
( \la^{1/d}(-\bx+ \Po) ) \cap B_K
\\
=
 \la^{1/d}((-\bx+ \Po)\cap B_{\la^{-1/d} K})
=
\la^{1/d}(-\bx+ ( \Po \cap B_{\la^{-1/d} K}(\bx)))
\\
=  \la^{1/d}(-\bx+ (\X'_m  \cap B_{\la^{-1/d} K}(\bx)))
\\
=  \la^{1/d}((-\bx+ \X'_m)  \cap B_{\la^{-1/d} K})
%\\
=  (\la^{1/d}(-\bx+ \X'_m) ) \cap B_{K}
\eean
so that $\dist(\H'_{\ka(\bx)}, \la^{1/d}(-\bx+ \X'_m) ) \leq 1/K$.
Hence, for any $K$ we have
 $$
P[\dist(\H'_{\ka(\bx)}, \la^{1/d}(-\bx+ \X'_m) ) > 1/K ] \to 0,
$$ 
Similarly, we have
 \bean
\max \{ P[\dist(\H'_{\ka(\bx)}, \la^{1/d}(-\bx+ \X'_\ell) ) > 1/K ], 
%\to 0,
%~~~
%\\ 
 P[\dist(\tH'_{\ka(\by)}, \la^{1/d}(-\by+ \X'_m) ) > 1/K ], 
%\to 0,
\\ 
 P[\dist(\H'_{\ka(\bx)}, \la^{1/d}(-\bx+ (\X'_\ell)^\by) ) > 1/K ], 
%\to 0,
%~~~
%\\ 
 P[\dist(\H'_{\ka(\bx)}, \la^{1/d}(-\bx+ (\X'_m)^\by) ) > 1/K ] ,
%\to 0,
\\ 
 P[\dist(\tH'_{\ka(\by)}, \la^{1/d}(-\by+ (\X'_m)^\bx) ) > 1/K ] 
\} \to 0.
\eean
Combining these, we have the required convergence in 
distribution. $\qed$ \\

For subsequent results, it is useful to define
the region
\bea
\Gamma_0 := \{\bx \in \inter(\Gamma): \ka(\bx) >0, \kappa(\cdot) 
\mbox{ continuous at } \bx \}
\lbl{gam1def}
\eea

\begin{lemm}
\lbl{couplem0a}
Suppose $(\bx, \by) \in \Gamma_0 \times \Gamma_0$, with
%$\ka(\bx)>0$, $\ka(\by)>0$ and
 $\bx \neq \by$.
Suppose also that
$R(\bx,T;\H_{\ka(\bx)})$ 
and $R(\by,T;\H_{\ka(\by)})$ 
are almost surely finite.
Suppose $(\la(m))_{m \geq 1}$ is a $(0,\infty)\times \N$-valued
sequence with $\la(m ) /m \to 1$  as $m  \to \infty$.
Then for Borel $A \subseteq \R^d$,
as $m \to \infty$ we have
\bea
 \xi_{\la(m)}(\bx,T;\X_m,\bx + \la^{-1/d} A ) 
\tod  \xi_\infty^{\bx,T} (\H_{\ka(\bx)},A),
\lbl{050401c}
\lbl{weak00}
\eea
%\bea
% \xi_{\la(m)}(\bx,T;\X_m,\R^d) 
%\tod
%  \xi_\infty^{\bx,T}(\H_{\ka(\bx)},\R^d),
%\lbl{weak00}
%\eea
and
\bea
%( \xi_{\la(m)}(\bx,T;\X_m,\bx + \la^{-1/d} A ), 
%\xi_{\la(m)}(\by,T';\X_m,\by + \la^{-1/d} A' ) )
( \xi_{\la(m)}(\bx,T;\X_m^{\by,T'},\bx+ \la^{-1/d} A), 
\xi_{\la(m)}(\by,T';\X_m^{\bx,T},\by+ \la^{-1/d} A) )
\nonumber
\\
\tod
% ( \xi_\infty^{\bx,T}(\H_{\ka(\bx)},A),
%\xi_\infty^{\by,T'}(\tH_{\ka(\by)},A') ) .
 ( \xi_\infty^{\bx,T}(\H_{\ka(\bx)},A),
\xi_\infty^{\by,T'}(\tH_{\ka(\by)},A) ) .
\lbl{050401b}
\eea
\end{lemm}
{\em Proof.}  
Given $A$,
 define the mapping $h_{A,\bx}:\MM \times \LL \to \R$ and
the mapping $h_A^2: (\MM \times \LL \times \MM \times \LL) \to  \R^2$ by
\bea
h_{A,\bx}(t,\X) =  \xi (\bx,t;\bx + \X, \bx + A) ;
\lbl{050404a}
\\
h_A^2(t,\X,t',\X') = (h_{A,\bx}(t,\X),h_{A,\by}(t',\X')).
\nonumber
\eea
%We assert that 
Since $R(\bx,T;\H_{\ka(\bx)}) <\infty$ a.s., 
the pair 
 $(T,\H_{\ka(\bx)})$ lies a.s.
at a continuity point of $h_{A,\bx}$, where the topology on $\MM \times \LL$
is the product of the discrete topology on $\MM$ and the topology
induced by our metric $\dist$ on $\LL$, defined at \eq{ppmet}.
Similarly, 
 $(T,\H_{\ka(\bx)}, T',\tH_{\ka(\by)})$ lies a.s.
at a continuity point of $h_A^2$.
We have by definition of $\xi_\la$ that
\bean
 \xi_{\la}(\bx,T;\X_m,\bx + \la^{-1/d} A ) 
%\\
= h_{A,\bx}(T,\la^{1/d}(-\bx+ \X_m) ) ;
\\
( \xi_{\la}(\bx,T;\X_m^{\by,T'},\bx+ \la^{-1/d} A),
\xi_{\la}(\by,T';\X_m^{\bx,T},\by+ \la^{-1/d} A ) )
\\
= h_A^2(T,\la^{1/d}(-\bx+ \X_m^{\by,T'}), 
T',\la^{1/d}(-\by+ \X_m^{\bx,T})).
\eean
By Lemma \ref{pplemweak},
we have
$
(T,\la^{1/d}(-\bx+ \X_m))
\tod  (T, \H_{\ka(\bx)})
$
so that \eq{050401c}  follows by  the Continuous Mapping Theorem.
%(\cite{Bill}, Chapter 1, Theorem 5.1). 
Also,
by Lemma \ref{pplemweak},
\bean
(T,\la^{1/d}(-\bx+ \X_m^{\by,T'}),T',\la^{1/d}(-\by+ \X_m^{\bx,T'}))
\tod  (T,\H_{\ka(\bx)},T', \tH_{\ka(\by)} )
\eean
so that \eq{050401b}  also follows by  the Continuous Mapping Theorem.
%
%The argument for \eq{050401c} is similar.
%Consider the mapping $h:\MM^2 \times \LL^2 \to \R^2$ given by
%\bean
%(t,t',\X,\X') \mapsto (\xi (\bx,t;\bx + \X, \R^d ),
%\xi (\by,t';\by + \X', \R^d) ),
%\\
%\eean 
%%We assert that with probability 1, the quadruple
%% $(T,T',\H_{\ka(\bx)},\H_{\ka(\by)})$ lies
%%at a continuity point of $h$, where the topology on $\MM^2 \times \LL^2$
%is the product of the discrete topology on both copies of $\MM$ 
%and the topology
%induced by our metric $\dist$ on both copies of $\LL$.  This 
%assertion follows from the 
%stabilization assumption that 
%$R(\bx,T;\H_{\ka(\bx)})$ 
%and $R(\by,T';\H_{\ka(\by)})$ 
%are almost surely finite.
%We have by definition of $\xi_\la$ that
%\bean
%( \xi_{\la}(\bx,T;\X_m,\R^d),
%\xi_{\la}(\by,T';\X_m,\R^d ) )
%= h(T,T',\la^{1/d}(-\bx+ \X_m; \la^{1/d}(-\by+ \X'_m))
%\eean
%and by Lemma \ref{pplemweak},
%we have
%(\cite{Bill}, Chapter 1, Theorem 5.1).
$\qed$ \\ 

The next lemma is one of several with the purpose
of comparing (with reference to a test function $f$)
 the measure $\xi_\la(\bx;\X,\cdot)$ to
the corresponding point measure $\xi_\la^*(\bx;\X,\cdot)$ so
as to derive results when Assumption A2 or A3, rather than A1,
is the case. In proving such results, we repeatedly use the notation
 \bea
\phi_\eps(\bx):= \sup\{|f(\by)-f(\bx)|: \by \in B_\eps (\bx)\}, ~~~
 \mbox{ for } \eps >0,
\lbl{phidef}
\eea
and for $f\in B(\R^d)$ we
  write $\|f\|_\infty$ for $\sup\{|f(\bx)|:\bx \in \R^d\}$.
Also, recall (see e.g. \cite{Pe}, \cite{Rudin}) that
 $\bx \in \R^d$ is 
a {\em Lebesgue point} of $f$ if $\eps^{-d}\int_{B_{\eps}(\bx)} |f(\by)-f(\bx)|
d\by$ tends to zero as $\eps \downarrow 0$, and that the
Lebesgue Density Theorem
 tells us that almost every $\bx \in \R^d$
is a Lebesgue point of $f$.

\begin{lemm}
\lbl{couplem0b}
Let $\bx \in \Gamma_0$, and
% with
%$\ka(\bx)>0$, and suppose
suppose
 that
$R(\bx,T;\H_{\ka(\bx)}) < \infty$ 
 almost surely.
Let
 $\by \in \R^d$ with $\by \neq \bx$.
Suppose that $f \in B(\R^d)$, and suppose either
 that $f$ is continuous at $\bx$,
 or that Assumption A$2$  holds and $\bx$ is a Lebesgue point of $f$.
Suppose $(\la(m))_{m \geq 1}$ is a $(0,\infty)\times \N$-valued
sequence with $\la(m ) /m \to 1$  as $m  \to \infty$.
Then as $m \to \infty$,
\bea
\langle f,
  \xi_{\la(m)}(\bx,T;\X_m)
-  
  \xi^*_{\la(m)}(\bx,T;\X_m) \rangle \toP 0
\lbl{050401d}
\eea
and
\bea
\langle f,
  \xi_{\la(m)}(\bx,T;\X_m^{\by,T'})
-  
  \xi^*_{\la(m)}(\bx,T;\X_m^{\by,T'}) \rangle \toP 0
\lbl{050401j}
\eea
\end{lemm}
{\em Proof.}
In this proof, 
we suppress the mark in the notation, writing $\bx$ for $(\bx,T)$,
and $\xi(\bx;\X,A)$ for 
 $\xi(\bx,T;\X,A)$,
 and so on.
Also, we write $\la$ for $\la(m)$.
The left side of \eq{050401d} is equal to
\bea
  \int_{\R^d} ( f(\bz) - f(\bx) )
\xi_\la(\bx; \X_m,d\bz ).
\lbl{50209a2}
\eea
Given $K>0$, we split the region of integration
in \eq{50209a2} into 
the complementary regions
 $B_{\la^{-1/d} K}(\bx)$ 
and  
 $\R^d \setminus B_{\la^{-1/d} K}(\bx)$. 
Consider the latter region first. 
%Given $K>0$, 
By
\eq{050401c} we have
% change for signed meas
\bean
 \left|
 \int_{\R^d \setminus B_{\la^{-1/d} K}(\bx)} ( f(\bz) - f(\bx) )
\xi_\la(\bx; \X_m,d\bz ) \right|
\leq 
2 \|f\|_\infty 
 \xi_\la(\bx;\X_m,\R^d \setminus B_{\la^{-1/d} K} (\bx)  ) 
\nonumber \\
\tod 
2 \|f\|_\infty 
 \xi_\infty^{\bx}(\H_{\ka(\bx)},\R^d \setminus B_{K}  ),
\eean 
where the limit is almost surely finite and converges in probability
to zero as $K \to \infty$. Hence for  $\eps >0$, we have
\bea
\lim_{K \to \infty} \limsup_{m \to \infty} P\left[
 \left|
 \int_{\R^d \setminus B_{\la^{-1/d} K}(\bx)} ( f(\bz) - f(\bx) )
\xi_\la(\bx; \X_m,d\bz ) \right| > \eps \right] =0.
\lbl{050401e}
\eea

Turning to the integral over
%in \eq{50209a2} to 
 $B_{\la^{-1/d} K}(\bx)$, we consider separately
the case where $f$ is continuous at $\bx$, and the case
where A2 holds and $\bx$ is a Lebesgue point of $f$. 
To deal with the first of these cases, 
observe that
%given $K>0$, 
% it is the case that
% change if signed
\bea
 \left| \int_{ B_{\la^{-1/d} K}(\bx)} ( f(\bz) - f(\bx) )
\xi_\la(\bx; \X_m,d\bz ) \right|
%\nonumber \\
\leq
\phi_{\la^{-1/d}K}(\bx)
 %|\xi_\la|(\bx;\X_m,\R^d)  |
 \xi_\la(\bx;\X_m,\R^d)  
\lbl{050401f}
\eea
and if $f$ is continuous at $\bx$, then
$\phi_{\la^{-1/d}K}(\bx) \to 0$, while 
 %$|\xi_\la|(\bx,\X_m,\R^d)  $ if signed 
 $\xi_\la(\bx,\X_m,\R^d)  $ 
converges in distribution
to the finite random variable
 %$|\xi_\infty^{\bx,T}|(\H_{\ka(\bx)},\R^d)  $  if signed
 $\xi_\infty^{\bx,T}(\H_{\ka(\bx)},\R^d)  $ 
by 
\eq{050401c}, and hence the right hand side of
\eq{050401f} tends to zero in probability as $m \to \infty$.
Combined with \eq{050401e}, this gives us 
\eq{050401d} in the case where $f$ is continuous at $\bx$. 

Under Assumption A2, 
for Borel $A \subseteq \R^d$,
 the change of variables $\bz = \bx +
\la^{-1/d}(\by-\bx)$ yields 
\bean
%\int_A \xi_\la(\bx;\X,d\by) = 
\xi_\la(\bx;\X,A)  = 
%\1_{\Gamma_\la}(\bx) 
%\xi 
%(\bx; \bx+\la^{1/d}(-\bx +\X),\bx + \la^{1/d}(-\bx + A) )
%\nonumber \\
%= 
%\1_{\Gamma_\la}(\bx)
 \int_{\bx + \la^{1/d} (-\bx +A)} \xi'
(\bx; \bx+\la^{1/d}(-\bx +\X),\by ) d\by
\nonumber  \\
= \la
% \1_{\Gamma_\la}(\bx) 
\int_{A} \xi'
(\bx; \bx+\la^{1/d}(-\bx +\X),\bx + \la^{1/d}(\bz -\bx) ) d\bz .
%\lbl{50210a}
\eean
Hence, under A2,
\bean
\left| \int_{ B_{\la^{-1/d} K}(\bx)} ( f(\bz) - f(\bx) )
\xi_\la(\bx; \X_m,d\bz ) \right|
\nonumber \\
= \la 
 \left| \int_{ B_{\la^{-1/d} K}(\bx)} ( f(\bz) - f(\bx) )
\xi'
(\bx; \bx+\la^{1/d}(-\bx +\X_m),\bx + \la^{1/d}(\bz -\bx) ) d\bz 
\right|
\nonumber \\
\leq K_0 \la 
  \int_{ B_{\la^{-1/d} K}(\bx)} | f(\bz) - f(\bx) | d\bz 
\eean
and if additionally $\bx$ is a Lebesgue point of $f$ then
this tends to zero.  
Combined with \eq{050401e}, this gives us 
\eq{050401d} in the case where A2 holds and $\bx$ is a Lebesgue point
of $f$.

The proof of \eq{050401j} is similar; we use \eq{050401b}
instead of \eq{050401c}. $\qed$ \\

By combining Lemmas \ref{couplem0a} and \ref{couplem0b},
we obtain the following, which is the main ingredient
in our proof of the Law of Large Numbers in Theorem \ref{thlln1}. 
\begin{lemm}
\lbl{couplem0c}
Suppose $(\bx, \by ) \in \Gamma_0 \times \Gamma_0$, with
%$\ka(\bx)>0$, $\ka(\by)>0$ and
 $\bx \neq \by$.
Suppose also that
$R(\bx,T;\H_{\ka(\bx)})<\infty$ 
and $R(\by,T;\H_{\ka(\by)})<\infty$, 
almost surely.
Let $f \in B(\R^d)$ and suppose that either A1 holds, or
A2 holds and $\bx$ is a Lebesgue point of $f$, or $\bx$
is a continuity point of $f$.
Suppose $(\la(m))_{m \geq 1}$ is a $(0,\infty)\times \N$-valued
sequence with $\la(m ) /m \to 1$  as $m  \to \infty$.
Then as $m \to \infty$,
\bea
\langle f, \xi_{\la(m)}(\bx,T;\X_m) \rangle 
\tod
f(\bx)  
  \xi_\infty^{\bx,T}(\H_{\ka(\bx)},\R^d)
\lbl{050401g}
\eea
and
\bea
\langle f, \xi_{\la(m)}(\bx,T;\X_m^{\by,T'}) \rangle 
\langle f, \xi_{\la(m)}(\by,T';\X_m^{\bx,T}) \rangle
\nonumber
\\
\tod
f(\bx) f(\by) 
  \xi_\infty^{\bx,T}(\H_{\ka(\bx)},\R^d)
  \xi_\infty^{\by,T'}(\tH_{\ka(\by)},\R^d).
\lbl{050401h}
\eea
\end{lemm}
{\em Proof.} Note first that
$$
%\langle f, \xi^*_{\la(m)}(\bx,t;\X_m^{\by,t'}) \rangle
\langle f, \xi^*_{\la(m)}(\bx,t;\X) \rangle
=  f(\bx)  
 %\xi_{\la(m)}(\bx,t;\X_m^{\by,T'},\R^d).
 \xi_{\la(m)}(\bx,t;\X,\R^d).
$$
Hence,
in the case where A1 holds (i.e., $\xi =\xi^*$),  
 \eq{050401g} is immediate from the case $A=\R^d$ of
\eq{weak00},
and similarly
\eq{050401h} is immediate from \eq{050401b}.

In the other two cases described, we have \eq{050401d} by
 Lemma \ref{couplem0b}. Combining this with
 \eq{050401g} 
for the case with $\xi=\xi^*$, we see by Slutsky's theorem
(see, e.g., \cite{Pe})
that \eq{050401g} still holds
in the other two cases. Similarly, since \eq{050401h}
holds  when $\xi= \xi^*$, by \eq{050401j} and Slutsky's
theorem we can obtain \eq{050401h} in the other cases too.
$\qed$  \\

The next two lemmas are key ingredients in proving Theorem
\ref{thm1} on convergence
of second moments.

\begin{lemm}
\lbl{lemweak1}
 Suppose $\bx \in \Gamma_0$.
%$\bx\in \inter(\Gamma)$ with
%$\ka(\bx) >0$ and $\ka(\cdot)$  continuous at $\bx$.
Suppose also that
 $\xi$ is $\ka(\bx)$-homogeneously stabilizing at $\bx$,
and that Assumption A4 or A5 holds.
Then  
for any  $\bz \in \R^d$, we have 
\bea
\xi_\la(\bx+ \la^{-1/d}\bz,T; \Po_\la ,\R^d)  \tod 
 \xi_\infty^{\bx,T} (\H_{\kappa(\bx)},\R^d) ~~{\rm as}~\la \to \infty.
\lbl{0216c}
\eea
Also, if $f \in B(\R^d)$ and $f$ is continuous at $\bx$,
%or if Assumption A1 holds,
%and $\bx$ is a Lebesgue point of $f$,
then
\bea
\langle f,\xi_\la(\bx+ \la^{-1/d}\bz,T; \Po_\la )   \rangle
\tod
 f(\bx) \xi_\infty^{\bx,T} (\H_{\kappa(\bx)},\R^d) ~~{\rm as}~\la \to \infty.
% \langle f,\xi^*_\la(\bx+ \la^{-1/d}\bz,T; \Po_\la )   \rangle
%\toP 0.
\lbl{50209f}
\eea
\end{lemm}
{\em Proof.}
The proof of \eq{0216c} is related to that of Lemma \ref{couplem0a}.
Given $\bx$ and $\bz$, set $\bv_\la := \bx + \la^{-1/d}\bz$.
For Borel $A \subseteq \R^d$,
define $g_A:\R^d \times \MM \times \LL \to \R$ by
$$
g_A(\bw,t,\X) = \xi( \bx+\bw,t; \bx + \bw + \X, \bx + \bw + A).
$$
%it is the case that
Then
\bean
\xi_\la(
\bv_{\la},T ; \Po_\la , \bv_{\la}
  + \la^{-1/d} A)
= g_A ( \la^{-1/d} \bz,T,\la^{1/d}(
-\bv_\la + \Po_\la ) ).
\eean
Taking our topology on 
$\R^d \times \MM \times \LL $ to be the product of the Euclidean
topology on $\R^d$, the discrete topology on $\MM$ and the 
topology induced by the metric $\dist$ on $\LL$ which was defined 
at \eq{ppmet}, we assert that  as $\la \to \infty$,
\bea 
 ( \la^{-1/d} \bz,T,\la^{1/d}(- \bv_\la
%\bx - \la^{-1/d} \bz 
+ \Po_\la ) )
\tod 
(\0, T, \H_{\ka(\bx)}).
\lbl{050404c}
\eea
To see this, for each $\la$, let $\Po'_\la$ and  $\H'_{\ka(\bx)}$
be the coupled realisations of  $\Po_\la$ and $\H_{\ka(\bx)}$
obtained by  the pivoted coupling
 %with pivot at $\bx+ \la^{-1/d}\bz$.
 with pivot at $\bv_{\la}$.
Then for $\eps >0$, by Lemma \ref{pivlem} we have
$P[\dist(\Po'_\la,\H'_{\ka(\bx)})>\eps] \to 0$ as
$\la \to \infty$. This gives us  \eq{050404c}.

If Assumption A4 (translation invariance) holds, 
then the functional $g_A(\bw,t,\X)$ does not depend on $\bw$, 
so that $g_A(\bw,t,\X) = g_A(\0,t,\X)$ and 
by the assumption  
 that $\xi$ is 
$\ka(\bx)$-homogeneously stabilizing at $\bx$, we have
that $(\0,T,\H_{\ka(\bx)})$ almost surely lies at a continuity
point of the functional $g_A$.

If, instead,  Assumption A5 (continuity) holds,
take $A=\R^d$ or $A =B_K$ or
 $A=\R^d \setminus B_K$, with $K >K_1$ and 
$K_1$ given in Definition \ref{defctsmeas}.
  Then by the assumption  that $\xi$ is 
$\ka(\bx)$-homogeneously stabilizing at $\bx$ (see \eq{homstab2}),
with probability 1 there exists
a finite (random) $\eta >0$ such that for $\dist(\X,\H_{\ka(\bx)}) < \eta$,
and for  $|\bw| < \eta$,
\bean
g_A(\bw,T,\X) =  \xi(\bx +\bw,T; \bx + \bw
 + (\H_{\ka(\bx)} \cap B_{1/\eta}), 
\bx + \bw + A)
\\
\to
\xi( \x,T; \bx +( \H_{\ka(\bx)} \cap B_{1/\eta} ), \bx + A) = 
g_A(\0,T, \H_{\ka(\bx)})
\mbox{ as } \bw \to \0.
\eean
Hence, 
 $(\0,T,\H_{\ka(\bx)})$ almost surely lies at a continuity
point of the mapping $g_A$ in this case too.

Thus, 
if $A$ is $\R^d$ or $B_K$ or
 $\R^d \setminus B_K$, for any $K$ 
under  A4 and for $K>K_1$ under A5,
 the mapping $g_A$ satisfies the
conditions for the Continuous Mapping Theorem, and this
with \eq{050404c}  gives us 
\bea
%\xi_\la(\bx+ \la^{-1/d}\bz,T; \Po_\la ,\bx + \la^{-1/d} \bz, A)  
\xi_\la(\bv_\la,T; \Po_\la , \bv_\la + \la^{-1/d} A)  
\tod 
 \xi_\infty^{\bx,T} (\H_{\kappa(\bx)},A) ~~{\rm as}~\la \to \infty.
\lbl{050404d}
\eea
Taking $A =\R^d$ in \eq{050404d} gives us \eq{0216c}.
 
%In the case where A1 holds, \eq{50209f} is immediate from
%\eq{0216c}. 
Now suppose that $f$ is continuous at $\bx$.
To derive \eq{50209f} in this case, note first that 
\bean
\langle f,
  %\xi_{\la}(\bx+ \la^{-1/d}\bz,T;\Po_\la)
  \xi_{\la}(\bv_\la,T;\Po_\la)
-  
  %\xi^*_{\la}(\bx+ \la^{-1/d}\bz,T;\Po_\la) \rangle 
  \xi^*_{\la}(\bv_\la,T;\Po_\la) \rangle 
%\nonumber \\
 = 
  %\int_{\R^d} ( f(\bw) - f(\bx + \la^{-1/d} \bz) )
  \int_{\R^d} ( f(\bw) - f(\bv_\la) )
%\xi_\la(\bx + \la^{-1/d} \bz,T; \Po_\la,d\bw ).
\xi_\la(\bv_\la,T; \Po_\la,d\bw ).
%\lbl{050404e}
\eean
%We shall choose $K>0$ and split the region of integration
%in \eq{050404e} into 
% $B_{\la^{-1/d} K}(\bx+ \la^{-1/d} \bz)$ 
%and  the complementary region
% $\R^d \setminus B_{\la^{-1/d} K}(\bx+ \la^{-1/d} \bz)$. 
%We consider the latter region first. 
Given $K>0$, 
by \eq{050404d}
 we have
\bean
 \left|
 %\int_{\R^d \setminus B_{\la^{-1/d} K}(\bx + \la^{-1/d} \bz)}
 \int_{\R^d \setminus B_{\la^{-1/d} K}(\bv_\la)}
 %( f(\bw) - f(\bx+ \la^{-1/d} \bz) )
 ( f(\bw) - f(\bv_\la) )
%\xi_\la(\bx+ \la^{-1/d}\bz; \Po_\la,d\bw ) \right|
\xi_\la(\bv_\la,T; \Po_\la,d\bw ) \right|
~~~~~~~~~~~~~~~~~~
~~~~~~~~~~~~~~~~~~
\\
\leq 
2 \|f\|_\infty 
 %\xi_\la(\bx+ \la^{-1/d}\bz,T;\Po_\la,
 \xi_\la(\bv_\la,T;\Po_\la,
\R^d \setminus B_{\la^{-1/d}K}(\bv_\la  ) )
%\R^d \setminus B_{\la^{-1/d}K}(\bx + \la^{-1/d} \bz  ) )
%\nonumber
% \\
\tod 
2 \|f\|_\infty 
 \xi_\infty^{\bx,T}(\H_{\ka(\bx)},\R^d \setminus B_{K}  ),
\eean 
where the limit is almost surely finite and converges in probability
to zero as $K \to \infty$. Hence for  $\eps >0$, we have
\bea
\lim_{K \to \infty} \limsup_{\la \to \infty} P\left[
 \left|
 %\int_{\R^d \setminus B_{\la^{-1/d} K}(\bx+ \la^{-1/d} \bz)}
 \int_{\R^d \setminus B_{\la^{-1/d} K}(\bv_\la)}
 %( f(\bw) - f(\bx+ \la^{-1/d} \bz) )
 ( f(\bw) - f(\bv_\la) )
%\xi_\la(\bx+ \la^{-1/d} \bz; \Po_\la,d\bw ) \right| > \eps \right] =0.
\xi_\la(\bv_\la,T; \Po_\la,d\bw ) \right| > \eps \right] =0.
~~~
\lbl{050404f}
\eea
Also,
given $K>0$, it is the case that
\bea
 %\left| \int_{ B_{\la^{-1/d} K}(\bx+ \la^{-1/d} \bz)} ( f(\bw) - f(\bx) )
 \left| \int_{ B_{\la^{-1/d} K}(\bv_\la)} ( f(\bw) - f(\bv_\la) )
%\xi_\la(\bx+ \la^{-1/d}\bz; \Po_\la,d\bw ) \right|
\xi_\la(\bv_\la,T; \Po_\la,d\bw ) \right|
\nonumber \\
\leq
2
\phi_{\la^{-1/d}(K+|\bz|)}(\bx)
 %|\xi_\la|(\bx,\X_m,\R^d)  |
 \xi_\la(\bv_\la,T;\Po_\la,\R^d)  
\lbl{050404g}
\eea
and by  continuity of $f$,
%ous at $\bx$, then
$\phi_{\la^{-1/d}(K+ |\bz|)}(\bx) \to 0$
%,
 while 
 %$|\xi_\la|(\bx,\X_m,\R^d)  $ if signed 
 % $\xi_\la(\bx+ \la^{-1/d} \bz,T;\Po_\la,\R^d)  $ 
 $\xi_\la(\bv_\la,T;\Po_\la,\R^d)  $ 
converges in distribution
to the finite random variable
 $\xi_\infty^{\bx,T}(\H_{\ka(\bx)},\R^d)  $ 
by 
\eq{0216c}, and hence the right hand side of
\eq{050404g} tends to zero in probability as $\la \to \infty$.
Combined with \eq{050404f}, this gives us 
\bea
\langle f,
  %\xi_{\la}(\bx+ \la^{-1/d}\bz,T;\Po_\la)
  \xi_{\la}(\bv_\la,T;\Po_\la)
-  
  %\xi^*_{\la}(\bx+ \la^{-1/d}\bz,T;\Po_\la) \rangle \toP 0
  \xi^*_{\la}(\bv_\la,T;\Po_\la) \rangle \toP 0.
\lbl{050404h}
\eea
%Under Assumption A2, 
%for Borel $A \subseteq \R^d$,
% the change of variables $\bw = \bx + \la^{-1/d} \bz
%+
%\la^{-1/d}(\by-\bx - \la^{-1/d} \bz )$ yields 
%\bean
%\xi_\la(\bx;\X,A)  = 
%%\1_{\Gamma_\la}(\bx) 
%%\xi 
%%(\bx; \bx+\la^{1/d}(-\bx +\X),\bx + \la^{1/d}(-\bx + A) )
%%\nonumber \\
%%= 
%%\1_{\Gamma_\la}(\bx)
% \int_{\bx + \la^{-1/d} \bz + \la^{1/d} (-\bx - \la^{-1/d} \bz +A)}
%d\by
%\nonumber \\
%\times
% \xi'
%(\bx+ \la^{-1/d} \bz; \bx+ \la^{-1/d}\bz+ \la^{1/d}(-\bx - \la^{-1/d} \bz
%+\X),\by ) 
%\nonumber  \\
%= \la
%% \1_{\Gamma_\la}(\bx) 
%\int_{A} \xi'
%(\bx + \la^{-1/d} \bz; \bx + \la^{-1/d} \bz+\la^{1/d}
%(-\bx - \la^{-1/d} \bz +\X),
%\nonumber \\
%\bx + \la^{-1/d} \bz + \la^{1/d}(\bw -\bx - \la^{-1/d} \bz) ) d\bw .
%%\lbl{50210a}
%\eean
%Hence, under A2,
%\bean
% | \int_{ B_{\la^{-1/d} K}(\bx+ \la^{-1/d} \bz)} ( f(\bw) - f(\bx) )
%\xi_\la(\bx+ \la^{-1/d} \bz; \Po_\la,d\bw ) |
%\nonumber \\
%= \la 
% | \int_{ B_{\la^{-1/d} K}(\bx+ \la^{-1/d} \bz)} ( f(\bz) - f(\bx) )
%\xi'
%(\bx+ \la^{-1/d} \bz; \nonumber \\
%\bx+ \la^{-1/d} \bz+\la^{1/d}(-\bx - \la^{-1/d}
% \bz +\Po_\la),\bx+ \la^{-1/d} \bz + \la^{1/d}(\bz -\bx - \la^{-1/d} \bw) )
% d\bw |
%\nonumber \\
%\leq K_0 \la 
% | \int_{ B_{\la^{-1/d} K}(\bx+ \la^{-1/d})} | f(\bw) - f(\bx) | d\bw |
%\eean
%and if additionally $\bx$ is a Lebesgue point of $f$ then
%this tends to zero.  
%Combined with \eq{050404f}, this gives us 
%\eq{050404h} in the case where A2 holds and $\bx$ is a Lebesgue point
%of $f$. 
%
%
%Under either A4 or A5, we have established \eq{050404h}.
Also, by \eq{0216c} and continuity of $f$ at $\bx$, we have
%in the case where $f$ is continuous at $\bx$. 
\bean
\langle f,
  \xi_{\la}^*(\bv_\la,T;\Po_\la) \rangle
\tod f(\bx) \xi_\infty^{\bx,T}(\H_{\ka(\bx)},\R^d),
\eean
and combined with \eq{050404h} this yields \eq{50209f}.
%Combining this with \eq{0216c} and using Slutsky's theorem,
%we obtain \eq{50209f}.
$\qed$ \\

\begin{lemm}
\lbl{lemweak2}
Suppose $\xi$ satisfies
 Assumption $A4$  or
%(translation invariance) or Assumption
 $A5$.
% (almost everywhere continuity). 
Then for Lebesgue-almost all $\bx\in \Gamma_0$ 
and all $\bz \in \R^d$,
as $\la \to \infty$ we have
\bea
 \xi_\la(\bx,T; \Po_\la^{\bx+ \la^{-1/d} \bz,T'},\R^d) 
 \xi_\la(\bx+ \la^{-1/d} \bz,T'; \Po_\la^{\bx,T},\R^d) 
\nonumber \\
\tod
 \xi_\infty^{\bx,T}( \H^{\bz,T'}_{\kappa(\bx)},\R^d )
\xi_\infty^{\bx,T'}(-\bz+  \H^{\0,T}_{\kappa(\bx)},\R^d ).
% ] 
\lbl{0426m} \eea
Also, for $f \in B(\R^d)$, if $f$ is continuous at $\bx$, then
% or Assumption A1 holds, then
\bea
 \langle f, \xi_\la(\bx,T;\Po_\la^{\bx + \la^{-1/d} \bz,T'}) \rangle
\times \langle f, \xi_\la(\bx+\la^{-1/d}\bz,T';\Po_\la^{\bx ,T}) \rangle
\nonumber
\\
%- \langle f, \xi^*_\la(\bx,T;\Po_\la^{\bx + \la^{-1/d} \bz,T'}) \rangle
%\toP 0;
%\lbl{top0a}
\tod f(\bx)^2 \xi_\infty^{\bx,T}(\H_{\ka(\bx)}^{\bz,T},\R^d)
 \xi_\infty^{\bx,T'}(-\bz + \H_{\ka(\bx)}^{\0,T},\R^d).
%- \langle f, \xi^*_\la(\bx+\la^{-1/d}\bz,T';\Po_\la^{\bx,T}) \rangle
%\toP 0.
\lbl{top0b}
\eea
\end{lemm}
{\em Proof.}
Again write $\bv_\la$ for $\bx+ \la^{-1/d} \bz$.
Let $A \subseteq \R^d$ be a Borel set.
Define the function $\tg_A:
\R^d  \times \MM \times \MM \times \LL \to \R^2$ by
\bean
\tg_A(\bw,t,t',\X) = (\xi(\bx,t;\bx +  \X^{\bz,t'},\bx + A),
\\
 \xi(\bx+ \bw,t';\bx + \bw -\bz + \X^{\0,T},\bx + \bw + A) ).
\eean
Then
\bean
%( \xi_\la(\bx,T;\Po_\la^{\bx + \la^{-1/d} \bz,T'},x + \la^{-1/d}A),
( \xi_\la(\bx,T;\Po_\la^{\bv_\la,T'},\bx + \la^{-1/d}A),
 %\xi_\la(\bx+\la^{-1/d}\bz,T';\Po_\la^{\bx ,T},\bx+ \la^{-1/d} \bz +
 \xi_\la(\bv_\la,T';\Po_\la^{\bx ,T},
%\bx+ \la^{-1/d} \bz
\bv_\la
 + \la^{-1/d} A) ) 
\\
= ( \xi(\bx, T;
\bx+ \la^{1/d}(-\bx + 
%\Po_\la^{\bx + \la^{-1/d}\bz,T'};\bx +A ),
\Po_\la^{\bv_\la,T'}),\bx +A ),
~~~~~~~~~~~~~~~~~~~~~
\\
%\times
 \xi(
%(\bx+\la^{-1/d}\bz,T'; \bx + \la^{-1/d} \bz 
\bv_\la,T'; \bv_\la 
+ \la^{1/d}( -\bx -\la^{-1/d} \bz + 
\Po_\la^{\bx ,T}),
%\bx+ \la^{-1/d} \bz +
\bv_\la + A) ) 
\\
= \tg_A ( \la^{-1/d} \bz,T,T',\la^{1/d}  (- \bx + \Po_\la) ).
\eean
Under A5, let us restrict attention to the case where $A$ is
$\R^d$, $B_K$ or $\R^d \setminus B_K$ with $K>K_1$.
Then under either A4 or A5,
by  similar arguments to those used in proving Lemma \ref{lemweak1},
$(\0,T,T',\H_{\ka(\bx)})$ lies almost surely
at a continuity point of $\tg_A$, and since 
$
\la^{-1/d} (-\bx + \Po_\la) \tod \H_{\ka(\bx)},$
the Continuous Mapping Theorem gives us
\bea
 %(\xi_\la(\bx,T;\Po_\la^{\bx + \la^{-1/d} \bz,T'},x + \la^{-1/d}A),
 (\xi_\la(\bx,T;\Po_\la^{\bv_\la,T'},\bx + \la^{-1/d}A),
 %\xi_\la(\bx+\la^{-1/d}\bz,T';\Po_\la^{\bx ,T},\bx+ \la^{-1/d} \bz +
 \xi_\la(\bv_\la,T';\Po_\la^{\bx ,T},\bv_\la  +
\la^{-1/d} A) )
\nonumber \\
\tod \tilde{g}_{A}(\0,T,T',\H_{\ka(\bx)})
= (\xi_\infty^{\bx,T}(\H_{\ka(\bx)}^{\bz,T'}, A), 
 \xi_\infty^{\bx,T'}(- \bz+ \H_{\ka(\bx)}^{\0,T}, A)) 
\lbl{050404i}
\eea
as $\la \to \infty$.
Taking $A= \R^d $ gives us \eq{0426m}.

Now suppose $f$ is continuous at $\bx$.
% and
% let $\phi_\eps(\bx):= \sup\{|f(\by)-f(\bx)|: \by \in B_\eps (\bx)\}$.  
%We need 
To prove  \eq{top0b} in this case,
 observe that for $K>0$,
we have
\bea
|\langle f, \xi_\la(\bx,T; 
%\Po_\la^{\bx + \la^{-1/d} \bz,T'}) 
\Po_\la^{\bv_\la,T'}) 
-
 %\xi^*_\la(\bx,T; \Po_\la^{\bx + \la^{-1/d} \bz,T'}) 
 \xi^*_\la(\bx,T; \Po_\la^{\bv_\la,T'}) 
\rangle |
%\nonumber \\
%\leq \left| \int_{B_{\la^{-1/d}K}(\bx)}
%  (f(\by) - f(\bx )  ) \xi_\la(\bx ,T;
%\tPo_\la^{\bx + \la^{-1/d} \bz,T'} ,d\by) \right| \nonumber \\
%+ \left| \int_{\R^d\setminus B_{\la^{-1/d} K}(\bx)} 
%( f(\by) - f(\bx )) \xi_\la(\bx , T;
%\tPo_\la^{\bx + \la^{-1/d} \bz,T'}
%,d\by) \right|. 
\leq \phi_{\la^{-1/d} K}(\bx)
 \xi_\la(\bx,T;\Po_\la^{\bv_\la,T'},\R^d)
\nonumber \\
+ 2 \|f\|_\infty
 \xi_\la(\bx,T; \Po_\la^{\bv_{\la},T'},
\R^d \setminus  B_{\la^{-1/d} K}(\bx)  )
~~~~~
\lbl{050323a}
\eea
The first term in the right hand
side of \eq{050323a}
% is bounded by
%\bean
%\phi_{\la^{-1/d} K}(\bx)
% \xi_\la(\bx,T;\Po_\la^{\bx + \la^{-1/d} \bz,T'},\R^d)
%\eean
%which
 tends to zero in probability for any fixed $K$,
by \eq{050404i} and the fact that
$\xi_\infty^{\bx,T}(\H_{\ka(\bx)}^{\bz,T'},\R^d)$ is almost surely 
finite. Also by \eq{050404i}, 
the second term in the right hand
side of \eq{050323a} 
converges in distribution, as $\la \to \infty$, to
% is bounded by 
%\bean
%2 \|f\|_\infty
% \xi_\la(\bx,T; \Po_\la^{\bx + \la^{-1/d}\bz,T'},
%\R^d \setminus  B_{\la^{-1/d} K}(\bx)  )
%\\
%\tod
$2 \|f\|_\infty
 \xi_\infty^{\bx,T}( \H_{\ka(\bx)}^{\bz,T'}, \R^d\setminus B_K),
$
%\eean
which tends to zero in probability as $K\to \infty$. Hence,
by \eq{050323a} we obtain
\bea
\langle f, \xi_\la(\bx,T; \Po_\la^{\bv_\la,T'})
- \xi^*_\la(\bx,T; \Po_\la^{\bv_\la,T'})
 \rangle
\toP 0.
\lbl{050324d}
\eea
We also have
\bea
|\langle f, \xi_\la(\bv_\la,T'; \Po_\la^{\bx ,T}) 
-  \xi^*_\la(\bv_\la,T'; \Po_\la^{\bx ,T}) 
\rangle |
%\\
%- f(\bx) \xi_\infty^{\bx,T}(-\bz + \tH_{\ka(\bx)}^{\0,T},\R^d) |
\nonumber \\
\leq \left|
 \int_{B_{\la^{-1/d}K}(\bv_\la)}
  (f(\by) - f(\bx )  ) \xi_\la(\bv_\la,T';
\Po_\la^{\bx ,T}
,d\by) \right|
 \nonumber \\
%+ \left| \int_{\R^d\setminus B_{\la^{-1/d} K}(\bx+ \la^{-1/d} \bz)} 
%( f(\by) - f(\bx )) \xi_\la(\bx + \la^{-1/d} \bz , T';
%\Po_\la^{\bx ,T}
%,d\by) \right|. 
+ 2 \|f\|_\infty
 \xi_\la(\bv_\la,T'; \Po_\la^{\bx,T},\R^d 
\setminus B_{\la^{-1/d} K}(\bv_\la) ).
\lbl{050331a}
\eea
By \eq{050404i} and the assumed continuity of $f$ at $\bx$,
the first term in the right side of \eq{050331a} tends to zero
in probability for any fixed $K$, 
%by along with and the assumption that 
%$\xi_\infty^{\bx,T}(\H_{\ka(\bx)},\R^d)$ is almost surely 
%finite. Also by \eq{050404i},
while the second term
 %in the right side of \eq{050331a}
 converges in distribution to
%is bounded by
$
%\bean
%\tod
2 \|f\|_\infty \xi_\infty^{\bx,T'}
(-\bz + \tH_{\ka(\bx)}^{\0,T},\R^d\setminus B_K),
$
%\eean
which tends to zero in probability as $K \to \infty$. Hence, 
as $\la \to \infty$ we have
\bea
\langle f, \xi_\la(\bv_\la,T'; \Po_\la^{\bx ,T}) 
-  \xi^*_\la(\bv_\la,T'; \Po_\la^{\bx ,T}) 
\rangle  \toP 0.
\lbl{050406a}
\eea
By continuity of $f$ at $\bx$, and the case $A=\R^d$ of
\eq{050404i}, we have
\bean
(\langle f, \xi_\la^*(\bx,T; \Po_\la^{\bv_\la,T'}) \rangle,
\langle f, \xi_\la^*(\bv_\la,T'; \Po_\la^{\bx,T}) \rangle)
~~~~~~
~~~~~~
~~~~~~
\\
\tod
( f(\bx) \xi_\infty^{\bx,T}(\H_{\ka(\bx)}^{\bz,T'},\R^d )  , 
%\toP 0 ; - 
f(\bx) \xi_\infty^{\bx,T'}(-\bz + \H_{\ka(\bx)}^{\0,T} ,\R^d)  ).
%\rangle \toP 0.
\eean
Combining this with
% \eq{050404i},
 \eq{050324d} and \eq{050406a} yields \eq{top0b}.
% the left
%side of \eq{050331a} tends to zero in probability.
%Combined with \eq{050324d}, 
% and \eq{0426m} this gives us  \eq{top0b}. (CHECK)
  $\qed$ \\ 

The  following lemma is a refinement of
Lemma \ref{couplem0c} and is proved in the same manner as that result.
 It will be  used for de-Poissonizing our
central limit theorems.
% i.e. for extending
% those  results for $\mu_\la^{\xi}$ to $\nu_{\la,n}^\xi$.
To ease notation, we do not mention the marks in the
notation for the statement and proof of this result.
\begin{lemm}
\lbl{lemweak3}
Let $(\bx,\by) \in \Gamma_0^2$ with $\bx \neq \by$,
 and let $(\bz,\bw) \in (\R^d)^2$. 
%Let $\bx$, $\by$, $\bz$ and $\bw$ in $\R^d$, with
%$\bx \in \inter(\Gamma)$,  $\by \in  \inter{\Gamma}$, and $\bx \neq \by$.
Suppose either that Assumption A1 holds, 
%or that Assumption A2 holds and
%$\bx$  and $\by$ are Lebesgue points of $f$, 
or that $\bx$
and $\by$ are continuity points of $f$.
%Also assume $\ka(\bx)>0$, $\ka(\by)>0$, and $\bx$ and $\by$
%are continuity points of $\ka$. Let $\H_{\ka(\bx)}$ and $\tH_{\ka(\by)}$
%be independent homogeneous Poisson point processes in $\R^d$
%with intensity $\ka(\bx)$ and $\ka(\by)$ respectivley.
%
Given  integer-valued functions $(\ell(\la), \la \geq 1)$ and
 $(m(\la), \la \geq 1)$  with $\ell (\la) \sim \la$ and $m(\la) \sim \la $
as $\la \to \infty$,
% for all $\la$,
we have convergence in joint distribution, as
$\la \to \infty$, of the $11$-dimensional
random vector
\bean
\left( \langle f, \xi_\la(\bx;\X_{\ell})\rangle ,
\langle f, \xi_\la(\bx;\X_\ell^{\by})\rangle, 
\langle f, \xi_\la(\bx;\X_\ell^{\bx+ \la^{-1/d}\bz})\rangle, 
\langle f, \xi_\la(\bx;\X_\ell^{\bx+ \la^{-1/d}\bz } \cup \{\by\} )\rangle, 
~^{~^{~^{~^{~}}}} 
\right.
%\nonumber \\
\\
\langle f,  \xi_\la(\bx;\X_{m})\rangle,
 \langle f, \xi_\la(\bx;\X_m^\by)\rangle,
 \langle f, \xi_\la(\bx;\X_m^\by \cup \{\bx + \la^{-1/d} \bz\})\rangle,
 \langle f, \xi_\la(\by;\X_m)\rangle,
%\nonumber \\
\\
 \langle f, \xi_\la(\by;\X_m^\bx)\rangle,
\langle f, \xi_\la(\by;\X_m^\bx \cup \{\bx + \la^{-1/d} \bz\})\rangle,
\\
\left.
~^{~^{~^{~^{~}}}} 
 \langle f, \xi_\la(\by;\X_m^\bx
\cup \{\bx + \la^{-1/d}\bz,\by+ \la^{-1/d} \bw \})\rangle  
%~^{~^{~^{~^{~}}}} 
\right)
\eean
to
\bean
%\nonumber \\
%\tod 
\left(
 f(\bx) 
\xi_\infty^\bx(\H_{\ka(\bx)},\R^d),
f(\bx) \xi_\infty^\bx(\H_{\ka(\bx)},\R^d),
f(\bx) \xi_\infty^\bx(\H_{\ka(\bx)}^{\bz},\R^d),
f(\bx) \xi_\infty^\bx(\H_{\ka(\bx)}^{\bz},\R^d),
~^{~^{~^{~^{~}}}} 
\right.
\nonumber \\
%\nonumber \\
f(\bx) \xi_\infty^\bx(\H_{\ka(\bx)},\R^d),
 f(\bx) \xi_\infty^\bx (\H_{\ka(\bx)},\R^d),
 f(\bx) \xi_\infty^\bx (\H_{\ka(\bx)}^{\bz},\R^d),
 f(\by) \xi_\infty^\by (\tH_{\ka(\by)},\R^d),
%~^{~^{~^{~^{~}}}} 
%\right. ~~~
\nonumber \\
%\left.
\left.
~^{~^{~^{~^{~}}}} 
 f(\by) \xi_\infty^\by (\tH_{\ka(\by)},\R^d),
 f(\by) \xi_\infty^\by (\tH_{\ka(\by)},\R^d),
 f(\by) \xi_\infty^\by (\tH_{\ka(\by)}^\bw,\R^d) \right).
%\lbl{bigweak}
\eean
\end{lemm}
{\em Proof.}
First, we assert that
\bean
\left(  \xi_\la(\bx;\X_{\ell},\R^d),
\xi_\la(\bx;\X_\ell^{\by},\R^d) ,
 \xi_\la(\bx;\X_\ell^{\bx+ \la^{-1/d}\bz},\R^d), 
 \xi_\la(\bx;\X_\ell^{\bx+ \la^{-1/d}\bz, } \cup \{\by\},\R^d ),
~^{~^{~^{~^{~}}}} 
\right.
\nonumber \\
 \xi_\la(\bx;\X_{m},\R^d),
 \xi_\la(\bx;\X_m^\by,\R^d),
  \xi_\la(\bx;\X_m^\by \cup \{\bx + \la^{-1/d} \bz\},\R^d),
% \langle f, \xi_\la(\by;\X_m)\rangle,
%\nonumber \\
\xi_\la(\by;\X_m,\R^d),
%~~~~~
\nonumber \\
\xi_\la(\by;\X_m^\bx,\R^d),
  \xi_\la(\by;\X_m^\bx \cup \{\bx + \la^{-1/d} \bz\},\R^d),
\\
\left.
~^{~^{~^{~^{~}}}} 
\xi_\la(\by;\X_m^\bx
\cup \{\bx + \la^{-1/d}\bz,\by+ \la^{-1/d} \bw \},\R^d)  
%~^{~^{~^{~^{~}}}} 
\right)
\eean
converges in distribution to
\bean
%\nonumber \\
%\tod 
\left(
\xi_\infty^\bx(\H_{\ka(\bx)},\R^d),
 \xi_\infty^\bx(\H_{\ka(\bx)},\R^d),
 \xi_\infty^\bx(\H_{\ka(\bx)}^{\bz},\R^d),
 \xi_\infty^\bx(\H_{\ka(\bx)}^{\bz},\R^d),
~^{~^{~^{~^{~}}}} 
\right.
%\nonumber \\
\nonumber \\
 \xi_\infty^\bx(\H_{\ka(\bx)},\R^d),
  \xi_\infty^\bx (\H_{\ka(\bx)},\R^d),
  \xi_\infty^\bx (\H_{\ka(\bx)}^{\bz},\R^d),
% f(\by) \xi_\infty^\bx (\tH_{\ka(\by)},\R^d),
  \xi_\infty^\by (\tH_{\ka(\by)},\R^d),
%~^{~^{~^{~^{~}}}} 
%\right. ~~~
\nonumber \\
%\left.
\left.
~^{~^{~^{~^{~}}}} 
  \xi_\infty^\by (\tH_{\ka(\by)},\R^d),
  \xi_\infty^\by (\tH_{\ka(\by)},\R^d),
  \xi_\infty^\by (\tH_{\ka(\by)}^\bw,\R^d) \right).
%\lbl{bigweak}
\eean
This is deduced from Lemma \ref{pplemweak} by a similar argument to
the proof of Lemma \ref{couplem0a}. For example, considering just
the third component,
 defining the mapping $h_{\R^d,\bx}$ on $\MM \times \LL$
by \eq{050404a}, 
%by
%$$
%h_4(t,\X) = \xi (\bx,t;\bx + \X, \R^d ),
%$$
we have
\bean
 \xi_\la(\bx;\X_\ell^{\bx+ \la^{-1/d}\bz},\R^d)
= h_{\R^d,\bx}(T,\la^{1/d}(-\bx + \X_\ell^{\bx + \la^{-1/d} \bz}) ) 
\\
= h_{\R^d,\bx}(T,\{\bz\} \cup \la^{1/d}(-\bx + \X_\ell) ),
\eean
and  by Lemma \ref{pplemweak}, $(T, \{\bz\} \cup \la^{1/d}(-\bx + \X_\ell) )$
converges in distribution to $(T, \H_{\ka(\bx)}^{\bz})$ which
is almost surely at a continuity point of $h_{\R^d,\bx}$. Similar arguments
apply for the other components and give us the assertion above.
This assertion implies that the result holds under A1, i.e.
when $\xi = \xi^*$. 

Now let us drop Assumption A1, but assume 
 that $\bx$ and $\by$ are continuity points of $f$.
Then by Lemma \ref{couplem0b},
\bean
\langle f, \xi_\la(\bx;\X_\ell) 
%\rangle - \langle f,
-
 \xi^*_\la(\bx;\X_\ell) \rangle \toP 0,
%\nonumber
%\lbl{050404b}
% \\
~~~
\langle f, \xi_\la(\bx;\X_\ell^\by) 
%\rangle - \langle f,
-
 \xi^*_\la(\bx;\X_\ell^\by) \rangle \toP 0,
\nonumber \\
\langle f, \xi_\la(\bx;\X_m) 
%\rangle - \langle f,
-  \xi^*_\la(\bx;\X_m) \rangle 
\toP 0, ~~~
\langle f, \xi_\la(\bx;\X_m^\by) 
%\rangle - \langle f,
-  \xi^*_\la(\bx;\X_m^{\by}) \rangle 
\toP 0, 
\\
\langle f, \xi_\la(\by;\X_m) 
%\rangle - \langle f,
-  \xi^*_\la(\by;\X_m) \rangle 
\toP 0, ~~~
\langle f, \xi_\la(\by;\X_m^\bx) 
%\rangle - \langle f,
-  \xi^*_\la(\by;\X_m^\bx) \rangle 
\toP 0,
%\nonumber \\
%\lbl{050331b}
\eean
and a similar argument to the proof of \eq{050324d} 
(working with $\X_m$ instead of $\Po_\la$) yields
\bean
\langle f, \xi_\la(\bx;\X_\ell^{\bx + \la^{-1/d} \bz})
-
 \xi^*_\la(\bx;\X_\ell^{\bx + \la^{-1/d} \bz}) \rangle  \toP 0.
%~~~
\eean
Very similar arguments (which we omit) yield
% to Lemma \ref{couplem0b} yield
\bean
%\nonumber \\
\langle f, \xi_\la(\bx;\X_\ell^{\bx + \la^{-1/d} \bz}\cup\{\by\}) 
%\rangle - \langle f,
-
 \xi^*_\la(\bx;\X_\ell^{\bx + \la^{-1/d} \bz} \cup \{\by\})
 \rangle  \toP 0,
\nonumber \\ 
\langle f, \xi_\la(\bx;\X_m^{\bx + \la^{-1/d} \bz}\cup\{\by\}) 
%\rangle - \langle f,
-
 \xi^*_\la(\bx;\X_m^{\bx + \la^{-1/d} \bz} \cup \{\by\})
 \rangle  \toP 0,
\nonumber \\
%\langle f, \xi_\la(\by;\X_m) \rangle -
%\langle f, \xi^*_\la(\by;\X_m) \rangle 
%\toP 0,
%\nonumber \\
\langle f, \xi_\la(\by;\X_m^\bx\cup \{\bx + \la^{-1/d} \bz\}) 
%\rangle - \langle f,
- 
 \xi^*_\la(\by;\X_m^\bx\cup \{\bx + \la^{-1/d} \bz\}) \rangle
\toP 0,
\nonumber \\
\langle f, \xi_\la(\by;\X_m^\bx\cup \{\bx + \la^{-1/d} \bz,
\by + \la^{-1/d} \bw\}) 
~~~~~
%\rangle 
\nonumber \\
-
%\langle f,
 \xi^*_\la(\by;\X_m^\bx\cup \{\bx + \la^{-1/d} \bz,
\by + \la^{-1/d} \bw\}) \rangle 
\toP 0.
%\lbl{50210a}
\eean
%All of these are proved in a similar manner to \eq{050331b},
%and  
Combining these eleven convergence in probability statements 
with the fact that we have
established our conclusion 
%\eq{bigweak} holds
 in the case 
where Assumption A1 
($\xi =\xi^*$) holds, and using Slutsky's theorem, we obtain our
conclusion 
% the result
%\eq{bigweak}
 in the other case as well. 
$\qed$ \\

\section{Proof of Theorem \ref{thlln1}}
\allco
Before proving Theorem \ref{thlln1}, we give general expressions
for the first two  moments of $\langle f, \mu_\la^\xi \rangle $ which
we shall use again later on.  
By Palm theory for the Poisson process
 (e.g. a slight generalization of Theorem 1.6 of \cite{Pe}), 
we have
\bea
\E  \langle f, \mu^\xi_\la \rangle = 
 \E  \sum_{(\bx,t) \in \Po_\la} \langle f, \xi_\la(\bx,t; \Po_\la) \rangle 
%\nonumber \\
= \la \E \langle
 f,  \xi_\la(\bX,T; \Po_\la) \rangle 
\lbl{0216a0}
\eea
and
\bea
\E [ \langle f, \mu^\xi_\la \rangle^2 ] 
%=
 %\E \left[ \left( \sum_{(\bY,T) \in \Po_\la} 
%\langle f,
% \xi_\la(\bY,T; \Po_\la)\rangle \right)^2
%\right] 
%\nonumber \\
=   
  \left( \E\sum_{(\bx,t) \in \Po_\la}
\langle f, \xi_\la(\bx,t; \Po_\la)\rangle^2 \right) 
\nonumber \\
+ 2 \E \sum_{\{(\bx,t),(\by,u)\} \subseteq \Po_\la} 
\langle f,
 \xi_\la(\bx,t; \Po_\la) \rangle
\langle f  ,\xi_\la(\by,u; \Po_\la) \rangle
\nonumber \\
= \la \E [ \langle f   \xi_\la(\bX,T; \Po_\la)\rangle^2 ] 
\nonumber \\ 
+ \la^2 
\E  [\langle f,\xi_\la(\bX,T; \Po_\la^{\bX',T'})\rangle 
  \langle f, \xi_\la(\bX',T'; \Po_\la^{\bX,T})\rangle ].
\lbl{02160}
\eea
Combining (\ref{0216a0}) and (\ref{02160}), we have
\bea
\la^{-2}\Var \langle f,\mu_\la^\xi\rangle 
= \la^{-1} \E [ \langle f, \xi_\la(\bX,T;\Po_\la)\rangle^2] 
\nonumber \\
+ \E[ \langle  f , 
\xi_\la(\bX,T;\Po_\la^{\bX',T'})\rangle 
\langle f,
\xi_\la(\bX',T';\Po_\la^{\bX,T}) \rangle ]
\nonumber \\
- ( \E [ \langle f,\xi_\la (\bX,T;\Po_\la) \rangle ] )^2.
\lbl{0826d}
\eea
Also, by  similar arguments,
\bea
n^{-1} \E \langle f, \nu_{\la,n}^\xi \rangle =
 \E \langle f, \xi_{\la}(\bX,T;\X_{n-1}) \rangle 
\lbl{050421}
\eea
and
\bea
n^{-2} \Var \langle f, \nu_{\lambda,n}^\xi \rangle
= n^{-1} \E [ \langle f, \xi_\la(\bX,T;\X_{n-1})\rangle^2]
\nonumber \\
+ \left( \frac{n-1}{n}\right) 
 \E [
\langle f, \xi_\la(\bX,T;\X_{n-2}^{\bX',T'} )\rangle
\langle f, \xi_\la(\bX',T';\X_{n-2}^{\bX,T} )\rangle ]
\nonumber \\
- ( \E[ \langle f,\xi_\la (\bX,T;\X_{n-1} ) \rangle] )^2.
\lbl{050331d}
\eea
Recall that by definition,
 $\xi_\la(\bx,t;\X,\R^d)=0$ for $\bx \in \R^d \setminus \Gamma_\la$,
and $(\skala,\la \geq 1)$ is a given nondecreasing
 family of Borel subsets of
$\R^d$ with limit set $\Gamma$ having Lebesgue-null boundary; in the  simplest
case $\skala = \R^d$ for all $\la$. \\ 

\noindent
{\em Proof of Theorem \ref{thlln1}.}
First we prove (i) for the case $q=2$.
Assume that \eq{mom} and \eq{mom2} hold for some $p >2$.
Let $\H_{\kappa(\bX)}$ denote a Cox point 
process in $\R^d \times \MM$, whose distribution, given $\bX=\bx$,
is that of $\H_{\kappa(\bx)}$. 
Set 
 $J := f(\bX )\xi_\infty^{\bX,T}(\H_{\ka(\bX)},\R^d) \1_{\Gamma}(\bX)$,
and let
$J'$ be an independent copy of $J$.
%of $f(\bX )\xi_\infty^{\bX,T}(\H_{\ka(\bX)},\R^d) \1_{\Gamma}(\bx)$.

For any bounded continuous test function $h$ on $\R$,
by \eq{050401g} from Lemma \ref{couplem0c}, 
as $\la \to \infty$
 we have
$\E[h(\langle f, \xi_\la(\bX,T;\Po_\la)\rangle )|\bX ] \to \E[J|\bX]$,
almost surely.  Hence, 
$\E[h(\langle f, \xi_\la(\bX,T;\Po_\la)\rangle )] \to \E[h(J)]$,
so that
%and moreover
% as $\la \to \infty$ we have
\bea
\langle f, \xi_\la(\bX,T;\Po_\la ) \rangle \toD J.
\lbl{050321b}
\eea
Similarly, using \eq{050401h} we obtain
%\eea
%and
\bea
 \langle  f,  
\xi_\la(\bX,T;\Po_\la^{\bX',T'}) \rangle 
\langle f,
\xi_\la(\bX',T';\Po_\la^{\bX,T}) \rangle 
\toD  J'J
\lbl{0826c}
\eea
Also, by \eq{mom} and \eq{mom2} the variables in the 
left side of \eq{050321b} and in the left hand 
side of \eq{0826c} are uniformly integrable so
we have convergence of means in both cases. Also, 
\eq{mom} shows that  the first term in the right side of
\eq{0826d} tends to zero. Hence we find that the expression
 \eq{0826d} tends to zero.
Moreover, by \eq{0216a0} and the convergence of expectations corresponding
to \eq{050321b}, $\la^{-1} \E\langle f, \mu_\la^\xi \rangle$
 tends  to $\E[J]$, and this gives us \eq{LLNpo} for
$q=2$, under the assumptions of part (i) of Theorem \ref{thlln1}.

Now consider the case $q=1$.
 Assume \eq{mom} and \eq{mom2} hold for some
$p >1$.
 We use a
 truncation argument.
% $t^+ = \max(t,0) $ and $t^- = (-t)^+$ as usual.
Define positive and negative parts $f^+$ and $f^-$ of the test function $f$ by
 $f^+(\bx) := \max(f(\bx),0) $ and $f^-(\bx):=\max(-f(\bx),0)$.
%Let $\xi^+(\bx,t;\X,A)$ be the positive
%part of the signed measure $\xi(\bx,t;\X,\cdot)$ restricted to $A$,
%and let
%$\xi^-(\bx,t;\X,A)$ be the negative part.
 For $K>0$, let 
% define 
%$\xi^{+,K}$ and $\xi^{-,K}$ to be a  
$\xi^{K}$ 
%and $\xi^{-,K}$ to be a  
 %truncated version of  the  positive measures $\xi^+$ and
 truncated version of  the measure $\xi^*$, 
%and $\xi^-$, respectively,  
defined
% under Assumption A1   
by 
\bean
%\xi^{K,+} (\bx,t;\X,A)
% := \min(   \xi^+(\bx,t;\X,\R^d), K) {\bf 1}_{A}(\bx); \\
\xi^{K} (\bx,t;\X,A) := \min( \xi(\bx,t;\X,\R^d), K) {\bf 1}_{A}(\bx).
% \\
%\xi^{K,-} (\bx,t;\X,A) :=
% \min(   \xi^-(\bx,t;\X,\R^d), K) {\bf 1}_{A}(\bx)
%~~~{\rm under~ A1}
\eean
%and under Assumption A2 by
%\bean
%\xi^{K} (\bx,t;\X,d\by) := \xi(\bx,\X,d\by) {\bf 1}_{B_{S}(\bx)} (\by) 
%%\xi^{K,+} (\bx,t;\X,d\by)
%% := \xi^+(\bx,\X,d\by) {\bf 1}_{B_{S^+}(\bx)} (\by) 
%%\\
%%\xi^{K,-} (\bx,t;\X,d\by) := \xi^-(\bx,\X,d\by) {\bf 1}_{B_{S^-}(\bx)}(\by)
%~~~{\rm under~ A2}
%% \max(  (f(\bx) \xi(x,\bX) )^-, K)    
%\eean
%where
% $S := \sup \{r: \xi^+(\bx,\X,B_r(\bx) ) \leq K \}$.
%%
%% $S^+ := \sup \{r: \xi^+(\bx,\X,B_r(\bx) ) \leq K \}$ and
%% $S^- := \sup \{r: \xi^-(\bx,\X,B_r(\bx) ) \leq K \}$.
%%Given $f$, we have a uniform bound on $\langle f^+,\xi^{K,+}$.
%%Then $\xi^{K,+}$ is a positive measure with
Then $\xi^{K}$ has 
%a positive measure with
 total measure bounded by
$K$.
By 
the case $q=2$ 
%which we have
 established above, 
%we have
\bea
%\la^{-1} \langle  f^+, \mu_\la^{\xi^{K}} \rangle \inL 
\la^{-1} \langle  f^+, \mu_\la^{\xi^{K}} \rangle \inL 
\E \left[ f^+(\bX)(\xi^{K})^{\bX,T}_\infty 
(\H_{\ka(\bX)},\R^d) \1_{\Gamma}(\bX)\right].
\lbl{trunc}
\eea
Then $\langle f^+ ,\xi^*(\bx,t;\X)\rangle = \lim_{K \to \infty} 
%Then $\langle f^+ ,\xi^+(\bx,t;\X)\rangle = \lim_{K \to \infty} 
%\langle f, \xi^{K,+}(\bx,t;\X)\rangle $.
\langle f^+, \xi^{K}(\bx,t;\X)\rangle $.
Using \eq{0216a0}, we have
\bean
0 \leq \E  [
\la^{-1} \langle  f^+, \mu_\la^{\xi^*}\rangle -  
\la^{-1} \langle  f^+, \mu_\la^{\xi^{K}}\rangle ] \\
= \E[ \langle f^+,\xi^*_\la(\bX,T;\Po_\la) - \xi_\la^{K}
(\bX,T;\Po_\la)\rangle]
\\
\leq \|f\|_\infty  \E [ 
\xi_\la(\bX,T;\Po_\la,\R^d)
{\bf 1}\{ \xi_\la(\bX,T;\Po_\la,\R^d) > K\} ]
\eean
which tends to zero as $K \to \infty$, uniformly in $\la$,
because the moments condition \eq{mom}, $p >1$, implies
that the random variables
% $\xi^+_\la(\bX,T;\Po_\la,\R^d)$ are uniformly integrable.
 $\xi_\la(\bX,T;\Po_\la,\R^d)$ are uniformly integrable.
Also, by monotone convergence, as $K \to \infty$ 
the right side of \eq{trunc}
converges to  
%$\E [ f^+(\bX) (\xi^{+})^{\bX,T}_\infty 
$\E [ f^+(\bX) \xi^{\bX,T}_\infty 
(\H_{\ka(\bX)},\R^d) \1_{\Gamma}(\bX)]$.
Hence, taking $K \to \infty$ in \eq{trunc} yields
\bean
\la^{-1} 
%\langle  f^+, \mu_\la^{\xi^{+}} \rangle 
\langle  f^+, \mu_\la^{\xi^*} \rangle 
\inL 
%\E [ f^+(\bX) (\xi^{+})^{\bX,T}_\infty (\H_{\ka(\bX)},\R^d) \1_{\Gamma}(\bX)],
\E [ f^+(\bX) \xi^{\bX,T}_\infty (\H_{\ka(\bX)},\R^d) \1_{\Gamma}(\bX)],
\eean
and a similar argument yields the equivalent statement
with 
%$\xi^+$ replaced by  $\xi^{-}$, and/or with 
$f^+$ replaced by
$f^-$. Combining these, and using linearity, we obtain 
\bea
\la^{-1} \langle f, \mu_\la^{\xi^*} \rangle 
= \la^{-1} ( \langle  f^+, \mu_\la^{\xi^*} \rangle 
 - \langle  f^-, \mu_\la^{\xi^*} \rangle 
%= \la^{-1} ( \langle  f^+, \mu_\la^{\xi^{+}} \rangle 
%- \langle  f^+, \mu_\la^{\xi^{-}} \rangle 
% - \langle  f^-, \mu_\la^{\xi^{+}} \rangle 
% + \langle  f^-, \mu_\la^{\xi^{-}} \rangle 
)
\nonumber \\
\inL
\E[( f^+(\bX) \xi^{\bX,T}_\infty (\H_{\ka(\bX)},\R^d)
-f^-(\bX) \xi^{\bX,T}_\infty (\H_{\ka(\bX)},\R^d)
%\E[( f^+(\bX) (\xi^{+})^{\bX,T}_\infty (\H_{\ka(\bX)},\R^d)
%- 
%f^-(\bX) (\xi^{+})^{\bX,T}_\infty (\H_{\ka(\bX)},\R^d)
%\\
%-f^-(\bX) (\xi^{-})^{\bX,T}_\infty (\H_{\ka(\bX)},\R^d)
%+ f^-(\bX) (\xi^{-})^{\bX,T}_\infty (\H_{\ka(\bX)},\R^d)
) \1_{\Gamma}(\bX)]
\nonumber \\
= 
\E [ \xi^{\bX,T}_\infty  ( \H_{\ka(\bX)},\R^d) f(\bX)\1_{\Gamma}(\bX)].
\lbl{050408a}
\eea
This gives us \eq{LLNpo} for $q=1$ when Assumption A1 hols.

Now suppose A2 or A3 holds. By Palm theory, analogously to
\eq{0216a0} we have
\bea
\E 
\la^{-1} \sum_{\bx \in \Po_\la} |\langle f,\xi_\la(\bx;\Po_\la)
- \xi_\la^*(\bx;\Po_\la) \rangle | 
~~~~~
\nonumber \\
= \int_{\Gamma}
 \ka(d\bx) \E|\langle
f, \xi_\la(\bx;\Po_\la) - \xi_\la^*(\bx;\Po_\la) \rangle 
|,
\lbl{050408}
\eea
and for almost every $\bx \in \Gamma_1$, by \eq{050401d} and \eq{mom},
the integrand tends to zero and is bounded so that we have \eq{LLN2po}.
Combining this with \eq{050408a} gives us \eq{LLNpo} for $q=1$ when
Assumption A2 or A3 holds, completing the proof of part (i).

%wh completing the proof of part (ii).

Next we turn to part (ii) with $q=2$.
Assume \eq{momlnb} holds for some $p >2$.
By Lemma \ref{couplem0c}, 
as $n \to \infty$ we have
\bea
\langle f, \xi_{\la(n)}(\bX,T; \X_{n-1})\rangle
\toD J;
\lbl{050331e}
\\
\langle f, \xi_{\la(n)}(\bX,T; \X_{n-2}^{\bX',T'})\rangle
\langle f, \xi_{\la(n)}(\bX',T'; \X_{n-2}^{\bX,T})\rangle
\toD J'J,
\lbl{050331f}
\eea
Also, by \eq{momlnb} and the Cauchy-Schwarz inequality,
 the variables in the 
left side of \eq{050331e} and in the left  
side of \eq{050331f} are uniformly integrable so
we have convergence of means in both cases. Likewise, 
\eq{momlnb} shows that  the first term in the right side of
\eq{050331d} tends to zero. Hence we find that the expression
 \eq{050331d} (with $\la = \la(n)$) tends to zero.
Also, by \eq{050421} and \eq{050331e}, $n^{-1} \E \langle f, \nu_{\la(n),n}
\rangle  \to \E[J]$, and this gives us \eq{LLNbi} for
$q=2$.
% under the assumptions of part (iii) of Theorem \ref{thlln1}. 

The case $q=1$ of part (ii) is deduced from the case $q=2$
%Part (iv)  is deduced from part (iii) 
in the same manner as in 
%part (ii) was deduced from
 part (i).
 $\qed$

\section{Proof of Theorems \ref{thm1} and \ref{thm2}}
\allco

Note that by definition,
 $\xi_\la(\bx,t;\X,\R^d)=0$ for $\bx \in \R^d \setminus \Gamma$. 
In the sequel, we fix a test function $f \in B(\R^d)$.
%and for $\la >0$, we set  $f_\la (\bx) := f(\bx) \1_{\skala}(\bx)$, and
Set
\bea
\alpha_\la : =  \int_{\Gamma}
%f_\la(\bx)^2
 \E[ \langle f,\xi_\la(\bx,T; \Po_\la)\rangle^2 ]\ka(\bx) d\bx
\lbl{alphadef}
\eea
and
\bea
\beta_\la :=
 \int_{\Gamma} \int_{\R^d}
(\E  [\langle f, \xi_\la(\bx,T; \Po_\la^{\bx+ \la^{-1/d} \bz,T'}) \rangle
 \langle f, \xi_\la(\bx+ \la^{1/d} \bz,T'; \Po_\la^{\bx,T})\rangle ]
\nonumber \\
 - \E[\langle f, \xi_\la(\bx,T; \Po_\la)\rangle] 
\E[\langle f, \xi_\la(\bx+ \la^{-1/d} \bz,T; \Po_\la)\rangle]) 
\nonumber \\
\times
\ka(\bx) \ka(\bx+ \la^{-1/d}\bz) d\bz d\bx;
\lbl{0420a}
\eea
\begin{lemm}
\lbl{abclem}
For $\la >0$, it is the case that
\bea
 \Var(\langle f,\mu_\la^\xi \rangle )=\la ( \alpha_\la + \beta_\la).
\lbl{abc}
\eea
\end{lemm}
{\em Proof.}
The first term in the right hand side of \eq{0826d}
 equals $\la^{-1}\alpha_\la$.
Thus, \eq{0826d} yields
\bea
\la^{-1} \Var(\langle f, \mu_\la^\xi\rangle ) - \alpha_\la 
 = \la \int_{\Gamma} \int_{\R^d}
(\E  [\langle f, \xi_\la(\bx,T; \Po_\la^{\by,T'}) \rangle
\langle f, \xi_\la(\by,T'; \Po_\la^{\bx,T}) \rangle ]
\nonumber
 \\
-  \E[\langle f, \xi_\la(\bx,T; \Po_\la)\rangle]
 \E[\langle f, \xi_\la(\by,T; \Po_\la)\rangle] )  
\ka(\bx) \ka(\by) d\by d\bx,
%- \gamma_\la
\lbl{0826b}
\eea
and the change of variables $\by=\bx+\la^{-1/d }\bz$ shows that
this equals $\beta_\la$ as given by 
(\ref{0420a}). $\qed$ \\

%Our next task is to establish the limiting behaviour of
Lemmas \ref{lemweak1} and \ref{lemweak2} establish limits
 in distribution
for the variables inside the expectations in 
the integrands in the expressions
 (\ref{alphadef}) and (\ref{0420a})   
for $\alpha_\la$ and $\beta_\la$.
%these integrands.
To prove Theorem \ref{thm1},
we need to take these limits outside the expectations and also
outside the integrals, which we shall do by a domination argument.
 It is in this step that we use
the condition of stabilization with respect to non-homogeneous
Poisson processes (Definition \ref{stab}), via the following lemma,
which  
is an estimate showing that the integrand 
in the definition (\ref{0420a}) of $\beta_\la$
is small for large $|\bz|$, uniformly in $\la$. 
To ease notation,
for $\bx \in \R^d$, $\bz \in \R^d$ and $\la >0$,
 we define random variables $X =X_{\bx,\bz,\la}$, 
$Z = Z_{\bx,\bz,\la}$, $X'=X'_{\bx,\bz,\la}$ and $Z' = Z'_{\bx,\bz,\la}$,
by
\bea
X:= \langle f, \xi_\la(\bx,T; \Po_\la^{\bx+ \la^{-1/d}\bz,T'})\rangle, ~~
Z:= \langle f,\xi_\la(\bx+ \la^{-1/d}\bz,T'; \Po_\la^{\bx,T}) \rangle,
\lbl{XZdef}
\\
X':= \langle f, \xi_\la(\bx,T; \Po_\la) \rangle , ~~
Z' := \langle f, \xi_\la(\bx+ \la^{-1/d}\bz,T'; \Po_\la) \rangle.
\lbl{XdZddef}
\eea
Similarly, we define random variables
 $X^* =X^*_{\bx,\bz,\la}$, 
$Z^* = Z^*_{\bx,\bz,\la}$, $X^{*'}=X^{*'}_{\bx,\bz,\la}$
 and $Z^{*'} = Z^{*'}_{\bx,\bz,\la}$,
by
\bea
X^*:=  \xi_\la(\bx,T; \Po_\la^{\bx+ \la^{-1/d}\bz,T'},\R^d), ~~
Z^*:= \xi_\la(\bx+ \la^{-1/d}\bz,T'; \Po_\la^{\bx,T},\R^d) ,
\lbl{XZsdef}
\\
X^{*'}:= \xi_\la(\bx,T; \Po_\la,\R^d) , ~~
Z^{*'} :=  \xi_\la(\bx+ \la^{-1/d}\bz,T'; \Po_\la,\R^d) .
\lbl{XdZdsdef}
\eea
We write $a\wedge b$ for $\min(a,b)$ in the sequel.
\begin{lemm}
\lbl{XZlem}
Suppose that $\xi$ satisfies \eq{mom} and \eq{mom2} for some $p>2$,
and is power-law stabilizing for $\kappa$ of order $q$ for some
$q > d p/(p-2)$. Then
there is a constant $C_1$, independent
of $\la$, such that for all $\la \geq 1$, 
$\bx \in \ska$ and $\bz \in \R^d$,
\bea
| \E[X_{\bx,\bz,\la}Z_{\bx,\bz,\la}]- \E[X'_{\bx,\bz,\la}]\E[Z'_{\bx,\bz,\la}]
| 
%\leq C_1 \exp(-|z|/C_1)
\leq C_1 ( |\bz|^{-d -(1/C_1)} \wedge 1) ; 
\lbl{0426k} \\
| \E[X^*_{\bx,\bz,\la}Z^*_{\bx,\bz,\la}]-
 \E[X^{*'}_{\bx,\bz,\la}]\E[Z^{*'}_{\bx,\bz,\la}] | 
%\leq C_1 \exp(-|z|/C_1)
\leq C_1 (|\bz|^{-d -(1/C_1)}\wedge 1)  .
\lbl{0426k2}
\eea
\end{lemm}
{\em Proof.}
Let $X:= X_{\bx,\bz,\la}$ and $Z:= Z_{\bx,\bz,\la}$.
Let $\tilde{X}= X{\bf 1}_{\{R_\la(\bx,T) \leq |\bz|/3\}}$ 
and let $\tilde{Z}= Z{\bf 1}_{\{R_\la(\bx+ \la^{-1/d} \bz,T')\leq |\bz|/3\}}$.
Then  $\tilde{X}$ and $\tilde{Z}$ are independent so that
$\E[\tilde{X}\tilde{Z}]= \E[\tilde{X}]\E[\tilde{Z}]$, and
\bea
\E[XZ]= \E[\tilde{X}] \E[\tilde{Z} ] + \E[\tilde{X}(Z - \tilde{Z})]
+ \E[(X- \tilde{X})Z]
\lbl{0424a}
\eea
while
\bea
\E[X']\E[Z']= \E[\tilde{X}]\E[\tilde{Z} ] + \E[\tilde{X}]\E[Z' - \tilde{Z}]
+ \E[X'- \tilde{X}]\E[Z'].
\lbl{0424b}
\eea
By 
 (\ref{mom2}) and H\"older's inequality,
and the assumed power-law stabilization of order $q > dp/(p-2)$,
there is a constant $C_2$ such that 
\bea
\E[ |(X-\tilde{X})Z|]
 = \E[ |XZ| 
{\bf 1}_{\{ R_\la(\bx) > |\bz|/3\}}  ]
\nonumber \\
\leq (\E[|X|^p])^{1/p}(\E[|Z|^p])^{1/p} 
(P[ R_\la(\bx,T) > |\bz|/3 ])^{1-(2/p)}
\nonumber \\
\leq C_2 ( |\bz|^{-d -(1/C_2)} \wedge 1)
\lbl{0424c}
\eea
and likewise
\bea
\E[ |X(Z-\tilde{Z})|] \leq
 C_2 ( |\bz|^{-d -(1/C_2)} \wedge 1).
\lbl{0424d}
\eea
By a  similar argument using (\ref{mom}), there is a constant $C_3$ such that
\bea
\max(\E[ |X'-\tilde{X}|], 
\E[ |Z'-\tilde{Z}|] ) < C_3 (|\bz|^{-d-(1/C_3)} \wedge 1).
\lbl{0424e}
\eea 
Subtracting (\ref{0424b}) from (\ref{0424a}) and using 
(\ref{0424c}), (\ref{0424d}) and (\ref{0424e})
along with the Cauchy-Schwarz inequality,
 we may deduce that there is a constant $C_1$, independent
of $\la$, such that for all $\la \geq 1$, (\ref{0426k}) holds. 
The argument for \eq{0426k2} is similar 
$\qed$ \\

%We now prove Theorem \ref{thm2}, 
%which is a corollary of Theorem \ref{thm1} together with 
%Lemma \ref{mainlem}

{\em Proof of Theorem \ref{thm1}.}
Let $f \in B(\R^d)$.
Let $\bx \in \Gamma_0 $ (as defined at \eq{gam1def}).
% the th be a continuity point of $\ka$.
Assume also, either that A1 holds, 
%or A2 holds and $\bx$ is a Lebesgue point of $f$,
 or that A3 holds and $\bx$
is a continuity point of $f$; also assume A4 or A5 holds.
By the case  $\bz=0$ of
(\ref{0216c}) when A1 holds, 
and the case $\bz=0$ of \eq{50209f} when A3 holds,
we have
\bea
\langle f, \xi_\la(\bx,T; \Po_\la )\rangle \tod
  f(\bx)  \xi^{\bx,T}_\infty(\H_{\kappa(\bx)},\R^d).
\lbl{050321c}
\eea
By (\ref{mom}), 
$\{\langle f, \xi_\la(\bx; \Po_\la )\rangle^2:
 \la \geq 1\}$ are uniformly integrable,
and hence the convergence in distribution \eq{050321c} 
 extends to convergence of second moments to a limit
which is bounded by
 (\ref{mom}) and 
 Fatou's Lemma. Hence, by the dominated convergence theorem,
$\alpha_\la$ given by (\ref{alphadef})
satisfies
\bea
\lim_{\la \to \infty} \alpha_\la  = 
\int_\Gamma 
( f(\bx)
\E [
  \xi^{\bx,T}_\infty(\H_{\kappa(\bx)},\R^d) ] )^2
 \kappa(\bx)d\bx < \infty.
\lbl{Aconv}
\eea

Next we show convergence of the expression $\beta_\la$
given by (\ref{0420a}). To this end, set 
\bean
g_\la(\bx,\bz) : = 
(\E[ X_{\bx,\bz,\la} Z_{\bx,\bz,\la}]- \E[X'_{\bx,\bz,\la} ]
\E[Z'_{\bx,\bz,\la}])
\ka(\bx+\la^{-1/d}\bz);
\\
g^*_\la(\bx,\bz) : = 
(\E[ X^*_{\bx,\bz,\la} Z^*_{\bx,\bz,\la}]- \E[X^{*'}_{\bx,\bz,\la} ]\E[Z^{*'}_{\bx,\bz,\la}])
\ka(\bx+\la^{-1/d}\bz).
\eean
%Suppose $\bx \R^d$ is  a continuity point of $\ka$.
%Then 
Suppose A3 holds. Then
for almost all $\bx\in \Gamma_0$ and 
all $\bz \in \R^d$,
 with $\ka(\bx)>0$
and $\ka $ continuous at $\bx$,
by Lemma \ref{lemweak2}
we have
 as $\la \to \infty$  that
\bea
X_{\bx,\bz,\la}Z_{\bx,\bz,\la} \tod
 f(\bx)^2
\xi_\infty^{\bx,T}( \H^{\bz,T'}_{\kappa(\bx)}, \R^d )
\xi_\infty^{\bx,T'}(-\bz+  \H^{\0,T}_{\kappa(\bx)}, \R^d ),
%\1_{\Gamma}(\bx)
\lbl{050321d}
\eea
and the variables in the left
side of \eq{050321d}
are  uniformly  integrable by (\ref{mom2}),
so that \eq{050321d} extends to  convergence of expectations.
Likewise by (\ref{mom}) and (\ref{50209f}), both 
   $\E[X^{'}_{\bx,\bz,\la}]$ and $\E[Z^{'}_{\bx,\bz,\la}]$ 
converge to
$f(\bx) \E[\xi_\infty^{\bx,T}( \H^\bz_{\kappa(\bx)},\R^d )]$, so we
 have under A3 that
\bea
\lim_{ \la \to \infty} (
g_{\la}(\bx,\bz) ) = g_\infty(\bx,\bz) , ~~~{\rm a.e.} ~(\bx,\bz)
\in \Gamma_0 \times \R^d,
\lbl{glim}
\eea
 with 
\bean
g_\infty(\bx,\bz) := f(\bx)^2
%\1_{\Gamma}(\bx)
\kappa(\bx)
 (
\E[ \xi_\infty^{\bx,T}( \H^{\bz,T'}_{\kappa(\bx)},\R^d )
\xi_\infty^{\bx,T'}(-\bz+  \H^{\0,T}_{\kappa(\bx)} ,\R^d)]
\\
- \E[ \xi_\infty^{\bx,T}( \H_{\kappa(\bx)} ,\R^d)]^2
) 
.
\eean
By our  assumptions on $\kappa (\cdot)$ and $\Gamma$,
\bea
\int_{\Gamma \setminus \Gamma_0} \ka(\bx) d\bx = 0.
\lbl{gambd}
\eea
By Lemma \ref{XZlem}, and the assumption that 
$\kappa$ is bounded, 
%and (\ref{glim})
 there is a constant $C$ such that
$
|g_\la(\bx,\bz) |
 \leq C ( |\bz|^{-d-(1/C)} \wedge 1), 
%~~~
%\lbl{expbd}
$
for almost every $(\bx,\bz) \in \Gamma_0 \times \R^d$, and $\la \geq 1$.
%\bea
%\eea
% for $1 \leq \la \leq \infty. 
Hence, by \eq{0420a}, \eq{glim}, 
%\eq{expbd}, 
\eq{gambd}, 
and the dominated convergence theorem,
we have
\bean
\beta_\la = \int_{\Gamma} \int _{\R^d} g_\la(\bx,\bz) \ka(\bx) d\bz d\bx 
\to 
 \int_{\Gamma} \int _{\R^d} g_\infty(\bx,\bz) \ka(\bx) d\bz d\bx  < \infty. 
\eean
Combined with \eq{Aconv} and \eq{abc}, this gives us \eq{BYeq} under A3.

Now suppose instead of A3 that A1 holds. Then by \eq{0426m}, 
for almost all $(\bx,\bz) \in \Gamma_0 \times \R^d$
% with $\ka(\bx)>0$ and $\ka $ continuous at $\bx$,
we have
 as $\la \to \infty$  that
\bea
X^*_{\bx,\bz,\la}Z^*_{\bx,\bz,\la} \tod
\xi_\infty^{\bx,T}( \H^{\bz,T'}_{\kappa(\bx)}, \R^d )
\xi_\infty^{\bx,T'}(-\bz+  \H^{\0,T}_{\kappa(\bx)}, \R^d ),
%\1_{\Gamma}(\bx),
\lbl{050321d2}
\eea
%by Lemma \ref{lemweak2}, 
and the variables in the left
side of \eq{050321d2}
are  uniformly  integrable by (\ref{mom2}),
so that \eq{050321d2} extends to  convergence of expectations.
Likewise by (\ref{mom}) and (\ref{0216c}), both 
   $\E[X^{*'}_{\bx,\bz,\la}]$ and $\E[Z^{*'}_{\bx,\bz,\la}]$ 
converge to
$\E[\xi_\infty^{\bx,T}( \H^\bz_{\kappa(\bx)},\R^d )]$, so we
 have under A1 that
\bea
\lim_{ \la \to \infty} (
g_{\la}^*(\bx,\bz) ) = g^*_\infty(\bx,\bz) , ~~~{\rm a.e.} ~(\bx,\bz)
\in \Gamma_0 \times \R^d
\lbl{glim2}
\eea
 with 
\bean
g^*_\infty(\bx,\bz) := 
%\1_{\Gamma}(\bx) 
\kappa(\bx) (
\E[ \xi_\infty^{\bx,T}( \H^{\bz,T'}_{\kappa(\bx)},\R^d )
\xi_\infty^{\bx,T'}(-\bz+  \H^{\0,T}_{\kappa(\bx)} ,\R^d)]
%\\
- \E[ \xi_\infty^{\bx,T}( \H_{\kappa(\bx)} ,\R^d)]^2
).
\eean
By Lemma \ref{XZlem}, the assumption that 
$\kappa$ is bounded, and (\ref{glim2})
 there is a constant $C$ such that
for almost every $(\bx,\bz)$ with $\ka(\bx)>0$,
\bea
| g_\la^*(\bx,\bz) | 
%: 1 \leq \la \leq \infty 
 \leq C ( |\bz|^{-d-(1/C)}\wedge 1),
 ~~~
 1 \leq \la \leq \infty. 
\lbl{expbd2}
\eea
If $\bx\in \R^d$ is a Lebesgue point
of $f$,
% (see e.g. \cite{Pe}),
 then
 for any $K>0$,
by \eq{0426k2} we have
$$
\lim_{\la \to \infty} \int_{B_K} g^*_\la(\bx,\bz)(f(\bx+\la^{-1/d}\bz) 
- f(\bx)) d\bz  =0,
$$
and combining this with (\ref{glim2})  and the dominated convergence
theorem gives us
\bea
\lim_{\la \to \infty} \int_{B_K} g^*_\la(\bx,\bz)f(\bx+\la^{-1/d}\bz)  d\bz
= \int_{B_K}
g^*_\infty(\bx,\bz) f(\bx) d\bz.
\lbl{0429a}
\eea
On the other hand, by (\ref{expbd2}) and the assumption that
$f$ is bounded, we have
$$
\lim_{K \to \infty} \limsup_{\la \to \infty}
 \int_{\R^d \setminus B_K}
|g^*_\la(\bx,\bz)
f(\bx+\la^{-1/d}\bz) 
 -g^*_\infty(\bx,\bz)f(\bx)| 
d\bz
=0
$$
and combining this with \eq{0429a} we have
$$
\lim_{\la \to \infty} 
 \int_{\R^d}
g_\la^*(\bx,\bz) f(\bx+\la^{-1/d}\bz) d\bz
= 
\int_{\R^d} g^*_\infty(\bx,\bz)f(\bx) 
d\bz.
$$
By the Lebesgue density theorem,
almost all every $\bx \in \Gamma_0$ is  a Lebesgue point
of $f$.
Hence, under A1,
by \eq{0420a}, \eq{gambd},
 the dominated convergence theorem,
\bean
\beta_\la = \int_{\Gamma} f(\bx)  \int_{\R^d}
f(\bx + \la^{-1/d} \bz) g_\la^*(\bx,\bz) \ka(\bx) d\bz
d\bx
%\\
%\lim_{\la \to \infty} \beta_\la^* =
\to
\int_{\Gamma}f(\bx)^2
\int_{\R^d} g^*_\infty(\bx,\bz) \ka(\bx) d\bz d\bx, 
%\\
%=  \int_{\Gamma}  \int_{\R^d} ( 
%\E[ \xi_{\infty}^{\bx,T} (H_{\ka(\bx)}^{\bz,T'},\R^d)
%\xi_{\infty}^{\bx,T'} (-\bz + H_{\ka(\bx)}^{\0,T},\R^d) ]
%- \E[ \xi_{\infty}^{\bx,T} (H_{\ka(\bx)},\R^d)]^2) 
%d\bz
%\\
%\times
%f(\bx)^2
%\ka(\bx)^2 
% d\bx,
\eean
and combined with (\ref{Aconv}) and (\ref{abc}), 
%and Lemma \ref{bstarlem},
this gives us (\ref{BYeq}) as required. $\qed$ \\

For the proof of Theorem \ref{thm2}
 (central limit theorem for random measures), 
we shall use results on
  normal approximation for
$\langle f, \overline{\mu}^{\xi}_{\la } \rangle $, suitably scaled.
In the case of point measures, these
 were proved by Stein's method in \cite{PY5},
and the method carries through to more general measures.
Let $\Phi$ denote the standard normal
distribution function, and let $\NN(0,\sigma^2)$ denote the
normal distribution with mean $0$ and variance $\sigma^2$
(if $\sigma^2>0$) or the unit point mass at $0$ if $\sigma^2 =0$.

\begin{lemm}\lbl{mainlem}
Suppose that $\ka $ has bounded support and $\|\ka\|_\infty < \infty$.
 Suppose that $\xi$ is exponentially stabilizing and
satisfies the moments condition $(\ref{mom})$ for some $p >2$.
Let $f \in B(\R^d)$, and $q \in (2,3]$ with $q <p$.
%and put $T_\la := \langle f, \mu^{\xi}_{\la } \rangle$.
  There exists
 a finite constant $C$ depending on $d, \xi$,  $\ka$, $	q$  and $f$,
such that for all $\la > 1$
 \be
\lbl{rates2} \sup_{t \in \R}  \left| P \left[ {
% T_\la - \E T_\la
 \langle f, \overline{\mu}^{\xi}_{\la } \rangle
\over ( \Var
 \langle f, \mu^{\xi}_{\la } \rangle
 )^{1/2}  } \leq t \right] - \Phi(t) \right|
\leq C ( \log \la )^{qd}  \la ( \Var
 %T_{ \la}
 \langle f, \mu^{\xi}_{\la } \rangle
 )^{-q/2}. \ee
\end{lemm}
{\em Proof.} In the case where $\xi = \xi^*$, i.e. $\xi(\bx,t;\X,\cdot)$
is always a point mass at $\bx$, this result is Theorem 2.1 of
\cite{PY5}. If we do not make this assumption on $\xi$,
the proof in \cite{PY5} carries through with little change, except
that the $T_\la$ and $T'_\la$ of Section 4.3 of \cite{PY5} should now 
be defined  (following notation of \cite{PY5}) by
\bean
T_\la = \sum_{i=1}^{V(\la)} \sum_{j=1}^{N_i} \langle f, \xi_\la
(X_{i,j},U_{i,j};\Po_\la) \rangle, 
\\
T'_\la = \sum_{i=1}^{V(\la)} \sum_{j=1}^{N_i} \langle f, \xi_\la
(X_{i,j},U_{i,j};\Po_\la) \rangle {\bf 1}_{E_{ij}}.  ~~~~ \qed
\eean
%and likewise for $T'_\la$. $\qed$\\

\begin{lemm}
\lbl{mainlem2}
Suppose $\|\kappa\|_\infty < \infty$. Suppose
 for some $p >3$ that
  $\xi$ is power-law stabilizing of order $q$
for some  $q > d(150 + 6/p)$, and
satisfies the moments condition $(\ref{mom})$.
Let $f \in B(\R^d)$.
%and put $T_\la := \langle f, \mu^{\xi}_{\la } \rangle$.
Suppose that
$\la^{-1}  \Var
\langle f, {\mu}^{\xi}_{\la } \rangle   
$ converges, as $\la  \to \infty$,
 to a finite limit $\sigma^2$. Then
 %for
 %some $\sigma^2 \geq 0$, 
$
\langle f, \la^{-1/2}\overline{\mu}_{\la }^{\xi}\rangle
$ converges in distribution, as $\la \to \infty$,
to the
 $\NN \left(0, \sigma^2 \right) $ distribution.
\end{lemm}
{\em Proof.} In the case where $\xi = \xi^*$, 
%i.e. $\xi(\bx,t;\X,\cdot)$ is always a point mass at $\bx$, 
this result is Theorem 2.2 of
\cite{PY5}. If we do not make this assumption on $\xi$,
the proof in \cite{PY5} carries through with
the same minor changes
 as indicated for Lemma \ref{mainlem} above. $\qed$\\

\noindent
{\em Proof of Theorem \ref{thm2}.}
%Suppose A6 and the other assumptions of Theorem \ref{thm2} hold.
Suppose $\|\kappa\|_\infty < \infty$
 and $\kappa$ has bounded support.
Suppose $\xi$ is almost everywhere continuous, and is
 $\ka(\bx)-$homogeneously stabilizing at $\bx$ for $\ka$-almost
all $\bx \in \R^d$.
Suppose $\xi$ satisfies either A4 or A5, and satisfies 
A6.

Let $f \in B(\R^d)$, and assume either A1 or A3 holds.
 By Theorem \ref{thm1}, $\la^{-1} \Var \langle f, \mu^{\xi}_{\la } \rangle $
 converges to the finite nonnegative limit 
%\linebreak
$\int_{\Gamma} f(\bx)^2 V^\xi(\bx,\ka(\bx)) \ka(d\bx)$.  If this
limit  is strictly positive, then the right hand side of (\ref{rates2})
tends to zero, so that by  
 Lemma \ref{mainlem}, 
$\la^{-1/2} \langle f, \overline{\mu}_\la^\xi \rangle$
is asymptotically centred normal with variance 
\linebreak
$\int_{\Gamma} f(\bx)^2 V^\xi(\bx,\ka(\bx)) \ka(d\bx)$.  
On the other hand, if this limiting variance is  zero,
then  it is immediate from Chebyshev's inequality that
$\la^{-1/2} \langle f, \overline{\mu}_\la^\xi \rangle$
converges in probability to zero. Hence 
for all $f\in \tilde{B}(\R^d)$, or under A1 for all $f\in B(\R^d)$,
 we obtain
\bea
\la^{-1/2} \langle f, \overline{\mu}_\la^\xi \rangle \tod \NN
\left(0,
\int_{\Gamma} f(\bx)^2 V^\xi(\bx,\ka(\bx)) \ka(\bx) d\bx\right).
\lbl{0518}
\eea
%where $\tod$ denotes convergence in distribution.
If A7 holds instead of A6, we obtain the same conclusion by using
Theorem \ref{thm1} and
Lemma \ref{mainlem2}; note that in A7, since $p>3$ the condition
$q> d(150+6/p)$ ensures that $q> dp/(p-2)$, so that Theorem
\ref{thm1} still applies here.

Now consider an arbitrary finite collection of
test functions $f_1,\ldots,f_J$, each of them in
$\tB(\R^d)$. For arbitrary real constants $b_1,\ldots,b_J$,
application of (\ref{0518}) to $f= \sum_{j=1}^J b_j f_j$
yields
\bean
\sum_{j=1}^J b_j \la^{-1/2} \langle f_j, \overline{\mu}_\la^\xi \rangle \tod \NN
\left(0,
\sum_{j=1}^J \sum_{k=1}^J 
b_j b_k
\int_{\Gamma} f_j(\bx)f_k(\bx) V^\xi(\bx,\ka(\bx)) \kappa(\bx) d\bx \right),
\eean
 so by the Cram\'er-Wold device (see e.g. \cite{Pe}, \cite{Bill}),
the variables $\la^{-1/d} \langle f_j, \overline{\mu}_\la^\xi \rangle, 1 \leq j
\leq J$, are asymptotically centred multivariate normal with covariance
matrix having $(j,k)$th entry 
$\int_\Gamma f_j(\bx)f_k(\bx) V^\xi(\bx,\ka(\bx)) \ka(\bx)d \bx$. This gives us 
the required convergence of random fields for $f_j \in \tB(\R^d)$.
If A1 holds, the same argument gives the required
convergence of random fields for $f_j \in B(\R^d)$.  $\qed$ \\

\section{Extension to the  non-Poisson case}
\lbl{secdepo}
\allco
In this section we prove Theorem \ref{thm3}.
%We represent points in $\R^d \times \MM$ as
%marked points in $\R^d$, to ease notation.
We assume here that all the point processes $\X_n$ are
coupled as described at \eq{050404}, in terms of a sequence
$((\bX_1,T_1),(\bX_2,T_2),\ldots)$ of independent 
random elements of $\R^d \times \MM$ with common distribution
 $\ka\times \mu_\MM$, 
Given $f \in B(\R^d)$ and $\la >0$,
 for each $m \in \N$ we   define
\bea
\ch_{m,\la} := \langle f, \nu_{m+1,\la}
%,f \rangle - \langle
-  \nu_{m,\la} \rangle 
= Y_{m+1,\la} + \sum_{i=1}^m \Delta_{i,m,\la} ,
\lbl{Rdef}
\eea
where we set
\bean
Y_{m+1,\la} := 
\langle f, \xi_\la(\bX_{m+1},T_{m+1} ; \X_m) \rangle;
\\
\Delta_{i,m,\la} :=
\langle f, 
 \xi_\la(\bX_i,T_i;\X_{m+1} ) - \xi_\la(\bX_i,T_i; \X_m)  \rangle.    
\eean
In this section we shall use the standard
 notation $\|X\|_p$
for the $L^p$-norm $\E[|X|^p]^{1/p}$ of a random variable $X$,
where $p\geq 1$.

\begin{lemm}
\lbl{fordepolem}
Suppose 
 $R(\bx,T;\H_{\ka(\bx)}) <\infty$ for $\ka$-almost
all $\bx \in \R^d$, and  
 Assumption A$1$ or A$3$ holds, along with Assumption A$6'$ or A$7'$.
Let $(h(\la))_{\la \geq 1}$ satisfy $h(\la)/\la \to 0$ as $\la \to \infty$.
Then with $\delta(\bx,\la)$ defined at \eq{deltadef},
\bean
\lim_{\la \to \infty}  \sup_{\la-h(\la)  \leq \ell < m  \leq \la+ h(\la) }
\left| \E \ch_{\ell,\la} \ch_{m,\la} - 
\left(\int_\Gamma \ka(\bx) f(\bx) \delta(\bx,\ka(\bx))d\bx \right)^2\right| =0.
\lbl{depoE12}
\eean
%\end{equation}
%and
%\begin{equation}
%\lim_{\la \to \infty}  \sup_{\la-h(\la)  \leq m   \leq \la+ h(\la) }
%|\E \ch_{m,\la}^2  | < \infty.
%\lbl{depoE2}
%\end{equation}
\end{lemm}
{\em Proof.}
We suppress mention of marks in the proof, writing simply $\bX$ 
for $(\bX,T)$.
 Suppose $(\ell(\la))_{\la \geq 1}$ and $(m(\la))_{\la \geq 1}$ 
  satisfy $\la-h(\la) \leq \ell(\la) < m(\la) \leq \la + h(\la)$. 
We shall show a sequential version of \eq{depoE12}.
To ease notation, we write $\ell$ for $\ell(\la)$, and $m$ for 
$m(\la)$, and  $Y_m$ for $Y_{m,\la}$ and $\Delta_{i,m}$
for $\Delta_{i,m,\la}$.

By using  \eq{Rdef}, expanding  and taking expectations,
we obtain
\bea
\E \ch_\ell \ch_{m}
 = \E \left[ \left( Y_{\ell +1} + \sum_{i=1}^\ell \Delta_{i,\ell} \right) \left(
 Y_{m+1} + \sum_{j=1}^\ell \Delta_{j,m} + \Delta_{\ell+1,m} + \sum_{
j=\ell+2}^m \Delta_{j,m} \right) \right] 
\nonumber \\
= \E[ Y_{\ell+1} Y_{m+1} ] + \ell E[ \Delta_{1,\ell} Y_{m+1} ]
+ \ell E[ Y_{\ell+1} \Delta_{1,m}]  + 
\ell(\ell-1) \E[\Delta_{1,\ell}\Delta_{2,m}]   
\nonumber \\
+ \ell \E[\Delta_{1,\ell} \Delta_{1,m} ]  + \E[Y_{\ell+1}\Delta_{\ell+1,m} ]
+ \ell \E[\Delta_{1,\ell} \Delta_{\ell+1,m} ] 
\nonumber \\
+ (m-\ell -1) \E[ Y_{\ell+1} \Delta_{\ell+2,m} ] 
+ \ell (m-\ell -1) \E[ \Delta_{1,\ell} \Delta_{\ell+2,m} ]. 
\lbl{eqexpand}
\eea
We shall establish the limiting behaviour of 
each term of \eq{eqexpand} in turn.
First we have
\bean
\E[Y_{\ell+1}  Y_{m+1}  ] =
 %\int_{\Gamma_\la} 
 \int 
\ka(\bx) (d \bx)
\int
%\int_{\Gamma_{\la}}
 \ka(\by) (d\by)  \E[
 \langle f,
 \xi_\la(\bx; \X_\ell ) \rangle
 \langle f,
 \xi_\la(\by;\X^\bx_{m-1})\rangle].
\eean  
Here and below, all domains of integration, when not specified,
are $\R^d$.
By % \eq{050401g} from
 Lemma \ref{lemweak3}, 
for almost all $(\bx,\by) \in \Gamma_0 \times \Gamma_0$,
we have
\bea
\langle f,
 \xi_\la(\bx ; \X_\ell ) \rangle
 \langle f,
 \xi_\la(\by;\X_{m-1}^\bx) \rangle 
%\nonumber \\
\tod f(\bx) f(\by)
%\1_{\Gamma} \1_{\Gamma}(\bX)
 \xi_\infty^\bx(\H_{\ka(\bx)},\R^d) \xi_\infty^\by(\tH_{\ka(\by)},\R^d).
%| \xi_\infty^\bx(\H_{\ka(\bx)})| |\xi_\infty^\by(\tH_{\ka(\by)})|,
\lbl{050404j}
\eea
%where $\tH_{\ka(\by})$ is a homogeneous point process in $\R^d$ of
%intensity $\ka(\by)$, independent of $\H_{\ka(\bx)}$. 
By \eq{mom4},
the variables $\langle f ,\xi_\la(\bx ; \X_\ell )
\rangle
\langle
f, \xi_\la(\by;\X_{m-1}^\bx) \rangle$
are uniformly integrable 
 so \eq{050404j} extends to  convergence
of expectations, and also the limit is bounded, so that
setting
\bea
\gamma_1 := \int_{\Gamma} 
f(\bx) \E[ \xi_\infty^{\bx}(\H_{\ka(\bx)} ,\R^d) ]
\ka(\bx)(d\bx),
\lbl{gam1}
\eea
 we have as $\la \to \infty$ that
\bea
\E[Y_{\ell+1}  Y_{m+1}  ] \to 
%\int_{\Gamma} \ka(d\bx) \int_{\Gamma} \ka(d\by)
%f(\bx) f(\by)
% \E[ \xi_\infty^{\bx}(\H_{\ka(\bx)} ,\R^d) ]
% \E[ \xi_\infty^{\by}(\H_{\ka(\by)} ,\R^d) ]
%\nonumber \\
%=
 \gamma_1^2.
\lbl{E1lim}
\eea
Next, observe that
%set
\bea
%E_2 : = \E \left[
%\langle f ,\xi_\la(\bX_{m+1} ; \X_m )\rangle
%\sum_{i=1}^\ell \langle
%f ,\xi_\la(\bX_i;\X_{\ell+1} ) - 
%\xi_\la(\bX_i;\X_\ell)  \rangle \right] 
%\\
 \ell \E[ \Delta_{1,\ell}Y_{m+1}]
=  \ell \E[\langle f, 
 \xi_\la(\bX_1;\X_{\ell+1} ) - \xi_\la(\bX_1; \X_\ell) \rangle   
\langle f, \xi_\la(\bX_{m+1};\X_m)\rangle ]
\nonumber \\
= \ell \int \ka(\bx) d\bx \int \ka(\by) d\by \int \ka(\bw) d\bw 
\E[  (
\langle f, \xi_\la(\bx; \X_{\ell-1}^\bw  )
- \xi_\la(\bx; \X_{\ell-1}  )) \rangle
\nonumber \\
\times
\langle f, \xi_\la(\by;\X_{m-2}\cup\{\bx,\bw\}  )\rangle ].  
\eea
%where $\X'_{m-\ell -1}$ is independent of $\X_{\ell-1}$.

Making the change of variables $\bz =  \la^{1/d}(\bw - \bx)$ we obtain
\bea
 \ell \E[ \Delta_{1,\ell}Y_{m+1}]
 = \frac{\ell}{\la} \int 
 \ka(\bx)(d\bx)
\int
 \ka(\by)(d \by)
 \int
 \ka(\bx + \la^{-1/d} \bz) d\bz
%f_\la(\bx) f_\la(\by)  
~~~~~~~
~~~~~~~
\nonumber \\
\times \E [ 
\langle f, \xi_\la(\bx; \X_{\ell-1}^{\bx + \la^{-1/d} \bz}  )
- \xi_\la(\bx; \X_{\ell-1}  ) \rangle
\times
\langle f, \xi_\la(\by;\X_{m-2}\cup\{\bx,\bx + \la^{-1/d} \bz\}) \rangle  ].  
%~~~~~~~~~
\lbl{0831c}
\eea
Suppose that A1 or A3 holds.
By Lemma \ref{lemweak3}, %\eq{bigweak},
for almost all
 $(\bx, \by) \in  \Gamma_0 \times \Gamma_0$,
we have that
\bea
\langle f, \xi_\la(\bx; \X_{\ell-1}^{\bx + \la^{-1/d} \bz}  )
- \xi_\la(\bx; \X_{\ell-1} )\rangle
\langle f, \xi_\la(\by;\X_{m-2}\cup\{\bx,\bx + \la^{-1/d} \bz\} ) \rangle  
\nonumber \\
\tod f(\bx) f(\by) 
(\xi_\infty^\bx(\H_{\ka(\bx)}^\bz,\R^d ) 
- \xi_\infty^\bx(\H_{\ka(\bx)},\R^d) )
 \xi_\infty^\by(\tH_{\ka(\by)},\R^d). 
%{\bf 1}_\Gamma(\bx) {\bf 1}_\Gamma(\by),
%\nonumber \\
\lbl{0831a}
\eea
By \eq{mom4}
and the Cauchy-Schwarz inequality, 
the variables in the left side of \eq{0831a} are uniformly integrable.
Therefore we have
convergence of expectations, 
so that 
%provided $\bx$ lies at a continuity point of $\ka$, 
the integrand in
\eq{0831c} tends to
\bea
\ka(\bx)^2 \ka(\by) f(\bx) f(\by) 
%\1_\Gamma(\bx) \1_\Gamma(\by)
\nonumber \\
\times
\E [ \xi_\infty^\bx(\H_{\ka(\bx)}^\bz,\R^d ) - \xi_\infty^\bx(\H_{\ka(\bx)},\R^d) ]
\E[ \xi_\infty^\by(\H_{\ka(\by)},\R^d)].
\lbl{0831g}
\eea
Also,
$\langle f, \xi_\la(\bx; \X_{\ell-1}^{\bx + \la^{-1/d} \bz}  ) 
%\rangle
 -
% \langle f, 
 \xi_\la(\bx; \X_{\ell-1}  )  \rangle$ is zero unless
 $|\bz| \leq R_{\la,\ell-1}(\bx) $,
defined at \eq{5311b}.
Also, by A$6'$ or A$7'$, we assume for some $p>2$ and $q>2dp/(p-2)$
that  
the moments condition \eq{mom4}
holds and we have binomial power law stabilization of order $q$
(in A$7'$, since $p>3$ the condition $q > d(150+6/p)$
 ensures that $q>2dp/(p-2)$).
Therefore the H\"older and Minkowski inequalities yield 
\bea
|
\E [\langle f, 
\xi_\la(\bx; \X_{\ell-1}^{\bx + \la^{-1/d} \bz}  )
- \xi_\la(\bx; \X_{\ell-1}  ) \rangle \langle f,
\xi_\la(\by;\X_{m-2}\cup\{\bx,\bx + \la^{-1/d} \bz\}  ) \rangle] |  
\nonumber \\
\leq 
(\|\langle f, \xi_\la(\bx; \X_{\ell-1}^{\bx + \la^{-1/d} \bz}  )\rangle
\|_p 
+
\|\langle f, \xi_\la(\bx; \X_{\ell-1} )\rangle
\|_p) 
\nonumber \\
\times
\| \langle f, \xi_\la(\by;\X_{m-2}\cup\{\bx,\bx + \la^{-1/d} \bz\}  \rangle 
\|_p  
P[ R_{\la,\ell -1}(\bx) > | \bz| ]^{1-2/p}
\nonumber \\
\leq \const ( |\bz|^{q(2-p)/p} \wedge 1).
\lbl{0831h}
\eea
Since  $q > dp/(p-2)$, this is integrable in $\bz$.
Set
\bea
\gamma_2 := 
%\la \E[\Delta_{1,\ell} Y_{m+1}]
% \to
 \int_\Gamma \ka(\bx)^2 d\bx f(\bx) \int_{\R^d} d\bz \E 
[\xi_\infty^\bx(\H_{\ka(\bx)}^\bz,\R^d)  
-\xi_\infty^\bx(\H_{\ka(\bx)},\R^d) ].
\lbl{gam2def}
\eea
By \eq{0831c} and the dominated convergence theorem we obtain
\bea
%E_2
\ell \E[\Delta_{1,\ell} Y_{m+1}]
 \to
% \int_\Gamma \ka(\bx)^2 d\bx f(\bx) \int_{\R^d} d\bz \E 
%[\xi_\infty^\bx(\H_{\ka(\bx)}^\bz,\R^d)  
%-\xi_\infty^\bx(\H_{\ka(\bx)},\R^d) ]
%\nonumber \\ 
%\times
\gamma_2 
\int_\Gamma
\ka(\by) d\by  f(\by) \E
[ \xi_\infty^\by (\H_{\ka(\by)},\R^d)]  = \gamma_1 \gamma_2.
\lbl{E2lim}
\eea
%We define $E_3$ by
%\bea
%E_3 := \E \left[
%\langle f, \xi_\la(\bX_{\ell+1} ; \X_\ell ) \rangle
%\sum_{j=1}^m \langle f, \xi_\la(\bX_j;\X_{m+1} ) - \xi_\la(\bX_j;\X_m) ) 
% \rangle
%\right]   
%\nonumber 
%\\
%= E_{3,1} + E_{3,2} + E_{3,3},
%\lbl{E3eq}
%\eea
%where we define $E_{3,i}, i =1,2,3$ by
%%First, we have
%\bean
%E_{3,1} : = \ell
% \E[ \langle f,
%\xi_\la(\bX_{\ell+1};\X_\ell) \rangle
%\langle f, \xi_\la (\bX_1;\X_{m+1}) - \xi_\la(\bX_1;\X_m)  \rangle ] ; \\
%E_{3,2} := (m-\ell -1) 
% \E[ \langle f, \xi_\la(\bX_{\ell+1};\X_\ell) \rangle
%\langle  f
%, \xi_\la (\bX_{\ell+2};\X_{m+1}) - \xi_\la(\bX_{\ell+2};\X_m)  \rangle ];
% \\
%E_{3,3} := 
% \E[ \langle f_\la,\xi_\la(\bX_{\ell+1};\X_\ell ) \rangle 
%\langle f, \xi_\la (\bX_{\ell+1};\X_{m+1}) 
%- \xi_\la(\bX_{\ell+1};\X_m) \rangle ]. 
%%\\
%%= (3a) + (3b) + (3c).
%\eean
Next, writing $\bX_1$ as $\bx$, $\bX_{\ell+1}$ as $\by$, and
$\bX_{m+1}$ as $\bx + \la^{-1/d} \bz$, we have
\bea
%E_{3,1}
\ell \E[Y_{\ell+1} \Delta_{1,m}] = \frac{\ell}{\la} \int
 \ka(\bx) d\bx
 \int
 \ka(\by) d\by
 \int
 \ka(\bx + \la^{-1/d}
\bz ) d\bz
% f_\la(\bx) f_\la(\by) 
\nonumber \\
\times
\E [ \langle f, \xi_\la(\by; \X_{\ell-1}^\bx) \rangle
 \langle f, \xi_\la (\bx;\X_{m-2}^\by \cup \{\bx + \la^{-1/d}\bz\}  )
 - \xi_\la (\bx;\X_{m-2}^\by    ) \rangle ].
\lbl{0831d}
\eea
By 
%\eq{bigweak},
Lemma \ref{lemweak3} (note that we do not assume that $\ell <m$ in that 
result), for almost all  $(\bx,\by) \in \Gamma_0 \times \Gamma_0$,
we have
\bea
 \langle f, \xi_\la(\by; \X_{\ell-1}^\bx) \rangle 
\langle f,
 \xi_\la (\bx;\X_{m-2}^\by \cup \{\bx + \la^{-1/d}\bz\}  ) 
-
 \xi_\la (\bx;\X_{m-2}^\by   ) 
\rangle
\nonumber \\
\tod f(\bx) f(\by) \xi_\infty^\by (\H_{\ka(\by)},\R^d)
(\xi_\infty^\bx ( \tH_{\ka(\bx)}^\bz,\R^d) 
%\1_{\Gamma}(\bx)\1_{\Gamma}(\by);
- \xi_\infty^\bx ( \tH_{\ka(\bx)},\R^d) ). 
\lbl{0831e}
\eea
Using \eq{mom4}, we obtain convergence of expectations
corresponding to \eq{0831e}.
 Hence, 
%when $\bx$ is
%at a continuity point of $\ka$,
  the integrand in
\eq{0831d} converges to the expression given at \eq{0831g}.
By a similar argument to the one used to establish
\eq{0831h}, the absolute value of this integrand
 is bounded by a constant
times $|\bz|^{q(2-p)/p} \wedge 1$,
 % where $q'$ is the order of assumed binomial power-law 
% stabilization (in Assumption A$7'$ or A$6'$),
 and this is integrable since % we take
 $q > dp/(p-2)$.
Hence, the dominated convergence theorem gives us
\bea
\ell \E[ Y_{\ell+1} \Delta_{1,m} ]
 \to \gamma_1 \gamma_2.
\lbl{E31lim}
\eea

Next, by
taking $\bX_1=\bx$, $\bX_2= \by$, $\bX_{\ell+1} = \bx + \la^{-1/d}\bz$
and $\bX_{m+1}= \by + \la^{-1/d} \bw$, we have
\bea
%E_{4,1} 
%\ell(\ell-1) 
\E[\Delta_{1,\ell} \Delta_{2,m}]
=
% \frac{\ell (\ell -1)}{\la^2}
\la^{-2}
 \int \ka(\bx)d\bx \int
 \ka(\by)d\by
\int  \ka(\bx + \la^{-1/d}\bz)d \bz \int \ka (\by + \la^{-1/d}   \bw)   
d \bw 
\nonumber \\ \times
\E [
\langle f,  \xi_\la(\bx; \X_{\ell -2}^\by \cup \{\bx + \la^{-1/d} \bz \} ) -
 \xi_\la(\bx; \X_{\ell -2 }^\y    ) \rangle
\nonumber \\
\times \langle f, 
 \xi_\la(\by; \X_{m-3}^\bx \cup
\{  \bx+ \la^{-1/d}\bz, \by + \la^{-1/d} \bw \}) 
- \xi_\la(\by; \X_{m-3}^\bx \cup \{\bx+ \la^{-1/d} \bz \}    )
\rangle ].
\nonumber \\
 ~~~~~~~
\lbl{0831m}
\eea
%By \eq{0921a}, \eq{0920g2}, \eq{0921b} and \eq{0920j2} respectively we obtain
For almost all $(\bx , \by) \in \Gamma_0\times \Gamma_0$,
 Lemma \ref{lemweak3} yields
\bean 
 \langle f, \xi_\la(\bx; \X_{\ell -2}^\by \cup \{\bx + \la^{-1/d} \bz \} ) 
 - \xi_\la(\bx; \X_{\ell -2}^\by  ) 
\rangle 
\\
\times 
\langle f,
 \xi_\la(\by; \X_{m-3}^\bx \cup
\{ \bx+ \la^{-1/d}\bz, \by + \la^{-1/d} \bw  \}) 
 - \xi_\la(\by; \X_{m-3}^\bx \cup
\{  \bx+ \la^{-1/d}\bz \}) 
\rangle 
\\
\tod 
f(\bx)
% \1_{\Gamma}(\bx) 
(\xi_\infty^\bx(\H_{\ka(\bx)}^{\bz},\R^d)
 - \xi_\infty^\bx(\H_{\ka(\bx)},\R^d ) ) 
%\\ \times 
f(\by)
%  \1_{\Gamma}(\by) 
(\xi_\infty^\by ({\tH}_{\ka(\by)}^\bw,\R^d) -
 \xi_\infty^\by ({\tH}_{\ka(\by)},\R^d ) ).
\eean
Also, the quantity on the left 
is uniformly integrable by 
the assumption that \eq{mom4}
holds for some $p>2$.
Hence we have corresponding convergence of expectations, so
the integrand
in \eq{0831m} converges to
\bean
f(\bx) f(\by) \ka^2(\bx) \ka^2(\by) 
%\1_{\Gamma}(\bx)\1_{\Gamma}(\by)
\E [\xi_\infty^\bx(\H_{\ka(\bx)}^{\bz},\R^d)
 - \xi_\infty^\bx(\H_{\ka(\bx)},\R^d ) ]
\\
\times \E 
[\xi_\infty^\by ({\H}_{\ka(\by)}^\bw,\R^d) - \xi_\infty^\by ({\H}_{\ka(\by)},\R^d ) ].
%\lbl{0831n}
\eean
Also, 
we have
$\langle f, \xi_\la(\bx; \X_{\ell -2}^\by \cup \{\bx + \la^{-1/d} \bz \} )
\rangle = \langle f,
 \xi_\la(\bx; \X_{\ell -2 }^\y  )\rangle  $ unless 
$R_{\la,\ell-2}(\bx;\{\by\}) > |\bz|$
 and
$ \langle f, \xi_\la(\by; \X_{m-3}^{\bx} \cup \{ 
\bx+ \la^{-1/d}\bz, \by + \la^{-1/d} \bw\}) \rangle
 = \langle f, \xi_\la(\by; \X_{m-3}^{\bx} 
  \cup \{\bx+ \la^{-1/d} \bz \}  )  \rangle$ unless
$R_{\la,m-3}(\by;\{\bx,\bx+ \la^{-1/d} \bz\}) \geq |\bw|$. 
Hence, H\"older's inequality shows
that the absolute value of the expectation in \eq{0831m} is at most
\bean
\| \langle f, \xi_\la(\bx; \X_{\ell -2}^\by \cup \{\bx + \la^{-1/d} \bz \} ) -
 \xi_\la(\bx; \X_{\ell -2 }^\y  ) \rangle  \|_p
\nonumber \\
\times
\|\langle f, \xi_\la(\by; \X_{m-3}^{\bx}  \cup
\{\bx+ \la^{-1/d}\bz,\by + \la^{-1/d} \bw\}) - \xi_\la(\by; \X_{m-3}^{\bx} 
  \cup \{\bx+ \la^{-1/d} \bz \}  )  \rangle \|_p 
\nonumber \\
\times
(P[ R_{\la,\ell-2}(\bx;\{\by\}) \geq |\bz|])^{(1/2)-1/p} 
(P[ R_{\la,m-3}(\by;\{\bx,\bx + \la^{-1/d} \bz\}) \geq |\bw|])^{(1/2)-1/p}. 
\eean 
By the assumption (A$6'$ or A$7'$)
 that  moments condition \eq{mom4} holds for
some $p>2$, and that
 $\xi$ is binomially power law stabilizing
of order $q > 2dp/(p-2)$, this is bounded by
a constant times
$$
(|\bz|^{q(2-p)/(2p)} \wedge 1)
(|\bw|^{q(2-p)/(2p)} \wedge 1)
$$
which is integrable in $(\bz,\bw)$.
% since we assume $q' > 2dp/(p-2)$.
Therefore the dominated convergence theorem applied to \eq{0831m} shows that
\bea
\ell(\ell-1) \E[ \Delta_{1,\ell} \Delta_{2,m} ]
 \to
% \int_\Gamma \int_\Gamma d\by d\bx f(\bx) \ka(\bx)^2 f(\by) \ka(\by)^2 
%\nonumber
%\\
% \times 
%\int_{\R^d} 
%[\xi_\infty^\bz ({\H}_{\ka(\bx)}^\bz,\R^d) - 
%\xi_\infty^\bz ({\H}_{\ka(\bz)},\R^d ) ]
%d \bz 
%\int_{\R^d} 
%[\xi_\infty^\by ({\H}_{\ka(\by)}^\bw,\R^d) -
% \xi_\infty^\by ({\H}_{\ka(\by)},\R^d ) ]
%d\bw
%\nonumber \\
%=
  \gamma_2^2.
%\left(\int_\Gamma d\bx f(\bx) \ka^2(\bx) \int_{\R^d} 
%[\xi_\infty^\bx ({\H}_{\ka(\bx)}^\bz,\R^d)
% - \xi_\infty^\bx ({\H}_{\ka(\bx)},\R^d ) ]
%d \bz \right)^2. 
\lbl{0831r}
\eea

%Turning to $E_{4,4}$, we take
Next, take
 $\bX_1=\bx$, $\bX_{\ell +1}= \by$, $\bX_{m+1}= \bx + \la^{-1/d}\bz$
to obtain
\bean 
%E_{4,4} 
\ell \E[ \Delta_{1,\ell} \Delta_{1,m}]
 = \frac{\ell}{\la} \int d\bx \int d\by  \int d\bz
 \ka(\bx) \ka(\by) \ka(\bx + \la^{-1/d} \bz)
 \nonumber \\
\times \E[ 
\langle f, \xi_\la(\bx;\X_{\ell-1}^{\by})
- \xi_\la(\bx;\X_{\ell-1})  \rangle \langle f,
%(\xi_\la(\by + \la^{-1/d} \bz; \X_{m-1}^\by \cup \{ \bx \} )
\xi_\la(\bx; \X_{m-2}^\by \cup \{ \bx + \la^{-1/d} \bz\} )
%\nonumber \\
- \xi_\la(\bx ; \X_{m-2}^\by  )  \rangle ].
\eean
By Lemma \ref{lemweak3} 
%(NEEDS YET MORE IN LEM \ref{lemweak3})
the quantity inside the expectation
tends to zero in probability for almost all $\bx,\by$ and all $\bz$.
Hence its expectation tends to zero as well, since it
is uniformly integrable by \eq{mom4}.
Also, the absolute value of 
this expectation is bounded by a constant times 
$|\bz|^{q(2-p)/p} \wedge 1$,
by a similar argument to \eq{0831h}.
Hence, dominated convergence yields
\bea
\ell \E[ \Delta_{1,\ell} \Delta_{1,m}] \to 0.
\lbl{e44lim}
\eea
%that given above for $E_{4,3}$. Hence,
%$E_{4,4}$ tends to zero.

Next we have
\bea
%E_{3,3}
\E[Y_{\ell+1}\Delta_{\ell+1,m}] 
= \int \ka(\bx)  \E[
\langle f,
\xi_\la(\bx;\X_{\ell})
\rangle
\langle f,
\xi_\la(\bx;\X_{m}) -\xi_\la(\bx;\X_{m-1})\rangle ].
\lbl{0831k}
\eea
By Lemma \ref{lemweak3},
for almost every $\bx \in \Gamma_0$,
we have
\bean
\langle f, \xi_\la(\bx;\X_{\ell}) \rangle \langle f ,
 \xi_\la(\bx;\X_{m}) \rangle \tod 
f(\bx)^2 (\xi_\infty^\bx( \H_{\ka(\bx)} ,\R^d))^2 
%\1_{\Gamma}(\bx);
\\
\langle f,
\xi_\la(\bx;\X_{\ell}) \rangle \langle f, \xi_\la(\bx;\X_{m-1}) \rangle \tod 
f(\bx)^2 (\xi_\infty^\bx( \H_{\ka(\bx)} ,\R^d))^2, 
%\1_{\Gamma}(\bx),
\eean
and using \eq{mom4}, we have the corresponding convergence of expectations
so that the integrand in \eq{0831k} tends to zero. Also by \eq{mom4},
this integrand is bounded, and thus
% $E_{3,3}$
\bea
\E[Y_{\ell+1}\Delta_{\ell+1,m}] \to 0. 
\lbl{e33lim}
\eea
% tends to zero.

Next, setting $\bX_1 =\by$, $\bX_{\ell+1} =\bx$, and
$\bX_{m+1} = \by +\la^{-1/d} \bz$, we find that
\bea
%E_{4,3}
\E[ \Delta_{1,\ell} \Delta_{\ell+1,m}] 
 = \frac{\ell}{\la} \int d\by \int d\bx  \int d\bz
 \ka(\by) \ka(\bx) \ka(\bx + \la^{-1/d} \bz)
%f_\la(\by + \la^{-1/d} \bx)
 \nonumber \\
\times \E[ 
\langle f, \xi_\la(\by;\X_{\ell-1}^{\bx})
- \xi_\la(\by;\X_{\ell-1})  \rangle \langle f,
\xi_\la(\bx; \X_{m-2}^\by \cup \{ \bx + \la^{-1/d} \bz\} )
\nonumber \\
- \xi_\la(\bx ; \X_{m-2}^\by  )  \rangle ].
\lbl{0831q}
\eea
%Under the assumptions of
By Lemma \ref{lemweak3}, 
%by that result 
%By \eq{bigweak},
%we have 
for almost all $(\bx,\by) \in \Gamma_0 \times \Gamma_0$ and all $\bz$,
 as $\la \to \infty$  we have
\bean
\langle f, \xi_\la(\by;\X_{\ell-1}^{\bx})
- \xi_\la(\by;\X_{\ell-1})  \rangle 
%\times
\langle f,
\xi_\la(\bx; \X_{m-2}^\by \cup \{ \bx + \la^{-1/d} \bz\} )
- \xi_\la(\bx ; \X_{m-2}^\by  )  \rangle 
\toP 0,
\eean
so that
%for each $\bx,\by,\bz$,
the quantity inside the 
 expectation in \eq{0831q} tends to zero in probability
and by \eq{mom4} it is uniformly integrable. Hence
the integrand in \eq{0831q} tends to zero. 
Also, by a similar argument to \eq{0831h},
the absolute value of this integrand is bounded by a constant times
$ |\bz|^{q(2-p)/p} \wedge 1$,
which is integrable since 
 $q > pd/(p-2)$.
Thus,
the integrand in \eq{0831q} is bounded by an integrable function
of $(\bx,\by,\bz)$ so the dominated convergence theorem
shows that %$E_{4,3}$ tends to zero. 
\bea
\ell \E[\Delta_{1,\ell} \Delta_{\ell+1,m} ] \to 0.
\lbl{e43lim}
\eea

Next,  write 
$\bX_{\ell+2}$ as $\bx$, 
$\bX_{\ell+1}$ as $\by$, 
and
$\bX_{m+1}$ as $\bx + \la^{-1/d} \bz$, to obtain
\bea
 %(m-\ell -1)
 \E[Y_{\ell+1} \Delta_{\ell+2,m}] 
 = \la^{-1}
% \frac{(m-\ell -1)}{\la}
\int \ka(\bx) d\bx \int \ka(\by)d\by \int  \ka(\bx+ \la^{-1/d} \bz )
d\bz
\nonumber \\ \times 
 \E[ \langle f, \xi_\la(\by;\X_\ell) \rangle
\langle f,   \xi_\la(\bx;\X_{m-2}^\by \cup \{\bx+\la^{-1/d}\bz \})
- \xi_\la(\bx;\X_{m-2}^\by  )\rangle].
\lbl{0831j}
\eea
By a similar argument to \eq{0831h},  the absolute value of
the expectation inside the
integral is bounded by a constant times $|\bz|^{q(2-p)/p} \wedge 1$,
which is integrable since $q > dp/(p-2)$.
Therefore, the triple integral in \eq{0831j} is bounded, and 
since $m-\ell-1 =o(\la)$, it follows that 
as $\la \to \infty$ we have
%$E_{3,2}$ tends to zero.
\bea
 (m-\ell -1) \E[Y_{\ell+1} \Delta_{\ell+2,m}] \to 0.
\lbl{e32lim}
\eea 

Next,  put $\bX_1 = \bx$, $\bX_{\ell+2}= \by$, 
$\bX_{\ell+1}= \bx + \la^{-1/d} \bz$, and 
$\bX_{m+1}= \by + \la^{-1/d}\bw$, to obtain
\bea
%\ell (m-\ell -1)
 \E[\Delta_{1,\ell}\Delta_{\ell+2,m}]
=
 %\frac{\ell (m-\ell -1 )}{\la^2}
 \la^{-2}
 \int  \ka(\bx)d\bx \int 
 \ka(\by)d\by
\int
d \bz  
 \int
d \bw 
\nonumber \\ \times
 \ka(\bx + \la^{-1/d}\bz) \ka (\by + \la^{-1/d}   \bw) 
\E [ \langle f,
 \xi_\la(\bx; \X_{\ell -1}^{\bx + \la^{-1/d} \bz } ) -
 \xi_\la(\bx; \X_{\ell -1 }  )  \rangle
\nonumber \\
\times \langle f,
 \xi_\la(\by; \X_{m -3}^{\bx} \cup
 \{  \bx + \la^{-1/d} \bz, \by + \la^{-1/d} \bw\} ) -
 \xi_\la(\by; \X_{m -3 }^{\bx} \cup \{  \bx + \la^{-1/d} \bz \} ) 
  \rangle ].
\nonumber \\
\lbl{0831p}
\eea
By the argument used in dealing with $\E [\Delta_{1,\ell} \Delta_{2,m}]$
 above,
the absolute value of the
 integrand in \eq{0831p} is bounded by a constant
times
$(|\bz|^{q(2-p)/(2p)} \wedge 1) (|\bw|^{q(2-p)/(2p)} \wedge 1)$, and
hence the integral in \eq{0831p} is bounded.
Since $\ell(m-\ell -1) = o(\la^2)$, this shows that 
\bea
\ell (m-\ell -1) \E[\Delta_{1,\ell}\Delta_{\ell+2,m}] \to 0.
\lbl{e42lim}
\eea
We have obtained
limiting expressions for the nine terms in
the right hand side of \eq{eqexpand}, namely
\eq{E1lim}, \eq{E2lim}, \eq{E31lim}, \eq{0831r}, \eq{e44lim}, \eq{e33lim},
\eq{e43lim}, \eq{e32lim} and \eq{e42lim}. Combining these, we find that
as $\la \to \infty$,
$\E [ \ch_\ell \ch_{m}]$ converges 
to $(\gamma_1 + \gamma_2)^2$.
%
%
%Thus, as $\la \to \infty$, we have that 
%$$
%E_4 \to 
%\left(\int_\Gamma d\bx f(\bx) \ka(\bx)^2 \int_{\R^d} 
%[\xi_\infty^\bx ({\H}_{\ka(\bx)}^\bz,\R^d) - 
%\xi_\infty^\bx ({\H}_{\ka(\bx)},\R^d ) ]
%d \bz \right)^2. 
%$$
%Combining this with our earlier limits \eq{E1lim}, \eq{E2lim} and \eq{E3lim}
% for the other terms $E_1$, $E_2$ and $E_3$,
%and using used the definition \eq{deltadef} of
%$\delta(\bx, \la)$, 
%we obtain that the expression \eq{0827a} converges to
%\bean
% \left( \int f(\bx) \ka(\bx) d\bx \E\xi_\infty^\bx(\H_{\ka(\bx)},\R^d )
%\right.
%\\
%\left.  + 
% \int_{\Gamma} \ka(\bx)^2 f(\bx) d\bx \int_{\R^d} d\bz
%\E[\xi_\infty^\bx(\H_{\ka(\bx)}^\bz,\R^d)-
% \xi_\infty^\bx(\H_{\ka(\bx)} ,\R^d)]\right)^2
%\\
%= \left( \int_{\Gamma} f(\bx) \ka(\bx) \delta(\bx, \ka(\bx)) d\bx  \right)^2.
%\eean
Since the choice of $\ell(\la),m(\la)$ is arbitrary, we then have
\eq{depoE12}.
$\qed$ 

\begin{lemm}
\lbl{lem0831}
Suppose the assumptions of Theorem \ref{thm3} hold.
Suppose $h(\la)$ is defined for $ \la \geq 1$  satisfying
$h(\la)>0$ and $h(\la)/\la \to 0$ as $\la \to \infty$. Then
$\ch_{m,\la}$ defined at \eq{Rdef} satisfy
\bea 
\lim_{\la \to \infty}  \sup_{\la-h(\la)  \leq m   \leq \la+ h(\la) }
\left( \E [ \ch_{m,\la}^2 ]  \right) < \infty.
\lbl{depoE2}
\eea
%It is the case that 
%\bean
%\sup_m \sup_{\la - \la^{-/3/4} \leq m \leq \la + \la^{3/4}}
%\E \ch_{\la,m}^2 < \infty.
%\eean
\end{lemm}
{\em Proof.}
We abbreviate notation as in the preceding proof.
Note that our assumptions (in particular A$6'$ or A$7'$) imply that
for some $p>2$ and $q> 2dp/(p-2)$,
 \eq{mom4} holds and $\xi$ is binomially power-law stabilizing
of order $q$.

Let $m= m(\la), \la \geq 1$, be defined to satisfy $m(\la) \sim \la$
as $\la \to \infty$. 
By a similar expansion to \eq{eqexpand}, we obtain
\bean
\E[\ch_{m,\la}^2 ] = \E[Y_{m+1}^2] + 2 (m -1) \E[Y_{m+1}\Delta_{1,m}]
+ m(m-1) \E[\Delta_{1,m} \Delta_{2,m}] + m  
 \E[\Delta_{1,m}^2].
\eean
 We consider these terms one by one. 
First,
 $\E[Y_{m+1}^2]$
is  bounded by  \eq{mom4}. Second, 
setting $\bX_1=\bx$ and $\bX_{m+1} = \bx + \la^{-1/d}\bz$ we have
that
\bean
2 m \E[ Y_{m+1} \Delta_{1,m}]
 = \frac{2m}{\la} \int d\bx  \int d\bz \ka(\bx) \ka(\bx + \la^{-1/d} \bz )
% f_\la(\bx) f_\la(\bx
%+ \la^{-1/d} \bz)
 \nonumber \\
\times \E [\langle f, 
\xi_\la(\bx + \la^{-1/d}\bz; \X_{m-1}^\bx)
\rangle \langle f,
\xi_\la(\bx; \X_{m-1}^{\bx + \la^{-1/d} \bz} ) -
\xi_\la(\bx; \X_{m-1} ) \rangle].
\eean
%Here and below, all integrals are over $\R^d$ unless otherwise specified.
Since 
$\langle f, \xi_\la(\bx; \X_{m-1}^{\bx + \la^{-1/d} \bz} )\rangle = \langle f,
\xi_\la(\bx; \X_{m-1} ) \rangle $ when $R_{\la,m-1}(\bx) > |\bz|$, use 
of H\"older's inequality,
followed by \eq{mom4} 
and the 
%assumption that $\xi$ is
 binomial power-law stabilizion,
%of order $q'$ 
shows that the absolute value of the 
expectation in the integrand is bounded by
\bea
\| \langle f, \xi_\la(\bx + \la^{-1/d}\bz; \X_{m-1}^\bx) \rangle \|_p
%\nonumber \\
%\times
\| \langle f,\xi_\la(\bx; \X_{m-1}^{\bx + \la^{-1/d} \bz} ) -
\xi_\la(\bx; \X_{m-1} ) \rangle  \|_p
\nonumber \\
\times
(P[R_{\la,m-1}(\bx)  \geq |\bz| ])^{1 -(2/p)}
\nonumber \\
\leq \const ( |\bz|^{q(2-p)/p} \wedge 1),
~~~~~~~~~~~\lbl{0831s}
\eea
which is an integrable function of $\bz$.
% by the assumption
% that  $q'> dp/(p-2)$.
 This shows that 
$2 m \E[ Y_{m+1} \Delta_m]$ is bounded.

Next, take $\bX_1 = \bx$, $\bX_{m+1} = \bx + \la^{-1/d} \bz$,
and $\bX_2 = \bx + \la^{-1/d} ( \bz +  \bw)$,  
to obtain
\bea
%$E_7
%m(m-1)
\E[\Delta_{1,m}\Delta_{2,m}]
 =
%\frac{ m(m-1)}{\la^2} 
\la^{-2} 
\int \ka(\bx) d \bx \int \ka(\bx + \la^{-1/d}\bx)
d \bz \int \ka(\bx + \la^{-1/d} (\bz + \bw) )   d \bw
%f(\bx) f(\bx+ \la^{-1/d} \bz ) f(\bx + \la^{-1/d} (\bz + \bw)) 
\nonumber \\
\times
\E[\langle f,
\xi_\la(\bx;\X_{m-2}^{\bx + \la^{-1/d}( \bz + \bw)}
 \cup \{ \bx + \la^{-1/d} \bz \} ) - 
\xi_\la(\bx;\X_{m-2}^{\bx + \la^{-1/d}( \bz + \bw)})  \rangle   
\nonumber \\ 
\times \langle f,
\xi_\la(\bx + \la^{-1/d} (\bz + \bw);\X_{m-2}^{\bx }
 \cup \{\bx + \la^{-1/d} \bz \} ) - 
\xi_\la(\bx + \la^{-1/d} (\bz + \bw);\X_{m-2}^{\bx }  )\rangle].
\nonumber \\ 
\lbl{0831t}
\eea
Inside the expectation, 
the first factor is zero if 
$R_{\la,m-2}(\bx;\{\bx+ \la^{-1/d}(\bz +\bw)\}) < |\bz|$, 
and the second factor is zero if
$R_{\la,m-2}(\bx + \la^{-1/d}(\bz + \bw);\{\bx\}) < |\bw|$. Hence 
 by H\"older's inequality, the assumption \eq{mom4},
and the assumption of binomial power-law stabilization of
order $q$, the expectation inside the right side of  \eq{0831t}
is bounded by
\bean
\| \langle f, \xi_\la(\bx;\X_{m-2}^{\bx + \la^{-1/d}( \bz + \bw)}
 \cup \{ \bx + \la^{-1/d} \bz \} ) - 
\xi_\la(\bx;\X_{m-2}^{\bx + \la^{-1/d}( \bz + \bw)})\rangle \|_p   
\nonumber \\
\times \| \langle f, \xi_\la(\bx + \la^{-1/d} (\bz + \bw);\X_{m-2}^{\bx }
 \cup \{\bx + \la^{-1/d} \bz \} ) - 
\xi_\la(\bx + \la^{-1/d} (\bz + \bw);\X_{m-2}^{\bx } ) \rangle \|_p
\nonumber \\
\times 
(P[R_{\la,m-2}(\bx;\{\bx + \la^{-1/d}\bz\}) > |\bz|])^{(1/2)-1/p}
\\
\times
(P[R_{\la,m-2}(\bx + \la^{-1/d} (\bz + \bw);\{\bx\}) > |\bw|])^{(1/2)-1/p}
\nonumber \\
\leq \const
(|\bz| ^{q(2-p)/(2p)} \wedge 1)\times
(|\bw|^{q(2-p)/(2p)} \wedge 1),
%~~~~~~~~
\eean
and since $q> 2dp/(p-2)$, this uniform bound is integrable
in $\bz,\bw$. This 
shows that $m(m-1) \E[\Delta_{1,m}\Delta_{2,m}]$ remains bounded. 

Finally, take $\bX_1=\bx$ and $\bX_{m+1} = \bx + \la^{-1/d}\bz$, to 
obtain
\bea
%E_8 
m \E[\Delta_{1,m}^2]
= \frac{m}{\la} \int d\bx \int d\bz \ka(\bx) \ka(\bx + \la^{-1/d} \bz)
%(f_\la(\bx))^2
 \nonumber \\
\times
\E[ \langle f, \xi_\la( \bx; \X_{m-1}^{\bx + \la^{-1/d} \bz } ) -
 \xi_\la( \bx; \X_{m-1}  )\rangle^2 ].
\eea
Since the quantity inside the expectation is zero unless
$R_{\la,m-1}(\bx) \geq |\bz|$, H\"older's inequality 
shows that this expectation is bounded by
\bean
(\E [ | \langle f,
 \xi_\la( \bx; \X_{m-1}^{\bx + \la^{-1/d} \bz } )
 - \xi_\la( \bx; \X_{m-1}^{\bx  } )
\rangle 
|^{p} ])^{2/p} 
(P[ R_{\la,m-1}(\bx) \geq | \bz|])^{1-2/p}
\nonumber \\
\leq \const ( |\bz|^{q(2-p)/p} \wedge 1),
\eean
which is integrable in $\bz$.
Hence, $m \E[\Delta_{1,m}^2]$ is also bounded. $\qed$ \\

{\em Proof of Theorem \ref{thm3}.} 
Suppose $\|\kappa\|_\infty < \infty$
 and $\kappa$ has bounded support.
Suppose $\xi$ 
is $\ka(\bx)-$homogeneously stabilizing at $\bx$ for $\ka$-almost
all $\bx \in \R^d$, 
satisfies Asuumption A$4$ or A$5$, and also satisfies
A$6'$ or A$7'$.
Let $f \in B(\R^d)$, and assume either A$1$ or A$3$ holds.
Suppose $(\la(n))_{n \geq 1}$ is a $(0,\infty)$-valued sequence
with $|\la(n)-n|= O(n^{1/2})$ as $n \to \infty$.

Let $H_n := \langle f, \nu_{\la(n),n}^\xi \rangle $
 and 
$H'_n := \langle f, \mu_{\la(n)}^\xi \rangle  $.
For this proof, assume that for all $n$,  $\X_n$ is given
by \eq{050404}  and that
$\Po_{\la(n)}$ is coupled to $\X_n$
by setting $\Po_{\la(n)} = \cup_{i=1}^{N_n} \{(\bX_i,T_i)\}$,
%$ \{\bX_{1}, \bX_{2},\ldots,\bX_{N_n} \}$ 
with
%$\{(\bX_{i},T_i)\}$ independent and $\kappa \times$-distributed on $\R^d$,
%and  
$N_n$ an independent Poisson variable with mean $\la(n)$.
Let
$$
\alpha := \int_\Gamma f(\bx) \delta(\bx,\ka(\bx)) \ka(\bx)  
d\bx.
$$ 
%By the proof
%of eqn (4.5) of \cite{PY1},
%%The first step is to prove that
%%as $n \to \infty$,
First we show that as $n \to \infty$,
\begin{equation}
\E \left[ (n^{-1/2} (H'_n - H_n - (N_n-n) \alpha))^2 \right] \to 0.
\lbl{fordepo}
\end{equation}
To prove this, note that the expectation in the left hand  side
is equal to
\bear
n^{-1} \sum_{m: |m-\la(n) |\leq n^{3/4} } \E \left[ 
%n^{-1}
( \langle f, \nu_{\la(n),m}  - \nu_{\la(n),n}  \rangle - (m-n)
 \alpha)^2 \right]
 P[N_n=m]
\nonumber
\\
 + n^{-1} \E \left[  \left(
 H'_n - H_n  - (N_n- n) \alpha \right)^2 {\bf 1}\{ |N_n - \la(n)|
> n^{3/4} \} \right].
\lbl{fordepo1}
\eear
Let $\eps >0$.
By (\ref{Rdef}) and
Lemmas \ref{fordepolem} and \ref{lem0831}, there exists $c>0$ such that for
large enough $n$
and all $m$ with $\la(n) \leq m \leq \la(n)+ n^{3/4}$,
$$
\E[ ( \langle f,\nu_{\la(n),m} - \nu_{\la(n),n} \rangle - (m-n) \alpha )^2]
= \E \left[  ( \sum_{\ell=n}^{m-1} (\ch_{\ell,\la(n)}- \alpha) )^2 \right]
\leq \eps(m-n)^2 + c(m-n),
$$
where the bound comes from expanding out the double sum
 arising from the expectation of the squared sum.
A similar argument applies when $\la(n)-n^{3/4} \leq m \leq n$, and
hence
the first term in (\ref{fordepo1}) is bounded by
the expression
\bean
n^{-1} \E[ \eps(N_n-n)^2 + c |N_n-n| ] 
\\
\leq n^{-1}( \eps(\la(n)-n)^2 + \eps
\E[(N_n-\la(n))^2] + c \E[|N_n-\la(n)|] + c |\la(n)-n| )
\\
\leq n^{-1} ( \eps(\la(n)-n)^2 + \eps \la(n) + c \la(n)^{1/2}
+ c|\la(n)-n|),
\eean
and so, since $\eps$ is arbitrary, the first
 term in \eq{fordepo1} tends to zero.

Our assumptions (A$6'$ or A$7'$) include
the moment bounds \eq{mom} and  \eq{mom4}
 for some $p>2$.
 Hence, choosing $p' \in (2,p)$ 
we can apply Lemma 4.3 of \cite{PY5}, 
 taking $\rho_\la = \la^{1/(2d)}$ in
that result, to bound the $L^{p'}$ norm of the contribution
to $H'$ from points in a cube of side $\la(n)^{-1/(2d)}$ by 
$O(\la(n)^{(p+1)/(2p)})$. The number of such cubes intersecting
$\ska$ is $O(\la(n)^{1/(2d)})$, so that by
Minkowski's inequality we obtain
\bea
\| H'_n \|_{p'} = O( \la^{(p+1)/(2p)} \times \la^{1/2} )
 = O( n^{(2p+1)/(2p)}  ).
\lbl{likepoly}
\eea
Also, the value of $H_n$ is the sum of contributions from 
the $n$ points of $\X_n$, and 
by \eq{mom4}, the $L^p$ norms of each contribution are
bounded, so that by Minkowski's inequality  $\|H_n\|_p = O(n)$
so that $\|H_n\|_{p'} = O(n)$ also. Moreover,
$\|N_n\|_{p'}= O(n)$.
Combining these facts  with \eq{likepoly}, we 
 may deduce that
$$
\|H'_n - H_n - (N_n-n) \alpha \|_{p'} = 
O(n^{(2p+1)/(2p)}).
$$
Hence, by H\"older's
inequality
the  second term in (\ref{fordepo1}) is bounded by
a constant times
$n^{-1}n^{(2p+1)/p}
   (P [|N_n-\la(n)| > n^{3/4}])^{1-(2/p')} $,
which tends to zero (see e.g. Lemma 1.4 of  \cite{Pe}).
 This completes the proof of (\ref{fordepo}).

Set  $\sigma^2:= \int_\Gamma f(\bx)^2 V^\xi(\bx,\ka(\bx))\ka(\bx)d\bx$.
and set $\tau_f^2 := \sigma^2 -\alpha^2$.
By Theorems \ref{thm1} and \ref{thm2}
we have as $n \to \infty$ that
$
%\limn n^{-1}
 \Var(H'_n) 
\to \sigma^2 
$ and
$
n^{-1/2}
(H'_n - \E H'_n )\tod
\NN(0,\sigma^2).
$
Using \eq{fordepo} and following page 1620 of \cite{PY1} verbatim,
we may deduce that 
$
\limn n^{-1} \Var(H_n) 
\to \tau_f^2 
$ and also
\bea
n^{-1/2} \langle f, \overline{\nu}^\xi_{\la(n),n} \rangle 
=
n^{-1/2} (H_n - \E H_n )\tod
\NN(0,\tau_f^2).
\lbl{norm1}
\eea
Since $\tau_f^2$ is
  the limiting variance in \eq{0824b}, we thus have \eq{0824b}.
Suppose that
 $f_1,\ldots,f_k$ are in $\tB(\R^d)$
or that 
 $f_1,\ldots,f_k$ are in $B(\R^d)$
and A1 holds. If
  $a_1,\ldots,a_k$ are
real constants, by \eq{norm1} we obtain convergence of $\sum_{i=1}^k
n^{-1}  a_i \langle f_i, \nu_{\la(n),n} \rangle$ 
in distribution to the centred normal with variance
\bean
\sum_{i=1}^k \sum_{j=1}^k
a_i a_j \left(
 \int_{\Gamma} f_i(\bx) f_j(\bx) V^{\xi}(\bx,\ka(\bx))
\ka(\bx) d\bx \right.
\\
- \left. \int_{\Gamma} f_i(\bx) 
\delta(\bx,\ka(\bx))
\ka(\bx)d\bx
 \int_{\Gamma} f_j(\bx)
 \delta(\bx,\ka(\bx)) 
\ka(\bx)d\bx
\right)
\eean
 with $V^\xi(\bx,a)$ given by (\ref{Vdef}). 
The convergence in distribution follows by the Cram\'er-Wold device.
$\qed$

%We prove convergence of $n^{-1}\Var(H_n)$. This follows from the identity
%$$
%n^{-1/2} H'_n = n^{-1/2} H_n + n^{-1/2} (N_n-n)\alpha
%+ n^{-1/2} (H'_n- H_n - (N_n-n)\alpha).
%$$
%In the right hand side, the third term has variance tending to zero
%by (\ref{fordepo}), while the second term has variance
% $\alpha^2$ and is independent of the first term. It follows that
%with $\sigma^2$ given by Theorem \ref{1Hclt},
%$$
%\sigma^2 =
%\limn n^{-1} \Var(H'_n) = \limn ( n^{-1} \Var(H_n) )+  \alpha^2,
%$$
%so that $\sigma^2 \geq \alpha^2$ and
% $n^{-1} \Var(H_n) \to \tau^2$, where we set $\tau^2 = \sigma^2
% -\alpha^2$.
%This gives us (\ref{varlim}).
%
%Theorem \ref{1Hclt} tells us that
%$
%n^{-1/2}
%(H'_n - \E H'_n )\tod
%\NN(0,\sigma^2).
%$
%Combined with (\ref{fordepo}) this gives us
%$$
%n^{-1/2}(H_n - \E H'_n + (N_n-n) \alpha) \tod
%\NN(0,\sigma^2),
%$$
%and since $n^{-1/2}(N_n-n)\alpha $ is independent of
%$H_n$ and is asymptotically normal with
%mean zero and variance $\alpha^2 $, it follows by considering
%characteristic functions that
%\begin{equation}
%n^{-1/2}(H_n - \E H'_n ) \tod
%\NN(0,\sigma^2-\alpha^2).
%\lbl{fordepo2}
%\end{equation}
%By (\ref{fordepo}), the expectation of
%$n^{-1/2} (H'_n - H_n  - (N_n-n) \alpha)$ tends to zero,
% so in (\ref{fordepo2})
% we can replace $\E H'_n$ by $\E H_n$, which gives
%us (\ref{2clteq}).  Save for showing that
%$\tau^2 $ is strictly positive, this completes the
%proof of Theorem \ref{2Hclt}.
%

\section{Applications}
\lbl{secexamples}

Many examples and applications of the general theory are described in
\cite{BY,PY1,PY2,PY4,PY5}, and here we discuss only
a few. The examples we consider have translation-invariant
$\xi$. There are interesting potential applications
of the theory with non translation-invariant $\xi$ 
to topics in  multivariate statistics
such as nonparametric density estimation and nonparametric
regression \cite{BB,BY2,EJ}, but these are not
easy to describe briefly in an already lengthy paper.
 The first example  discussed here
illustrates the application of the Law of Large Numbers
in Theorem \ref{thlln1}.

\subsection{Voronoi coverage}
\lbl{secvoro}
\allco

For finite $\X \subset \R^d$ and $\bx \in \X$, let
 $C(\bx;\X)$ denote the Voronoi cell 
with nucleus $\bx$ for the Voronoi tessellation induced by $\X$.
Given a density function $\ka$ on $\R^d$, and $\la >0$, let
$\Po_\la$ denote a Poisson point process in $\R^d$ with intensity
function $\la \kappa(\cdot)$.
With a view to potential applications 
in nonparametric statistics and in image analysis,
Khmaladze and Toronjadze \cite{KT} ask whether,
 for an arbitrary bounded Borel set $A \subset \R^d$,  
the total volume of bounded cells $C(\bx;\Po_\la)$
 with nuclei at points $\bx \in \Po_\la \cap A$
converges almost surely  to the Lebesgue measure of $A$, as $\la \to \infty$.
They also ask whether
$$
\left| A\triangle \cup_{\bx \in \Po_\la \cap A} C(\bx;\Po_\la) \right| \to 0,
$$
where $|\cdot |$ here denotes Lebesgue measure. They answer these
questions affirmatively for the case
 $d=1$ only.  
If one is satisfied with $L^1$ convergence, one can use
our law of large numbers (Theorem \ref{thlln1})  to answer the first question
 affirmatively for general
$d$,
and to partially answer the second question also. 
We also have corresponding results for $\X_n$.
It does not seem possible to achieve these results using
only the results of \cite{PY4} or \cite{BY}.

% consider the Voronoi tessellation of $\R^d$ induced
%by $\Po_\la$.

%Actually, the question mentioned above is raised in the abstract
%of \cite{KT}, but in body of \cite{KT}
%it is asked if, 
%Setting $C(\bx;\Po_\la)$ to be the Voronoi cell 
%with nucleus $\bx$ for the Voronoi tessellation induced by $\X$, we have

To put these questions in our framework,
define the set
 $V(\bx;\X)$ to be the set
 $ C(\bx;\X) $ if this set  is bounded, and to be the empty set otherwise.
For finite $\X \subset \R^d$ and $\bx \in \X$,
let $\xi(\bx;\X,\cdot)$ be the restriction of Lebesgue measure
to $V(\bx;\X)$, and let $\xi^*(\bx;\X,\cdot)$
be the corresponding point measure, defined at \eq{eq2.2a}
Note that in this case, $\xi$ is translation-invariant and
points do not carry marks.

Our choice of $\xi$  
has the homogeneity property of order $d$, which says
that $\xi(a \bx;a \X,\R^d) = a^d  \xi(\bx;\X,\R^d)$ for any $a >0$.
Hence, for Borel $A \subseteq \R^d$, using \eq{TIxila}
 we have
\bean
\la^{-1} \mu_\la^{\xi^*} (A) = \la^{-1} 
\sum_{\bx \in A \cap \Po_\la} \xi_{\la}(\bx;\Po_\la,\R^d)
%\\
= \la^{-1} \sum_{\bx \in A \cap \Po_\la} \xi( \la^{1/d} \bx;\la^{1/d} \Po_\la,\R^d)
\\
=  \sum_{\bx \in A \cap \Po_\la} \xi(  \bx; \Po_\la,\R^d)
=  \sum_{\bx \in A \cap \Po_\la} |V(\bx;\Po_\la)|.
%\lbl{kh1}
\eean
Similarly,
\bean
\la^{-1} \nu^{\xi^*}_{\la,n} (A) = \sum_{\bx \in A \cap \X_n} \xi(  \bx; \X_n).
%\lbl{kh2}
\eean
By arguments in \cite{PY1}, the measure $\xi$ satisfies
$R(\0;\H_\la) <\infty$ almost surely for all $\la >0$.
Also, 
$$
\E [ \xi_\infty^\bx ( \H_\la,\R^d) ] = 1/\la
$$
since the average volume of Voronoi cells in a homogeneous Poisson process
of intensity $\la$ must be $1/\la$ (one can use Theorem \ref{thlln1}
to show this rigorously).
Moreover, by arguments in \cite{PY1} the measure $\xi$
 satisfies the moments conditions \eq{mom} and \eq{mom2},
for example if
$\ka$ is supported by the unit cube in $\R^d$ and bounded away
from zero on its support.  Hence, under these
conditions on $\ka$, by setting $f$ to be the indicator function of
$\1_A$,
% and writing $|\cdot|$ for Lebesgue measure,
 we can apply Theorem \ref{thlln1} to deduce that  
$$
\la^{-1} \mu_{\la}^{\xi^*}(A)
 \inLL \int_{A \cap \supp(\ka)} (1/\ka(\bx) ) \ka(\bx)d\bx = |A \cap \ska|
$$
and likewise if $\la(n) \sim n$ as $n \to \infty$, then 
$n^{-1} \nu^{\xi^*}_{\la(n),n}(A) \inLL |A \cap \ska|$ for
arbitrary Borel $A \subseteq \R^d$.
This answers the first question raised above, in the $L^1$ sense.

If  $C(\bx,\X)$
is bounded, then  with $f=\1_A$ we have
\bean
\langle f,\xi_\la (\bx;\X)  -  \xi_\la^*(\bx;\X)\rangle 
= \la( |C(\bx,\X)\cap A| - |C(\bx,\X) | \1_A(\bx) ) 
\\
= \la (|C(\bx,\X)\cap A|\1_{A^c}(\bx)  - |C(\bx,\X) \setminus A|
\1_{A}(\bx)) 
\eean
so that if there are no unbounded Voronoi cells  having
non-empty intersection with $A$, then
\bean
 |A \triangle \cup_{\bx \in \Po_\la \cap A} C(\bx; \Po_\la)|
=
\la^{-1} \sum_{\bx \in \Po_\la}
|\langle f,\xi_\la (\bx;\Po_\la)- \xi^*(\bx;\Po_\la) \rangle| 
\eean
%\inL 0
which converges in $L^1$ to zero by Theorem \ref{thlln1}.
In cases where the closure of $A$ is contained in the interior of
$\ska$, the probability of there being any unbounded Voronoi cells
intersecting $A$ tends to zero, so the above partially
 answers the second question raised above in an $L^1$ sense.
However, we do not
 fully deal here with sets touching the boundary of $\ska$.

\subsection{Germ-grain models}
\lbl{secgerm}
Germ-grain models are a fundamental model of random
sets in stochastic geometry; see for example \cite{Hall,MR,SKM}.
In the germ-grain model, a random subset of $\R^d$ is
generated as the union of sets $(X_i + T_i)$ where $\{X_i\}$ (the
germs) are the points
of a point process, and $\{T_i\}$ (the grains)
are independent identically distributed
random compact subsets of $\R^d$.  Our results
can be applied to obtain limit theorems for random measures
associated with germ-grain models, in the case where the
point process of germs  is  $\Po_\la$ or $\X_n$, and where
the grains are scaled by a factor of $\la^{-1/d}$ as $\la \to \infty$.
%For ease of presentation, we assume the grains are balls
%centred at  the origin, but it should be clear that our arguments
%carry over to a much wider class of distributions on the grains.

%Set $\MM = [0,\infty)$, and for finite $\X \subset \R^d \times \MM$,  set
Let $\MM$ denote the space of compact subsets of $\R^d$.
For finite $\X \subset \R^d \times \MM$, and $\la >0$,  set
%$\X^\la := $ and if  
%for a set $\X = \{(\bx_i,t_i)\}\subset \R^d \times \MM$  we set
%$\X^\la := \{(\bx_i, \la^{-1/d} t_i) \}$,
$\X^\la := \{(\bx, \la^{-1/d} t): (\bx,t) \in \X \}$,
and set
$$
%\Xi_\la(\X) = \cup_{(\bx,t) \in \X} B_{\la^{-1/d} t}(\bx).   
\Xi_\la(\X) = \cup_{(\bx,t) \in \X} ( \bx +\la^{-1/d} t).   
$$
When $\X$ is $\Po_\la$ or $\X_n$, the set $\Xi_\la(\X)$ is a germ-grain
model with germs given by a Poisson process or binomial process
and grains scaled by a factor of $\la^{-1/d}$.
We can  apply our general results to the volume measure of
 $\Xi_\la(\Po_\la)$
(i.e., the restriction of Lebesgue measure to $\Xi_\la(\Po_\la)$)
%after re-scaling,
% converges
%to a white noise, 
and  the surface measure of $\Xi_\la(\Po_\la)$
(i.e., the restriction of $(d-1)$-dimensional 
Hausdorff measure to the boundary of $\Xi_\la(\Po_\la)$), and likewise 
for $\X_n$.

For $t \in \MM$, i.e. for $t$ a compact set in $\R^d$, let $|t|:= \max\{|\bz|:
\bz \in t\}$ and let $\|t\|$ denote the Lebesgue measure of $t$.

\begin{theo}
\lbl{thvollln}
Suppose that for $q=1$ or $q=2$, and for some $p>q$, we have
%$\ka is bounded$ 
$ E[\|T\|^{p}]< \infty$. Then for $f \in B(\R^d)$ the integral
$\int_{\Xi_\la(\Po_\la)} f(\bx) d\bx $
 converges in $L^q$ to a finite non-random limit. 
If $\la(n) \sim n$ as $n \to \infty$, the integral
$\int_{\Xi_{\la(n)}(\X_n)} f(\bx) d\bx $
converges in $L^q$ to the same limit.
\end{theo}  

We sketch the proof.
For finite $\X \subset \R^d \times \MM$,
%let us view $\X$ as a subset of points $\R^d$
let $\pi(\X)$ denote the projection of $\X$ onto $\R^d$, i.e. the
subset of $\R^d$ obtained if we ignore the marks carried by points of $\X$.
Also, for each $\bx \in \pi(\X)$ let $T(\bx)$ denote the mark carried by
$\x$, i.e. the value of $t$ such that $(\bx,t) \in \X$.
%For $\by \in \Xi_1(\X)  = \cup_{(\bx,t) \in \X}B_t(\bx)$, 
For $\by \in \Xi_1(\X)  = \cup_{(\bx,t) \in \X}(\bx + t)$, 
let $\NN_\X(\by)$ denote 
the nearest point $\bx \in \pi(\X)$   to $\y$ such
that $ \by \in \bx + T(\bx)$
 (in the event of a tie when seeking the `nearest
point', use the lexicographic ordering as a tie-breaker).  
Take $\xi(\bx,t;\X,\cdot)$ to be the restriction of Lebesgue measure
%to the set of  $\by \in B_{T(\bx)}(\bx)$ such that 
to the set of  $\by \in \bx+ t$ such that 
$\bx = \NN_\X(\by)$.

Then, since $\NN_\X(\by)$ is unique 
 for each $\by \in \Xi_1(\X)$, 
$\sum_{(\bx,t) \in \X} \xi(\bx,t;\X,\cdot) $ is precisely
the volume measure of $\Xi_1(\X)$.
Also, $\xi$ is translation-invariant, so
by \eq{TIxila},
\bean
\xi_\la(\bx,t;\X,A) = \xi( \la^{1/d} \bx,t;\la^{1/d}\X,\la^{1/d} A )
= \la \xi(  \bx, \la^{-1/d} t ; \X^\la, A),
\eean
so that
$\la^{-1} \sum_{(\bx,t) \in \X} \xi_\la(\bx,t;\X,\cdot) $ is 
the volume measure of  $\Xi_1(\X^\la)$, which is the same
as the volume measure of $\Xi_\la(\X)$. Hence, with this choice of
$\xi$, we have that
\bea
\la^{-1}\mu_\la^\xi (d\bx) =  \1_{\Xi_\la(\Po_\la)} (\bx) d\bx ; 
\lbl{0415a}
\\
\la^{-1}\nu_{\la,n}^\xi (d\bx) =  \1_{\Xi_\la(\X_n)} (\bx) d\bx . 
\lbl{0415b}
\eea
The measure  $\xi(\bx,t;\X,\cdot)$  is  supported
by $\bx +t$, 
%by $B_t(\bx)$ 
and this measure is unaffected by changes to
$\X$ outside $B_{2|t|}(\bx) \times \MM$. This is because for any
 $\by \in \bx + t$
it is the case that $|\by - \NN_\X(\by)| \leq |t|$ so by the
triangle inequality,
 $ \NN_{\X}(\by)$ lies in $B_{2|t|}(\bx)$.  
Hence, $2|t|$ serves as a radius of stabilization, which
is almost surely finite since $\MM$ is the space of compact
subsets of $\R^d$. Also, $\xi(\bx,t;\X,\R^d)$ is bounded by
$\|t\|$, and the conditions (\ref{mom}) and \eq{mom2} follow. 

We can then apply Theorem \ref{thlln1}   to obtain the result.
The measure $\xi(\bx,t;\X,\cdot)$ satisfies Assumption A2,
so we can take test functions in $B(\R^d)$. The limit is
given by the right hand side of \eq{LLNpo}.
%$\int_{\R^d} (1-\exp(-\ka(\bx)\E\|T\|))\ka(\bx) d\bx$.

By applying Theorems \ref{thm2} and \ref{thm3} to the above
choice of $\xi$, we obtain the following 
 result on convergence
to a Gaussian field  for the volume measure on $\Xi_\la(\Po_\la)$ or 
on $\Xi_{\la(n)}(\X_n)$. 
We state the result in terms of the measures $\mu_\la^\xi$
and $\nu_{\la,n}^\xi$, which translate into statements
about the volume measures by \eq{0415a} and \eq{0415b}.

\begin{theo} 
\lbl{thvolclt}
Suppose $\ka $ is bounded and has bounded support. 
Suppose
% either that 
%$\E[\|T\|^p]<\infty $
%for some $p >2$ 
% and there exists $C>0$ such that
%$ P[|T|>s] < C \exp(-s/C)$ for all $s$, or
for some $p>3 $
that $\E[\|T\|^p] <\infty$
 and there exists $C>0$  and $q >d(150+6/p)$ such
that $P[|T|>s ]< C s^{-q} $ for all $s$.

Then with $\xi$ as given above,
%if we set $\rho_\la$ to be  the volume measure on $\Xi_\la(\Po_\la)$
%and set  $\overline{\rho}_\la$ to be $\rho_\la 
%- \E  \rho_\la$ then  
the finite dimensional distributions of the random
field  $\la^{-1/2} \langle f, \overline{\mu}_\la^\xi \rangle$, $
f \in \tB(\R^d)$,
converge  to those of a centred Gaussian random field
with covariances given by 
$\int_{\R^d} f_1(\bx) f_2(\bx) V^{\xi}(\bx, \ka(\bx))
d\bx$, with $V^{\xi}$ given by \eq{Vdef} and $\xi$ as defined above.
Likewise, if $|\la(n)-n|= O(n^{-1/2})$, then 
the finite dimensional distributions of the random
field  $\la^{-1/2} \langle f, \overline{\nu}_{\la,n}^\xi \rangle$, $
f \in \tB(\R^d)$,
converge  to those of a centred Gaussian random field
with covariances given by 
the right hand side of \eq{0209}.
\end{theo}

%Hence, if $\mu_\MM$ has an exponentially decaying tail, our $\xi$
%will satisfy exponential stabilization. If
% $\mu_\MM$ has a polynomially decaying tail, our $\xi$
%will satisfy polynomial  stabilization of the same order.
%
Theorem \ref{thvolclt} 
%demonstrates
%Using all this, we can apply our general results to obtain
adds
  to the results for germ-grain models
in (\cite{BY}, Section 3.3) in several ways. In particular, in \cite{BY}
 it is assumed that
the distribution of $|T|$  is supported
by a compact interval, whereas here we need only power-law decay 
%some bound on
of the tail of this distribution. Also,
 in $\cite{BY}$ 
the term  `volume measure' is used in a non-standard way to
refer to an atomic measure supported
 by the points of $\X$.
  Our usage of the terminology `volume
measure' seems more natural, and is also in agreement
with the standard usage found, for example, in
\cite{HM,SKM}.
It is not clear whether the general results in \cite{HM}
 can be applied to the volume measure. \\

%Third, as mentioned earlier, our central limit theorem
%for the random fields acts
%on a larger class of test functions than those considered
%in \cite{BY} where only  continuous test functions are considered.

We now consider the surface measure of $\Xi_\la(\X)$. 
We assume here that with probability 1, 
each grain is a finite  union of bounded convex sets.
%
%As in the volume measure case, we would like to express the 
%surface measure as a sum of contributions from germs, 
%with the contribution
%from each germ determined by the configuration of germs
% within a constant multiple of the associated grain radius from that germ.
%With some effort, this can be achieved as follows.
For $(\bx,t) \in \X$, 
%let $B^o_t(\bx)$ and $\partial B_t(\bx)$
let $(\bx+t)^o$ and $\partial (\bx + t)$
denote the interior and boundary, respectively, of the 
set $\bx +t$.
 %$B_t(\bx)$.
For $\bz \in \cup_{(\bx,t) \in \X}(\bx+t)^o$, let $\NN^*_\X(\bz)$
denote the closest point $\bx \in \pi(\X)$ to $\bz$ such that
$\bz \in (\bx +T(\bx))^o$, using lexicographic ordering as
a tie-breaker.

Define the set $NC(\bx,t)$ 
(i.e., the
points for which $\bx$ is the `nearest covering' germ) by
$$
%NC(\bx,t) = \{ \bz \in  B^o_t(\bx): 
%\nexists (\by,u) \in \X \mbox{ with } |\bz -y | < \min (t,u) \}
NC(\bx,t) =
%    B^o_t(\bx) \setminus \cup_{ (\by,u) \in \X} B_{\min (t,u)}(\by) 
 \{ \bz \in (\bx+t)^o: \bx = \NN^*_\X(\bz) \}
$$
%which is the set of points interior to $B_t(\bx)$
which is the set of points interior to $\bx+ t$
which are not covered by a set with germ closer than $\bx$. 
Define the set $CC(\bx,t)$ (the points
of  $\partial(\bx + t)$ with `closer cover') by
$$
%CC(\bx,t) = \{ \bz \in \partial B_t(\bx): 
%\exists (\by,u) \in \X \mbox{ with } |\bz -y | < \min (t,u) \}
CC(\bx,t) = 
 %\partial B_t(\bx) \cap 
%\left( \cup_{ (\by,u) \in \X } B_{ \min (t,u) }(\by) \right) .
\cup_{(\by,u) \in \X\setminus \{(\bx,t)\}} \{\bz \in \partial(\bx+t)
\cap (\by + u)^o
: \by = \NN_\X(\bz) \}.
$$
Let us take $\xi(\bx,t;\X,\cdot)$ to be the following
{\em signed}  measure:
\begin{itemize}
\item
 Let $\xi^+(\bx,t;\X,\cdot)$ be  
the  restriction
 of $(d-1)$-dimensional Hausdorff measure
to $\partial(\bx+ t) \setminus CC(\bx,t)$.
%to the those $\bz$ on the boundary $\partial B_{t}(\bx)$ of $B_t(\bx)$ 
%which are not
% covered by any ball $B_{u}(\by)$ 
%with $(\by,u) \in \X$ and $\y$ closer than
%$\bx$ to $\bz$.
\item
 Let $\xi^-(\bx,t;\X,\cdot)$ be the  restriction
 of $(d-1)$-dimensional Hausdorff measure to the set
$$
NC(\bx,t)
\cap
\left( \cup_{(\bz,u) \in \X \setminus \{(\bx,t)\}} 
 \partial (\bz+ u) \setminus CC(\bz,u) \right)
.
$$ 
%which is the set 
%the set of points $\bw$ such that  
%for some $(\bw,v)$ we have (i) $\bz \in \partial B_v(\bw)$;
%(ii) $\bz$ is not  
% covered by any ball $B_{u}(\by)$ 
%with $(\by,u) \in \X$ and $|\y-\bz| <v$ 
%(iii) $\bz$ IS covered by $B_t(\bx)$, and 
%(iv) $\bz$ is NOT covered by any $B_s(\bx')$ with
%$v < |\bx' -\bz| < \bx - \bz|$.
%Finally
\item Let $\xi (\bx,t;\X,\cdot)$ be the signed measure $\xi^+(\bx,t;\X,\cdot) - \xi^-(\bx,t;\X,\cdot)$.
\end{itemize}

The signed measure 
$\xi(\bx,t;\X,\cdot)$ is supported by $\bx+t$ and is unaffected
by changes to $\X$ outside 
$B_{2|t|}(\bx) \times \MM$. Achieving this is the purpose of the
 definition of $\xi$ used here, since it ensures $|T|$ serves as
a radius of stabilization.

We assert that $\sum_{(\bx,t) \in \X} \xi(\bx,t;\X,\cdot)$ is precisely
the surface measure of $\Xi_1(\X)$.  
To see this, suppose $\bz$ lies on the surface
of $\bx+t$, but is covered by some other $(\by+u)^o$ with
$(\by,u) \in \X$ (take the closest such $\by$ to $\bz$).
If $|\by - \bz| <|\bx - \bz|$, 
then $\bz \in CC(\bx,t)$ so that
$\xi^+(\bx,t;\X,d\bz) = \xi^-(\bx,t;\X,d\bz) = 0$.
 If $|\by -\bz | > |\bx - \bz|$, then  $\bz \notin CC(\bx,t)$ so that
$\bz \in NC(\by,u) \cap \partial (\bx+t) \setminus CC(\bx,t)$, so that
  $\xi^+(\bx,t;\X,d\bz)$ is 
the surface measure of $\partial (\bx+t)$, and $\xi^-(\by,u;\X,d\bz)$
is also the surface measure of $\partial (\bx+t)$, and these cancel out.
If also $\bz \in (\bw+v)^o$ with
 $(\bw,v) \in \X$ and $|\bw -\bz|> |\by-\bz|$,
then $\bz \notin NC(\bw,v)$ so that $\xi^-(\bw,v;\X,d\bz)=0$.

Since $\xi$ is translation-invariant, 
by \eq{TIxila} we have
\bean
\xi_\la(\bx,t;\X,A) = \xi( \la^{1/d} \bx,t;\la^{1/d}\X,\la^{1/d} A )
= \la^{(d-1)/d} \xi(  \bx, \la^{-1/d} t ; \X^\la, A)
\eean
and hence, 
$\la^{(1-d)/d} \mu_\la^\xi$ is the surface measure of $\Xi_\la(\Po_\la)$,
while
$\la^{(1-d)/d} \nu_{\la,n}^\xi$ is the surface measure of $\Xi_\la(\X_n)$.

To apply our general results here, we need the moments conditions
such as \eq{mom} to apply to both the positive and negative
parts of the measure $\xi$. We write $|\partial T|$ for the
$(d-1)$-dimensional Hausdorff measure of the boundary of $T$. Clearly  this is
an upper bound for  
$\xi^+(\bx,T;\X,\R^d)$. To estimate the negative part,
observe that all contributions
to  
$\xi^-(\bx,T;\X)$ come from the boundaries of sets associated with germs
within distance at most $2|T|$ from $\bx$. Hence,
% conditionally on $T$,
we have that
\bean
\E [ \xi^-(\bx,T;\Po_\la,\R^d )^p |T] \leq \E \left[  \left(\sum_{i=1}^N X_i 
\right)^p 
\right]
\eean
where  $X_i$ are independent copies of $|\partial T|$, and $N$ is
Poisson with mean $\theta (2|T|)^d$, with $\theta$
denoting the volume of the unit ball. 
By Minkowski's inequality, the above expectation is bounded by
$\E[N^p] \E[|\partial T|^p]$, and 
hence the condition
\bea
\E[|T|^{dp}] \E[|\partial T|^p] <\infty
\lbl{momsurf}
\eea
suffices to give us all the moments conditions \eq{mom}, \eq{mom2},
\eq{momlnb} and
\eq{mom4}, for both the positive and the negative parts of $\xi$.

Thus, with this choice of $\xi$ we may apply Theorem
\ref{thlln1}
 (with test functions $f \in \tB(\R^d)$)
 if for $q=1$ or $q=2$ we have
\eq{momsurf} for some $p>q$. 
We may apply Theorems \ref{thm2} and \ref{thm3} 
 (again with test functions $f \in \tB(\R^d)$), either
if \eq{momsurf} holds for some $p>2$ and 
 $P[|T|>r]\leq Ce^{-r/C}$ for some $C >0$ and all $r>0$, or
if for some $p>3$ and some $q >d(150+6/p)$, \eq{momsurf} holds and 
$P[|T|>r] \leq C r^{-q}$ for some $C>0$ and all $r>0$.

\subsection{Random  packing measures}
The random packing measures discussed in Section 3.2
of 
 \cite{BY} are obtained by particles (typically balls) 
being deposited
in space at random times, according to a space-time
Poisson process. Particles have non-zero volume and
(in some versions of the model)
may grow with time, but deposition and 
growth are limited by an excluded volume effect. 
In \cite{BY}, the measures associated with 
these packing processes are obtained
as a sum of unit point masses, with one point for each particle.
As in the case of germ-grain models,
it is quite natural instead to consider the
 the volume  measure associated with the random set obtained
as  the union of particles (balls), or even the surface measure
of this random set.  The setup of this paper enables us to do this,
but we do not give details.

%Mathew D. Penrose, 
Department of Mathematical Sciences, 
University of Bath, Bath BA2 7AY, United Kingdom.
\\
Email: {\texttt mathew.penrose@durham.ac.uk} \\
URL: \texttt{http://www.maths.bath.ac.uk/}$\sim$\texttt{masmdp}
%\Comment{do we need to give our addresses twice??}

%\vskip 1pc

%J.E. Yukich, Department of Mathematics, Lehigh University,
%Bethlehem PA, USA 18015:
%\\
%{\texttt joseph.yukich@lehigh.edu}

\end{document}